\documentclass[final]{siamart171218}

\pdfoutput=1

\usepackage[utf8]{inputenc}
\usepackage{amssymb}
\usepackage{amsmath}
\usepackage{algorithm}
\usepackage{algorithmic}
\usepackage{bm}

\usepackage{multirow}
\usepackage{color}
\newtheorem{remark}{Remark}
\newtheorem{property}{Property}

\usepackage{geometry}
\usepackage{subcaption}

\usepackage{mathtools}

\definecolor{myRed}{rgb}{1.0, 0.0, 0.0}

\definecolor{myBlue}{rgb}{0.2, 0.2, 1.0}

\definecolor{myPurp}{rgb}{0.5, 0.0, 1.0}

\renewcommand{\chi}{\mathcal{X}}
\newcommand{\He}{\mathcal{H}}
\newcommand{\gF}{{\bm F}}
\newcommand{\gC}{{\bm C}}
\newcommand{\bigO}{\mathcal{O}}
\newcommand{\incr}[1]{ {\delta_#1 } }

\title{The Characteristic Mapping Method for the Linear Advection of Arbitrary Sets}

\author{Olivier Mercier\thanks{The Department of Mathematics and Statistics, McGill University, Montreal, Canada  H3A 0B9.}
\and Xi-Yuan Yin\footnotemark[1]
\and Jean-Christophe Nave\footnotemark[1]  \thanks{Corresponding author. \newline \hspace*{13pt} E-mail addresses: \email{oli.mercier@gmail.com} (OM). \email{xi.yin@mail.mcgill.ca} (XYY). \email{jcnave@math.mcgill.ca} (JCN).}}

\begin{document}

\maketitle

\begin{abstract}
%% Text of abstract
% !TEX root = CMmethod.tex
%\begin{abstract}
We present a new numerical method for transporting arbitrary sets in a velocity field. The method computes a deformation mapping of the domain and advects particular sets by function composition with the map. This also allows for the transport of multiple sets at low computational cost. Our strategy is to separate the computation of short time advection from the storage and representation of long time deformation maps, employing appropriate grid resolution for each of these two parts. We show through numerical experiments that the resulting algorithm is accurate and exhibits significant reductions in computational time over other methods. Results are presented in two and three dimensions, and accuracy and efficiency are studied.
%\end{abstract}

\end{abstract}


%% main text
% !TEX root = CMmethod.tex
\section{Introduction} \label{sec:introduction}
The transport of objects, densities, or data under a given velocity field is ubiquitous in computational mathematics. Many problems ranging from computational physics to computer graphics require solving the linear advection equation. For example, tracking the evolution of a temperature distribution usually requires solving a transport equation. Additionally, the linear advection equation has been extensively used to evolve closed curves and surfaces using the level-set framework \cite{osher2003level}. However, the evolution of more complicated objects such as open curves and open surfaces \cite{leung2009grid}, or even more general sets (e.g. sets with poor regularity or even fractal structures) remains a challenge.

There are two classes of methods for the advection of surfaces: implicit methods represent the surface as the zero-level-set of a function, while explicit methods represent the surface as a collection of sample points or as a moving mesh.

Explicit methods consist in discretizing advected quantities into Lagrangian particles whose locations follow the velocity field. The computation of particle trajectories is fast, accurate, and highly parallelizable, but these particles are not guaranteed to be well-distributed in the domain at all times. Fully-Lagrangian mesh-free approaches, such as the smoothed-particle hydrodynamics method \cite{gingold1977smoothed}, can deal with complex moving boundaries. However, they require complicated projections and mesh-free interpolation routines onto a fixed spacial grid, for instance when visualizing the fluid, or computing a Laplacian operator for viscous fluids. More recently, Bowman et al. \cite{bowman2015fully} proposed a Fully-Lagrangian mesh-free method with a conservative particle-to-Eulerian-grid projection step. This method allows for the conservation of all Casimir invariants and eliminates numerical dissipation. Other explicit methods, for instance \cite{misztal2012multiphase, clausen2013simulating}, use a moving Lagrangian mesh. In multiphase fluid simulations, these methods naturally preserve fluid volumes. Moreover, fluid boundaries and interfaces have a direct meshed representation, and complex solid boundaries are more easily resolved. However, these meshes can become distorted and complex, and costly remeshing routines are needed to maintain spacial resolution and accuracy.

Implicit methods take an Eulerian approach where the locations of discretized elements are stationary. In particular, level-set methods represent advected surfaces as the zero contour of a level-set function. These are most naturally used when it is important to query whether a point is inside or outside the surface, as indicated by the sign of the level-set function. A summary of standard level-set methods can be found in \cite{osher2003level,osher2001level,sethian1999level,gibou2018review}. Additionally, particle tracking has been used by Enright et al. \cite{enright2002hybrid} to improve accuracy of the standard level-set approach when subgrid structures must be resolved. More recently, Gradient-Augmented level-set (GALS) and Jet Schemes (JS) methods were proposed in \cite{nave2010gradient,chidyagwai2011comparative,seibold2011jet}. These methods couple time integration of velocities with spacial interpolation of advected quantities to maintain a functional representation of the solution at all times. This improves subgrid accuracy while retaining stencil compactness and overall computational efficiency. Further applications of these methods can be found in \cite{lee2014narrow,bockmann2014gradient,kolahdouz2013semi}.

The approach presented in this paper adopts an implicit point of view. We propose an Eulerian method to advect arbitrary sets under a given velocity field. Note that for divergence-free velocity fields, topological changes in the solution are impossible. Coalescence and breaking dynamics are the result of non-linear behaviors and will require significant extension to the current approach. This will be the subject of follow-up work. For the linear advection equation, one can check that the deformation map induced by the velocity, together with the initial condition, contains the full information for the solution of the advection equation. Our method discretizes the evolution of the deformation map on a regular Cartesian grid. The solution of the advection equation at any time can then be obtained as the pullback of the initial condition by this map. We call this approach the Characteristics Mapping (CM) method. This method is based on the work of Kohno and Nave \cite{kohno2013new} and relies on the Gradient-Augmented level-set (GALS) \cite{nave2010gradient} approach to evolve the map. Similar ideas have been used in various other contexts, most recently by Pons et al. \cite{pons2006maintaining} for point correspondence, Kamrin et al. \cite{kamrin2012reference} in the context of elasticity, and Ying and Cand\`{e}s \cite{ying2006fast} for the computation of geodesics.

The CM framework enables one to simultaneously evolve multiple advected quantities since a characteristic map captures the deformation generated by a given velocity field, and hence can be use to transport any quantity evolving under this flow. As a consequence, the evolution of multiple fluid interfaces only requires the computation of a single characteristic map, resulting in a fast and internally coherent transport of surfaces. Furthermore, since the initial condition is not involved in the discretization and is simply rearranged by the characteristic map, the method achieves arbitrary fine spatial resolution due to the diffeomorphism property of the map: even fractal domains can be transported and arbitrary subgrid structures can be resolved. In this paper, we also make the fundamental observation that the velocity field and the global time deformation map do not contain the same scales of fine features. That is, even though a velocity field might be well represented on some given grid, the deformation map induced by this velocity can still develop arbitrarily fine features over time. This has led other methods \cite{mirzadeh2016parallel,kolomenskiy2016adaptive, deiterding2016comparison, li2015adaptive,behrens1996adaptive} to employ adaptive meshes in an explicit context (e.g. Quadtree/Octree structures) to represent fine evolving features. We instead apply this insight in an implicit framework by using a coherent two-grid strategy. The semigroup structure of the deformation maps allows us to separate the two scales present in the problem by using a coarse grid for the local-time advection and a fine grid of dynamic resolution for the global time representation of the deformation. As will be demonstrated, this allows for an efficient CM method that can advect arbitrary sets (closed, open, irregular, and even fractal) with high precision, arbitrary subgrid resolution, and advantageous computational times.
%Other methods employing dynamic spacial resolution include \cite{kolomenskiy2016adaptive, deiterding2016comparison, li2015adaptive,behrens1996adaptive} where adaptive meshes are used to represent fine features.

This paper is organized as follows: in section 2 we derive the mathematical formulation of the CM method, in section 3 we present the numerical implementation and provide pseudo-codes, in section 4 we perform several tests (including closed, open, and fractal initial conditions), and thoroughly examine accuracy and efficiency, and in section 5 we make some concluding remarks and propose future directions of work.

% !TEX root = CMmethod.tex
%\renewcommand{\chi}{\mathcal{X}}
\section{Mathematical Formulation} \label{sec:math}

\subsection{Characteristic mapping for advection}
In this section, we present the mathematical formulation for the Characteristic Mapping (CM) method. First we consider the advection equation, written in its strong form:
\begin{gather} \label{eq:advEqn}
	\partial_t \phi + \vec{v} \cdot \nabla \phi = 0 , \quad \phi(\vec{x}, 0) = \phi_0(\vec{x})
\end{gather}
where we assume $U \subset \mathbb{R}^d$ is a rectangular domain in $d$ dimensions and $\phi  : U  \times [0, T] \to  \mathbb{R}$ is the advected quantity. Throughout the paper, we assume the given velocity field $\vec{v}(\vec{x}, t)$
\begin{gather}
	\vec{v} : U \times [0, t] \to  \mathbb{R}^d  \label{eq:defV}
\end{gather}
to be smooth and divergence-free.

We apply the method of characteristics to (\ref{eq:advEqn}) . Let $\vec{x}$ be the solution to the following initial value problem (IVP):
\begin{gather} \label{eq:charCurve}
\frac{d}{dt} \vec{x}(t) = \vec{v}(\vec{x}(t), t),  \quad \vec{x}(0) = \vec{x}_0
\end{gather}
where $\vec{x}$ are the characteristic curves of the equation generated by the velocity $\vec{v}$. 

The method of characteristics consists in realizing that \eqref{eq:advEqn} simplifies to the following equations \eqref{eq:genTransportEqnDt} and \eqref{eq:genTransportEqn}. We call \eqref{eq:genTransportEqn} the \emph{generalized} advection equation for some given velocity field $\vec{v}$ generating the characteristics $\vec{x}(t)$. This is more general in that \eqref{eq:advEqn} is ill-posed when regularity assumptions on $\phi$ are dropped; we can therefore use \eqref{eq:genTransportEqn} to transport arbitrary functions:
\begin{gather} \label{eq:genTransportEqnDt}
\frac{d}{dt} \phi(\vec{x}(t), t) = 0
\end{gather}
or equivalently
\begin{gather} \label{eq:genTransportEqn}
\phi(\vec{x}(t), t) = \phi_0(\vec{x}_0)
\end{gather}
for every time $t$ and every initial condition $\vec{x}_0$.

The CM method consists in computing a backward-in-time solution operator to the IVP (\ref{eq:charCurve}). That is, we compute a map $\vec{X} : U \times [0, T] \to U$ such that for any curve $\vec{x}$ satisfying the IVP (\ref{eq:charCurve}) with some arbitrary initial condition $\vec{x}_0 \in U$, we have that 
\begin{gather}
\vec{X}(\vec{x}(t), t) = \vec{x}_0 .
\end{gather}
Existence and uniqueness for $\vec{X}$ under mild regularity assumptions have been established in \cite{filippov2013differential}.

One can check that $\vec{X}$ is the solution to the vector valued advection equation
\begin{gather} \label{eq:mapsEvolve}
\partial_t \vec{X} + (\vec{v} \cdot \nabla) \vec{X} = 0, \quad \vec{X}(\vec{x}, 0) = \vec{x} .
\end{gather}
We call $\vec{X}$ the \emph{characteristic map} associated to the velocity $\vec{v}$. Under suitable regularity assumptions on $\vec{v}$, $\vec{X}$ is known to be a diffeomorphism. 

We see that the map $\vec{X}(\cdot, t)$ takes a point $\vec{x}(t)$ at time $t$ on its characteristic curve and returns the initial position $\vec{x}_0$ of this curve. This allows us to obtain the solution of the advection equation (\ref{eq:advEqn}) as
\begin{gather} \label{eq:pullBackInit}
\phi(\vec{x}, t) = \phi_0 \left(\vec{X}(\vec{x}, t) \right)
\end{gather}

This framework provides solutions to a generalization of equation \eqref{eq:advEqn}. Any $\phi$ satisfying equation \eqref{eq:genTransportEqn} can be expressed using \eqref{eq:pullBackInit}. This allows for the transport of more general sets, since $\phi$ is allowed to be any function defining a transported set. The lack of regularity assumptions on $\phi$ implies that arbitrary sets, even fractal ones, can be advected under this framework. For instance, taking $\phi_0$ to be the indicator function of an initial set $S_0$, the pullback $\phi(\cdot, t) = \phi_0 \circ \vec{X}(\cdot, t)$ is the indicator function of
\begin{gather}
S_t = \vec{X}^{-1} \left( S_0, t \right)
\end{gather}
where $\vec{X}^{-1}$ denotes preimage and is valid even without the diffeomorphism assumption on $\vec{X}$.
\begin{remark}
Under the assumption that the map $\vec{X}$ generated by $\vec{v}$ is injective, the equation \eqref{eq:genTransportEqn} is always well-posed. Treatment of cases where $\vec{v}$ is non-smooth and of measure valued $\phi$ can be found in \cite{ambrosio2008existence}.
\end{remark}
\begin{figure}[h]
    \begin{center}
            \includegraphics[width=0.5\linewidth]{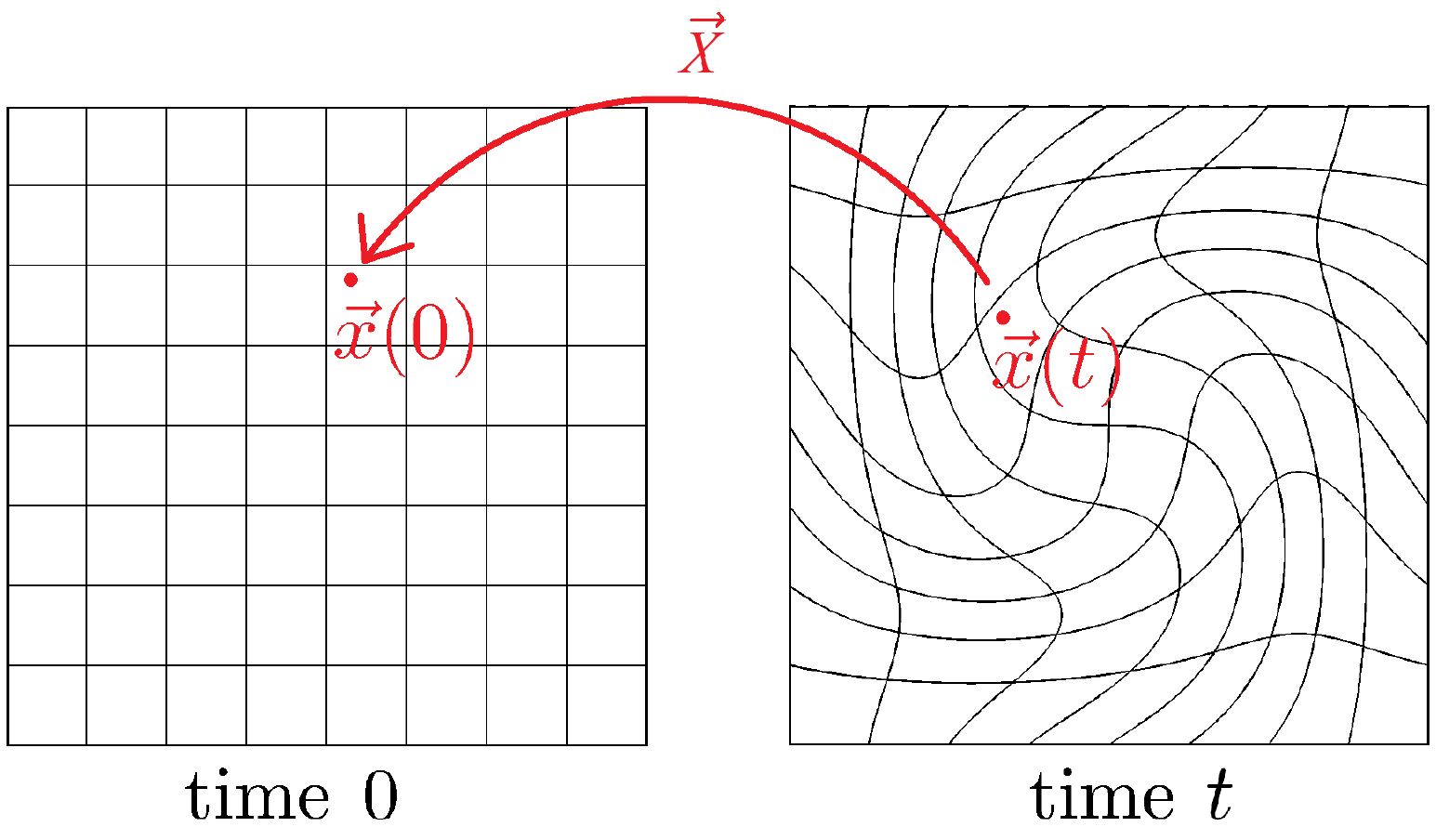}
    \end{center}	
    \caption{The diffeomorphism $\vec{X}$ takes a points $\vec{x}$ at time $t$ and returns its original position at time $t=0$.}
\label{fig:whatIsChi}
\end{figure}

The characteristic map $\vec{X}$ possesses a semi-group property which we exploit to significantly improve the numerical efficiency of the method. For this purpose, one can solve equation \eqref{eq:mapsEvolve} in an arbitrary subinterval $[\tau_{i-1}, \tau_{i}] \subseteq [0, T]$. The characteristic map for this subinterval satisfies
\begin{subequations} \label{eqGroup:submapsEvolve}
\begin{align}  
& \partial_t \vec{X}_i + \left(\vec{v} \cdot \nabla \right) \vec{X}_i = 0 \quad  \text{for } t \in [\tau_{i-1}, \tau_{i}]  \label{eq:chi_1} \\
& \vec{X}_i(\vec{x}, \tau_{i-1}) = \vec{x}
\end{align}
\end{subequations}

We denote by $\vec{X}_{[\tau_{i-1}, \tau_{i}]} : U \to U$ the solution of the above map at time $t = \tau_{i}$. That is, $\vec{X}_{[\tau_{i-1}, \tau_{i}]}$ represents the characteristic map generated by the velocity in the time interval $[\tau_{i-1}, \tau_{i}]$. The semi-group structure allows us to decompose the characteristic map for $[0, T]$ into arbitrarily many maps of subintervals. This presents a numerical advantage since it gives us an extra degree of freedom for adaptivity in time. The decomposition is given by:
\begin{gather} \label{eq:decomposition}
\vec{X}_{[0,T]} = \vec{X}_{[0, \tau_1]} \circ \vec{X}_{[\tau_1, \tau_2]} \circ \cdots \circ \vec{X}_{[\tau_{m-2}, \tau_{m-1}]} \circ \vec{X}_{[\tau_{m-1}, T]}
\end{gather}

We emphasize that the computation of the characteristics map in each individual subinterval $[\tau_i, \tau_{i+1}]$ start with the identity map as initial condition, and depends only on the velocity in this time interval. Hence, the computation is independent of any previous maps and can be solved individually.

\begin{remark}
Note that the time intervals $[\tau_i, \tau_{i+1} ]$ should be distinguished from the time step $\incr{t}$ of the numerical discretization. Whereas the latter is meant to be an approximation of the $\incr{t} \to 0$ limit, the time intervals of the submaps $\vec{X}_{[\tau_i, \tau_{i+1}]}$ are not small. The length of these intervals will be selected dynamically to distribute the task of representing the spacial features in $\vec{X}_{[0, T]}$ through the composition of several better-behaved submaps.
\end{remark}

\subsection{Characteristic map evolution using GALS framework}
For each subinterval, we compute the map by solving \eqref{eq:mapsEvolve} numerically. For this, we employ the Gradient Augmented level-set (GALS) method \cite{nave2010gradient}, an unconditionally stable semi-Lagrangian method that couples high order Hermite interpolation with Runge-Kutta time stepping schemes of matching order. Since the deformation map in CM is computed using GALS, and the advected quantity is obtained by pullback whenever a solution is required, it is a direct corollary that the CM method for linear advection is also unconditionally stable.

In this paper, we use Hermite cubic polynomials for spatial interpolation. We denote by $\He_\gC$ the projection operator from $C^1 (U)$ functions to Hermite cubics on some grid denoted by $\gC$ with cell widths $\incr{x}$. To construct the Hermite cubic $\He_\gC[f]$, we evaluate $f$, $\nabla f$ and mixed partial derivatives of order at most one in each spacial dimension at grid points of $\gC$; $\He_\gC[f]$ is then defined in each cell to be the Hermite cubic interpolant with these grid data. For time discretization we use a fixed time step of $\incr{t}$. Time integration is performed using Runge-Kutta 3.

We denote by $\vec{\chi}_i(\vec{x}, t)$ the numerical approximation of $\vec{X}_i$ represented as a piecewise Hermite cubic. The Hermite cubic interpolant is updated using GALS stepping: we compute a one-step map $\vec{\Psi}$ which maps a point $\vec{x}$ at time $t$ to its position at time $t- \incr{t}$ following its characteristic curve in $[t- \incr{t}, t]$, then we compose $\vec{\Psi}$ with the characteristic map at time $t - \incr{t}$ and project the composition on the space of Hermite cubics. We denote by $\vec{\gamma}_{\vec{x}}(s)$ the characteristic curve satisfying \eqref{eq:charCurve} with $\vec{\gamma}_{\vec{x}}(t) = \vec{x}$ and compute $\vec{\Psi}$ by integrating the velocity along this curve:
\begin{subequations} \label{eqGroup:footpoints}
\begin{align} 
& \vec{\Psi}(\vec{x}) = \vec{x} - \int_{t- \incr{t}}^{t} \vec{v}(\vec{\gamma}_{\vec{x}}(s), s) ds \label{eq:footLoc} \\
& \vec{\chi}_i(\vec{x}, t) = \He_\gC \left[ \vec{\chi}_i(\vec{\Psi}(\vec{x}), t- \incr{t} ) \right] \label{eq:GALSStep}
\end{align}
\end{subequations}
The time integrations for $\vec{\gamma}_{\vec{x}}$ and $\vec{\Psi}$ are approximated using RK3. The one-step map $\vec{\Psi}$ is evaluated when updating $\vec{\chi}_i$ in \eqref{eq:GALSStep}; the spacial derivatives at grid points required for the definition of the Hermite interpolant are obtained by direct chain rule: derivatives of $\vec{\chi}_i(\vec{x}, t - \incr{t})$ are computed from its Hermite interpolant definition and those of $\vec{\Psi}$ are directly evaluated from \eqref{eq:footLoc}.

Our strategy is to use the above time-stepping to evolve submaps $\vec{\chi}_i$ over some time interval $[\tau_{i-1}, \tau_{i}]$ which is either selected \emph{a priori} or determined dynamically. Once time $\tau_{i}$ is reached, we stop the submap evolution for $\vec{\chi}_{[\tau_{i-1}, \tau_{i}]}$ and inductively store the composition $\vec{\chi}_{[0, \tau_{i-1}]} \circ \vec{\chi}_{[\tau_{i-1}, \tau_{i}]}$ as a single Hermite cubic $\vec{\chi}_{[0, \tau_i]}$ of high enough resolution; this is called the \emph{remapping step}. We then repeat with the next submap $\vec{\chi}_{[\tau_i, \tau_{i+1}]}$.

To improve efficiency, we compute the submaps $\vec{\chi}_{[\tau_{i-1}, \tau_{i}]}$ on a coarse grid, since submap updates happen often, and a low spatial resolution is sufficient to represent the small deformations. On the other hand, the spatial resolution required to represent the global map is high since large deformations can accumulate over time. We therefore store the global map $\chi_{[0, \tau_i]}$ on a fine grid. Updating the global map is costly, however these updates only happen at the remapping steps. Hence we must select the submap time intervals appropriately as to not incur too many remapping steps and remap only when the coarse grid is no longer able to accurately resolve the submap deformations. The interplay between the submap time intervals and the coarse and fine grid resolutions will be analyzed in the following sections.

To summmarize:
\begin{enumerate}
	\item[] Given $0 = \tau_0 < \tau_1 < \ldots < \tau_{m-1} < \tau_m = T$. 
	\item[1.] Initialize $\vec{\chi}_{[0, \tau_0]} (\vec{x}) = \vec{x}$. For $i = 1$ to $m$, repeat the following 2 steps:
	\item[2.] Compute $\chi_{[\tau_{i-1}, \tau_{i}]}$ on a coarse grid $\gC$ using GALS.
	\item[3.] Compose and project  on a fine grid $\gF$, as Hermite cubic:
	\begin{gather} \label{eq:remappingStep}
	\vec{\chi}_{[0, \tau_{i}]} = \He_{\gF} \left[ \vec{\chi}_{[0, \tau_{i-1}]} \circ \vec{\chi}_{[\tau_{i-1}, \tau_{i}]} \right]
	\end{gather}
\end{enumerate}

\begin{remark}
At any intermediate time $t \in [\tau_{i-1}, \tau_i]$, the advected quantity can be evaluated by 
\begin{gather}
	\phi(\vec{x},t) = \phi_0 \left( \vec{\chi}_{[0, \tau_{i-1}]} \left( \vec{\chi}_i(\vec{x}, t) \right) \right)
\end{gather}
\end{remark}
\begin{remark}
Final time $T$ does not need to be known in advance. The evolution can continue for as long as desired.
\end{remark}

The computation of $\chi_{[\tau_{i-1}, \tau_{i}]}$ uses the GALS method. In each GALS step, the method produces two errors, one due to the approximate time integration and the other due to interpolation :
\begin{gather}
\text{GALS local truncation error} = \underbrace{\bigO(\incr{t}^4)}_\textrm{time integration} + \underbrace{\bigO(\incr{t}^2\incr{x}^2)}_\textrm{Hermite interpolation}
\end{gather}

Note that since Hermite cubics match the function value and $1^{\text{st}}$ derivative at grid points, the error scales quadratically with the distances between the query point and both the neighbouring grid points in each spacial dimension. For $\vec{x}$ a grid point, the evaluation of $\vec{\chi}_i(\vec{\Psi}(\vec{x}), t- \incr{t} )$ in \eqref{eq:GALSStep} occurs at query points located $| \vec{\Psi}(\vec{x}) - \vec{x}| = \bigO ( \incr{t} )$ and $\bigO (\incr{x})$ away from the closest grid points. Hence the $\bigO(\incr{t}^2\incr{x}^2)$ interpolation error.

For each $\chi_{[\tau_{i-1}, \tau_{i}]}$, we perform $\frac{\tau_{i} - \tau_{i-1}}{ \incr{t} }$ GALS steps. With each step accumulating the above error, the overall error for each submap is
\begin{gather}
\text{Submap global truncation error} = \bigO((\tau_{i}-\tau_{i-1})\incr{t}^3)+ \bigO((\tau_{i}-\tau_{i-1}) \incr{t} \incr{x}^2)
\end{gather}

After composition of the submaps, the global error due to submap evolution is 
\begin{gather}
\text{Global map error due to submap evolution} = \bigO(T\incr{t}^3)+ \bigO(T \incr{t} \incr{x}^2)
\end{gather}

we see that for any number of time subdivisions, the composition of the submaps yields a globally third-order accurate method.

In the next section, we will discuss the numerical properties specific to this third-order method. Namely, we will look at the effect of the choice of coarse and fine grids, as well as the dynamic adaptation of time interval length and fine grid resolution.

% !TEX root = CMmethod.tex
\section{Numerical Implementation}
In this section, we detail the numerical implementation of the method described in section \ref{sec:math}. The first subsection describes adaptive decomposition of $[0, T]$ into subintervals, the second subsection deals with the use of dynamic grid resolution for the storage of $\vec{\chi}_{[0, \tau_i]}$, and the third subsection summarizes the CM method in pseudo-code for ease of implementation.

\subsection{Adaptive time subdivision} \label{sec:adaptiveTime}

In practice, the subdivision times $\tau_i$ are chosen adaptively and not predetermined. We define a measure of the error of the submap computation $\vec{\chi}_i$ and take $\tau_{i}$ to be the first time this error exceeds some predetermined tolerance.

To evaluate this error, we use Lagrangian particles $p \in \{1, \dots, m\}$, initially distributed uniformly in $U$ at positions $\vec{y}_p^{\,0}$. These particles are independently evolved under $\vec{v}$ by solving the set of ODEs
%\begin{subequations} \label{eqGroup:particle}
%\begin{align}
%	\frac{\partial \vec{y}_p(t)}{\partial t} &= \vec{v}(\vec{y}_p(t), t), \quad \forall \, t \geq \tau_i \quad \forall \, p \in \{1, \dots, m\} \label{eq:particle}\\
%	\vec{y}_p(\tau_i) &= \vec{y}_p^{\,0}
%\end{align}
%\end{subequations}
\begin{gather} \label{eq:particle}
\frac{\partial \vec{y}_p(t)}{\partial t} = \vec{v}(\vec{y}_p(t), t), \quad \vec{y}_p(\tau_i) = \vec{y}_p^{\,0}
\end{gather}
We use RK3 integration scheme to solve \eqref{eq:particle} for particles $\vec{y}_p$. A measure of the accumulated interpolation error is given by
\begin{gather}
M_1(\vec{\chi}_i, t) := \max_p ||\vec{\chi}_i(\vec{y}_p(t), t) - \vec{y}_p^{\,0}|| \label{eq:normParticles}
\end{gather}

Equipped with the error measure $M_1$, we define $\tau_{i}$ as the first time at which the evolution of $\vec{\chi}_i$ induces a representation error greater than some predefined tolerance $\mathcal{E}_1$, that is
\begin{gather}
M_1(\vec{\chi}_i, \tau_{i} - \incr{t}) \leq \mathcal{E}_1 < M_1(\vec{\chi}_i, \tau_{i}) \label{eq:tolerance_1_1_1}
\end{gather}

At this time, we stop evolving $\vec{\chi}_i$ and start computing the next submap. Notice that choosing $\mathcal{E}_1$ sufficiently small ensures that the higher frequencies omitted by the coarse grid representation do not become large. Indeed, even though the grid values of $\chi_i$ may be very precise (for instance, by taking very small $\incr{t}$), the Hermite cubic representation can still be bad due to a lack of resolution. In this case, projecting and composing on a fine grid leads to high error for the fine grid values. This is why the error measure we use takes into account the entire domain and controls error at off-grid points. Consequently, we may oversample $\vec{\chi}_{[\tau_{i-1}, \tau_{i}]}$ and compose it with $\vec{\chi}_{[0, \tau_{i-1}]}$ onto a finer grid with controlled error.

\begin{remark}
To ensure that particles are always distributed uniformly over the domain $U$, we reinitialize their position at every remapping step, that is $\vec{y}_p(\tau_{i}) = \vec{y}_p^{\,0}$.
\end{remark}

\subsection{Dynamic grid resolution} \label{sec:dynamicGrid}
The entire algorithm utilizes two grids: A coarse grid $\gC$ of $N_c$ points per dimension where the $\vec{\chi}_i$ submaps are evolved, and a fine grid $\gF$ of $N_f$ points per dimension where the accumulated global maps $\vec{\chi}_{[0, \tau_i]}$ are stored.

As seen in the previous subsection, $\mathcal{E}_1$ controls the importance of the omitted higher frequencies in the $\vec{\chi}_{[\tau_{i-1}, \tau_{i}]}$ submaps. These frequencies which cannot be represented on a $N_c$ grid arise either from the projection of the velocity onto the coarse grid or from advection effects on existing lower frequencies (e.g. shear effects). Similarly, when we perform the composition \eqref{eq:remappingStep},
%\begin{gather} \label{eq:globMapUpdate}
%\vec{\chi}_{[0, \tau_{i}]} = \He_{\gF} \left[ \vec{\chi}_{[0, \tau_{i-1}]} \circ \vec{\chi}_{[\tau_{i-1}, \tau_{i}]} \right],
%\end{gather}
advection effects will generate higher frequencies than what is originally present in $\vec{\chi}_{[0, \tau_{i-1}]}$. In this case, the original resolution of $N_f$ might not be enough to resolve these frequencies. Therefore, we dynamically adjust the fine grid resolution $N_f$ to control the error made by omitting higher frequencies.

To do this, we define another error measure. This one is intended to evaluate the capacity of a grid $\gF$ to properly represent some $C^1$ function $g$:
\begin{gather} \label{eq:M2GlobMap}
M_2(\gF, g) := \left\| g - \He_{\gF}[g]   \right\|_\infty = \left\| \left( \mathcal{I} -  \He_{\gF} \right) [g]   \right\|_\infty
\end{gather}
where the sup-norm is approximated by sampling the error at off-grid points.

We can then decide if a grid $\gF$ is fine enough to accurately represent the updated global map $\vec{\chi}_{[0, \tau_{i-1}]} \circ \vec{\chi}_{[\tau_{i-1}, \tau_i]}$ by evaluating $M_2(\gF, \vec{\chi}_{[0, \tau_{i-1}]} \circ \vec{\chi}_{[\tau_{i-1}, \tau_i]})$. If $M_2 \geq \mathcal{E}_2$, for some predefined tolerance $\mathcal{E}_2$ we will refine the fine grid, replacing $\gF$ by a grid $\gF^f$ of $2N_f$ points per dimension. We will compute $\vec{\chi}_{[0, \tau_i]}$ from \eqref{eq:remappingStep} using $\gF^f$ as $\gF$.

Conversely, we want to take advantage of situations where the fine grid could be coarsened. To detect such situations, we define a coarser grid $\gF^{c}$ of $N_f/2$ points per dimension and compute the $M_2$ error according to \eqref{eq:M2GlobMap}. If we have $M_2 \leq \mathcal{E}_2$, it means that the coarser grid is sufficient to represent $ \vec{\chi}_{[0, \tau_{i-1}]} \circ \vec{\chi}_{[\tau_{i-1}, \tau_i]} $. In this case, we coarsen the fine grid, replacing $\gF$ by a grid $\gF^c$ of $N_f/2$ points per dimension and compute $\vec{\chi}_{[0, \tau_i]}$ from \eqref{eq:remappingStep} using $\gF^c$ as $\gF$.

Using such a dynamic grid has two main advantages. By refining the grid when needed, we ensure that the transformation is always well represented, and by coarsening it when possible, we reduce the computational time required to do a remapping step. Note that the redefinition of $\vec{\chi}_{[0, \tau_i]}$ on a different grid is easily performed at low cost due to the Hermite cubic interpolant structure.

\subsection{Pseudo-code algorithm} \label{sec:algo}
We summarize here the CM algorithm with some minor adjustments at the implementation level. We use a single variable $\vec{\chi}_0$ for all global maps $\vec{\chi}_{[0, \tau_i]}$, as the old global map is overwritten every time it is updated. We also use a single variable $\vec{\chi}$ for all local submaps $\vec{\chi}_{[\tau_{i-1}, \tau_{i}]}$, as old submaps are no longer useful after the global map is updated. In the following, algorithm \ref{alg:GASM} describes the GALS method for evolving submaps, algorithm \ref{alg:CM} updates the global map by submap composition on a dynamically selected fine grid.
\begin{algorithm}[h]
\caption{GASM : GALS method for submap evolution} \label{alg:GASM}
\begin{algorithmic} 
\REQUIRE $\vec{v}, \tau_s, \gC, \incr{t}$, $\mathcal{E}_1$, $\vec{y}_p^{\,0}$
\ENSURE $\{ \vec{\chi}, \tau_f \}$
\STATE Initialization:  $\vec{\chi} \gets \He_\gC[\vec{x}]$ identity map, $t \gets \tau_s$, $\vec{y}_p \gets \vec{y}_p^{\,0}$
\WHILE{$\max_p  \left| \vec{\chi}(\vec{y}_p) - \vec{y}_p^{\,0} \right| < \mathcal{E}_1$}
\STATE $\vec{\Psi} \gets \vec{x} - \int^{t+ \incr{t}}_t \vec{v} (\vec{\gamma}(s), s) ds$ according to \eqref{eq:footLoc}
\STATE $\vec{\chi} \gets \He_\gC [ \vec{\chi} \circ \vec{\Psi} ]$ according to \eqref{eq:GALSStep}.
\STATE $\vec{y}_p \gets  \vec{y}_p +  \int_t^{t+ \incr{t}} \vec{v} (\vec{y}_p(s), s) ds$ according to \eqref{eq:particle}
\STATE $t \gets t+ \incr{t}$
\ENDWHILE
\STATE $\tau_f \gets t$
\RETURN $\{ \vec{\chi}, \tau_f \}$ 
\end{algorithmic}
\end{algorithm}
\begin{algorithm}[h]
\caption{The Characteristic Mapping Method} \label{alg:CM}
\begin{algorithmic} 
\REQUIRE $\vec{v}, \gC,  \gF, \incr{t}$, $\mathcal{E}_1$, $\mathcal{E}_2$, $\vec{y}_p^{\,0}$, $T$
\ENSURE $ \vec{\chi}_0$ 
\STATE Initialization: $\vec{\chi}_0 \gets \He_\gF[\vec{x}]$ identity map, $t \gets 0$
\WHILE{$t < T$}
\STATE $\tau_s \gets t$
\STATE $\{ \vec{\chi}, \tau_f \} \gets \text{GASM} \left[  \vec{v}, \tau_s, \gC, \incr{t}, \mathcal{E}_1, \vec{y}_p^{\,0}  \right] $ 
\STATE $\gF^f \gets $ grid 2$\times$ finer than $\gF$, $\gF^c \gets $ grid 2$\times$ coarser than $\gF$
\IF{$M_2(\gF, \vec{\chi}_0 \circ \vec{\chi}) > \mathcal{E}_2$}
\STATE $\gF \gets \gF^f$
\ELSIF{$M_2(\gF^c, \vec{\chi}_0 \circ \vec{\chi}) < \mathcal{E}_2$}
\STATE $\gF \gets \gF^c$
\ENDIF
\STATE $\vec{\chi}_0 \gets \He_{\gF} [\vec{\chi}_0 \circ  \vec{\chi}]$
\STATE $t \gets \tau_f$
\ENDWHILE
\RETURN $ \vec{\chi}_0$ 
\end{algorithmic}
\end{algorithm}
\section{Numerical Examples}

In this section, we evaluate the accuracy and performance of the CM method described previously. We use a grid of $N_c$ cells per dimension for the advection of $\vec{\chi}$ and a grid of $N_f$ cells per dimension to store $\vec{\chi}_0$. We will often compare the CM method to the GALS method, which will use a single grid having $N_g$ cells per dimension.

Our numerical experiments demonstrate the following numerical advantages of the CM method:
\begin{property}
The CM method is always initialized with the identity map, which has a much better behavior than the initial level-set function for a specific advection problem. As a result, sharp features generated by the deformation take longer to appear. This improves the performance of CM over GALS and other direct methods for comparable grid resolutions, even without remapping.
\end{property}
\begin{property}
The CM method optimizes efficiency through the separation of a coarse and fine grid. Short time deformation can be accurately resolved with low grid resolution. Using fine grids for frequent submaps updates does not significantly reduce error and is extremely costly and wasteful. The global time map contains sharp features formed from composition of short time maps. A dynamic fine grid allows us to resolve the global map correctly. To exploit this separation of scales, the CM method provides 4 parameters, $N_c$, $N_f$, $\mathcal{E}_1$ and $\mathcal{E}_2$ which, when picked correctly, can improve the efficiency of the algorithm.
\end{property}
\begin{property}
The interpolation structure of the algorithm implies that the evolution of the deformation map happens in the space of $C^1$ diffeomorphisms. This means that the map is available everywhere in the domain, allowing for arbitrary subgrid resolution. Since the method computes the solution operator of the advection equation instead of acting directly on the transported quantity, we can preserve the subgrid information for the initial condition and resolve arbitrarily fine features for all time. In practice, this means that the spacial resolution of the map is dissociated from the resolution of advected quantity. Since the initial condition of the advection is simply rearranged by the characteristic map during the simulation, a high resolution initial condition can be transported using a map of lower resolution while maintaining nonetheless the same high resolution in the solution at all times: the characteristic map only needs to resolve the dynamics of the flow, its spacial features is independent of the initial conditions of the advection equation.
\end{property}

We present standard benchmark tests in 2D and 3D in sections \ref{sec:swirl} and \ref{sec:3D}. We then apply the CM method to more complicated sets in section \ref{sec:complex}. Finally, we present accuracy and efficiency results for the method in section \ref{sec:CPU}.

\subsection{2D swirl test} \label{sec:swirl}

We apply the characteristic mapping method to the following 2D vector field taken from \cite{leveque1996high}
\begin{gather} \label{eqGroup:swirl_def}
\vec{v}(x, y, t) = \left( \begin{matrix*}[r]
\cos\left(\frac{\pi t}{A}\right) \left( \sin\left(\pi x\right) \right)^2\sin\left(2 \pi y\right) \\
-\cos\left(\frac{\pi t}{A}\right) \left( \sin\left(\pi y\right) \right)^2 \sin\left(2 \pi x\right)
\end{matrix*}  \right)
\end{gather}
with $A=16$ in the domain $[0,1]\times[0,1]$. This velocity creates a swirl centered at $(0.5, 0.5)$. The swirl reaches its maximal deformation at $t=8$ and then unwinds itself. At $t = 16$, the deformation induced in the first half of the period are completely undone and the transformation returns to identity.

The initial set is a circle of radius $0.15$ centered at $(0.5, 0.75)$ represented by a level-set function. From times $t=0$ to $t=16$, this flow will stretch the circle into a thin swirl and return it back to its original state.

\begin{figure}[h]
	\centering
%	\begin{center}
		\begin{subfigure}{0.2\linewidth}
			\includegraphics[width=\textwidth]{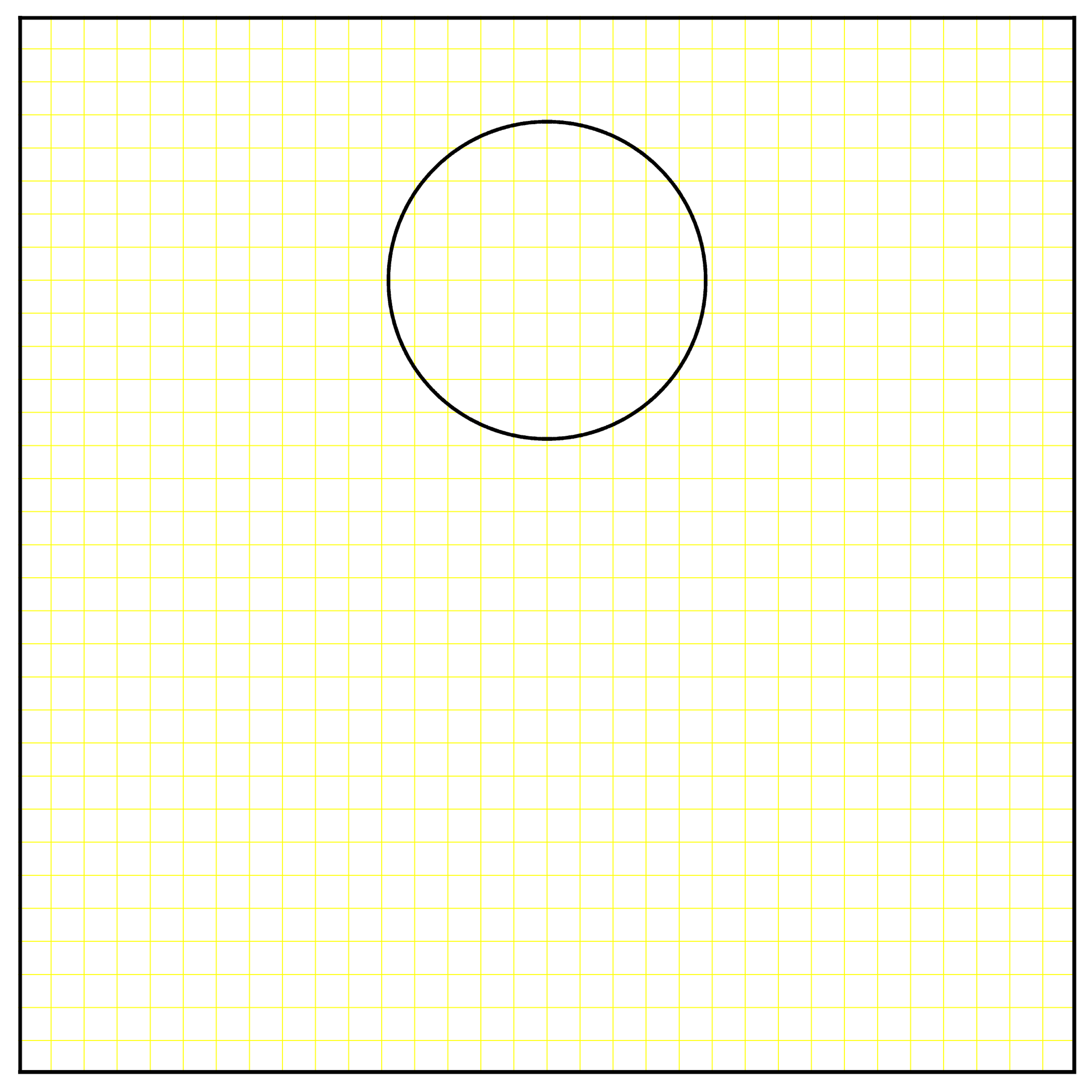}
			\caption{$t=0, N_f = 32$}	
		\end{subfigure}
		\begin{subfigure}{0.2\linewidth}
			\includegraphics[width=\textwidth]{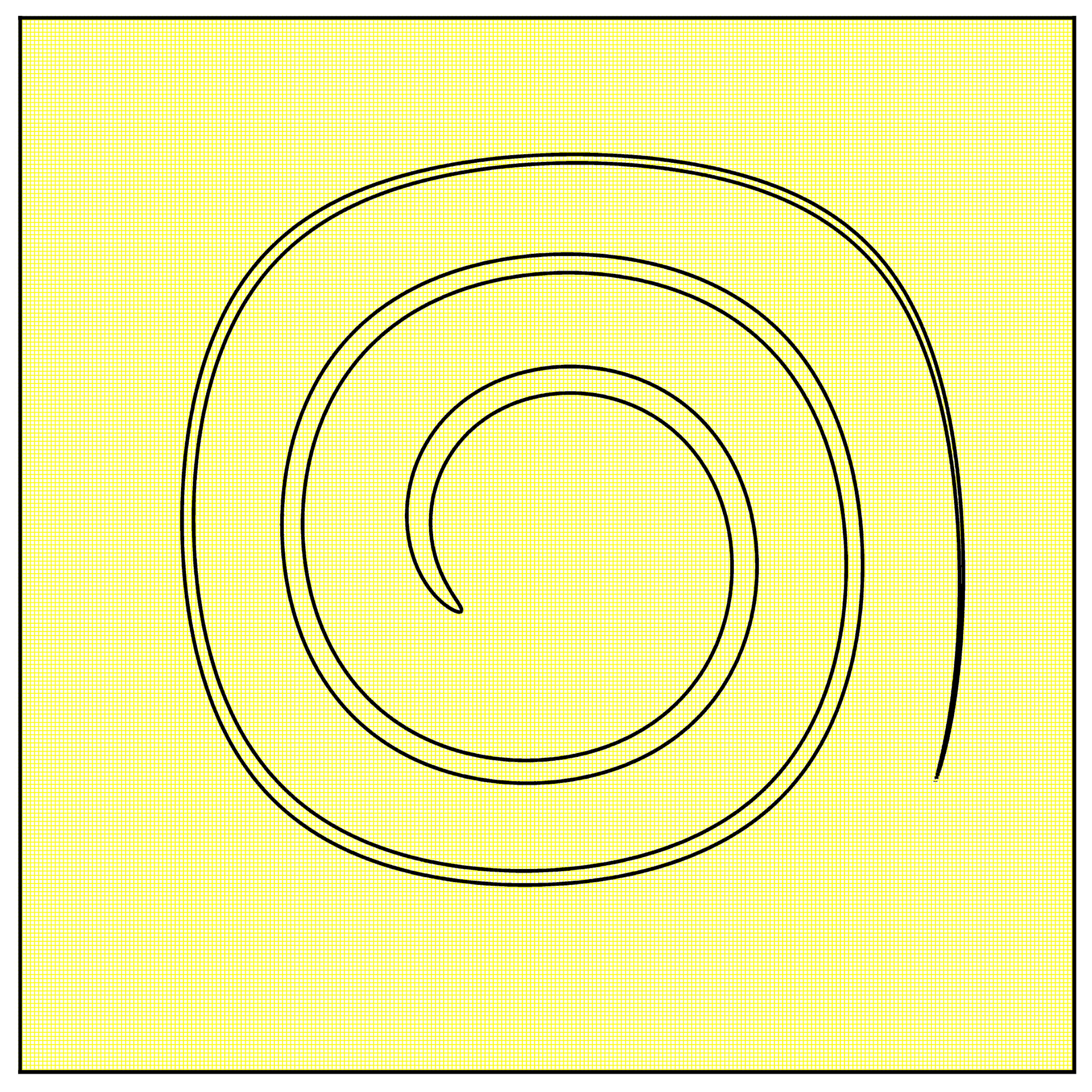}
			\caption{$t=4, N_f = 256$}	
		\end{subfigure}
		\begin{subfigure}{0.2\linewidth}
			\includegraphics[width=\textwidth]{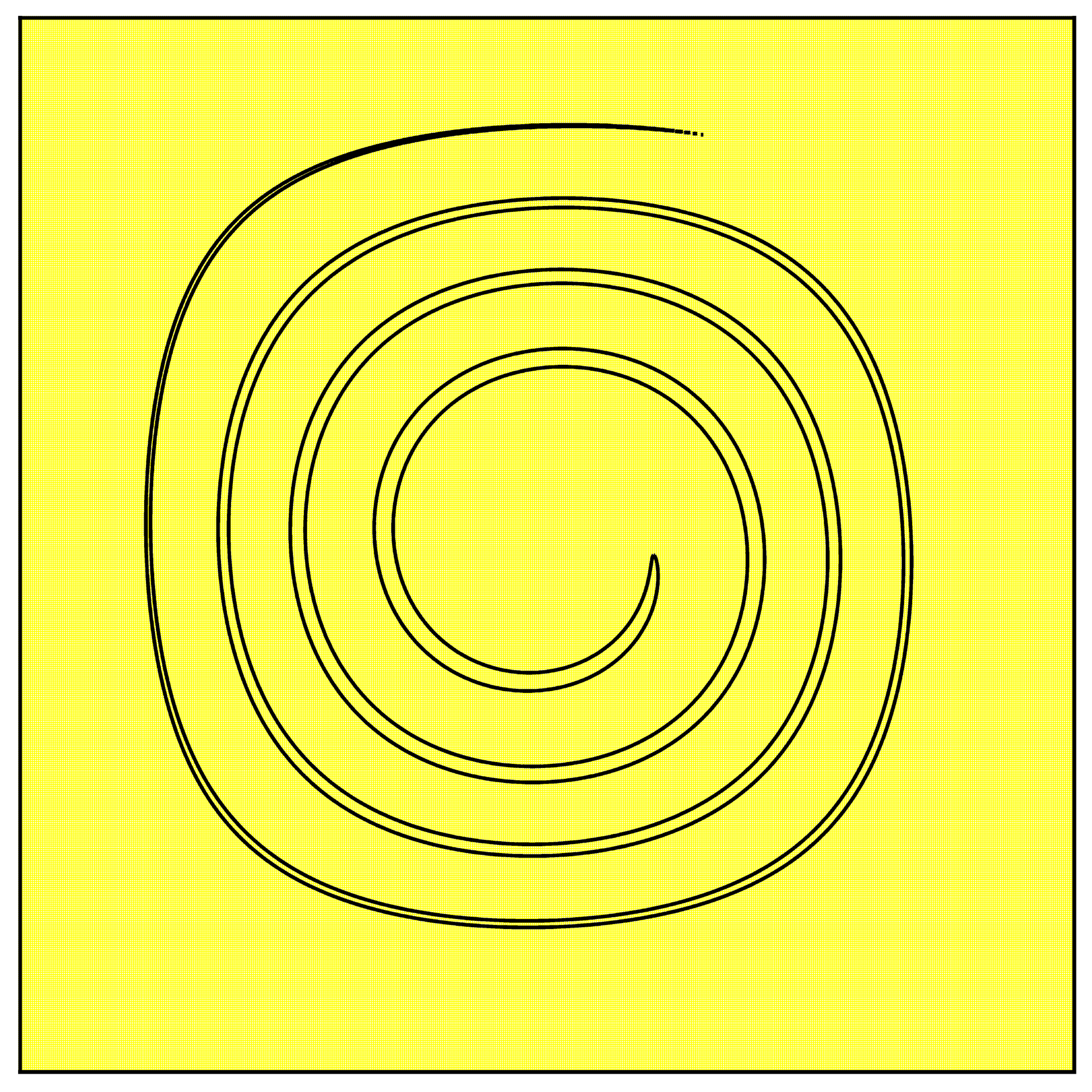}
			\caption{$t=8, N_f = 512$}	
		\end{subfigure} \\
		\begin{subfigure}{0.2\linewidth}
			\includegraphics[width=\textwidth]{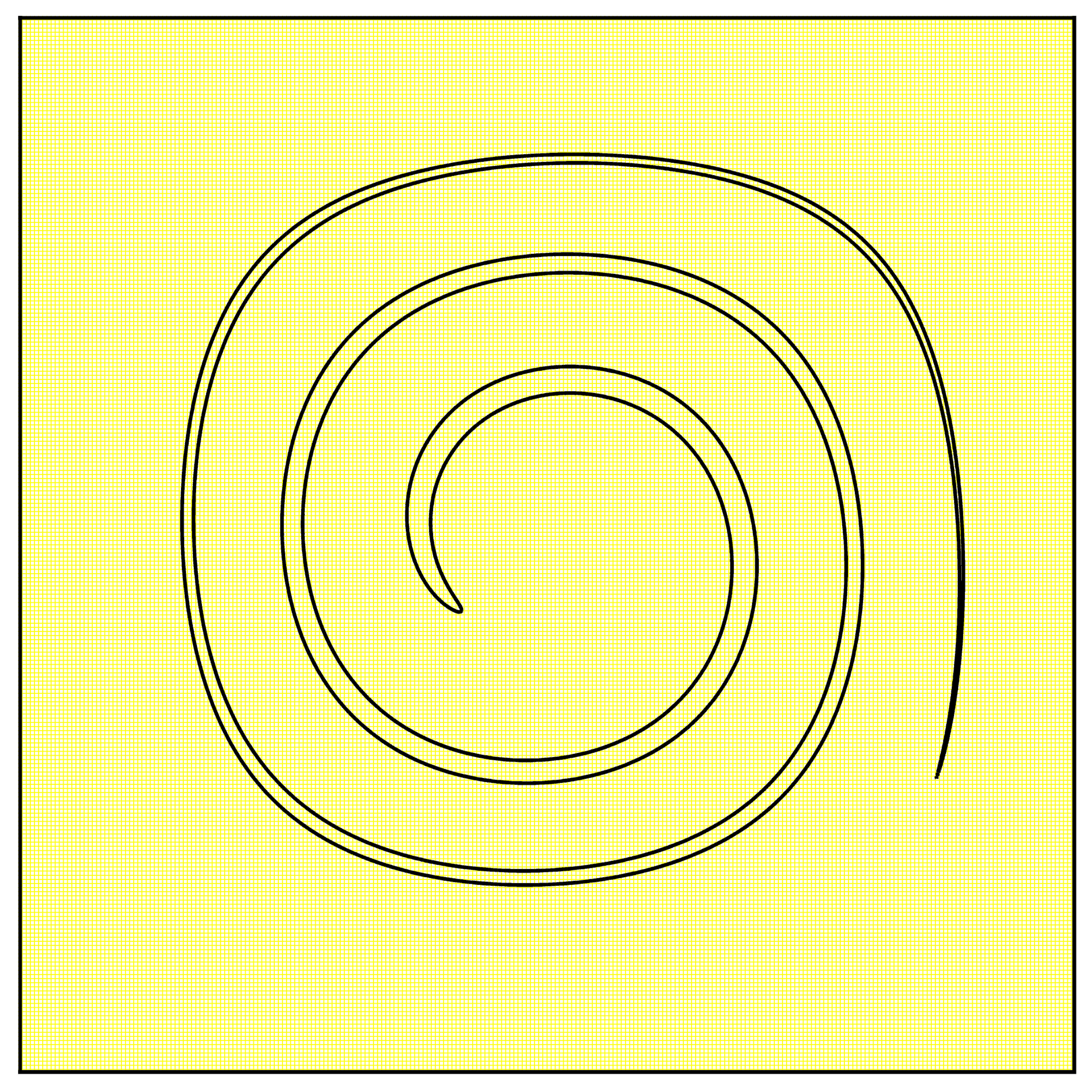}
			\caption{$t=12, N_f = 256$}	
		\end{subfigure}
		\begin{subfigure}{0.2\linewidth}
			\includegraphics[width=\textwidth]{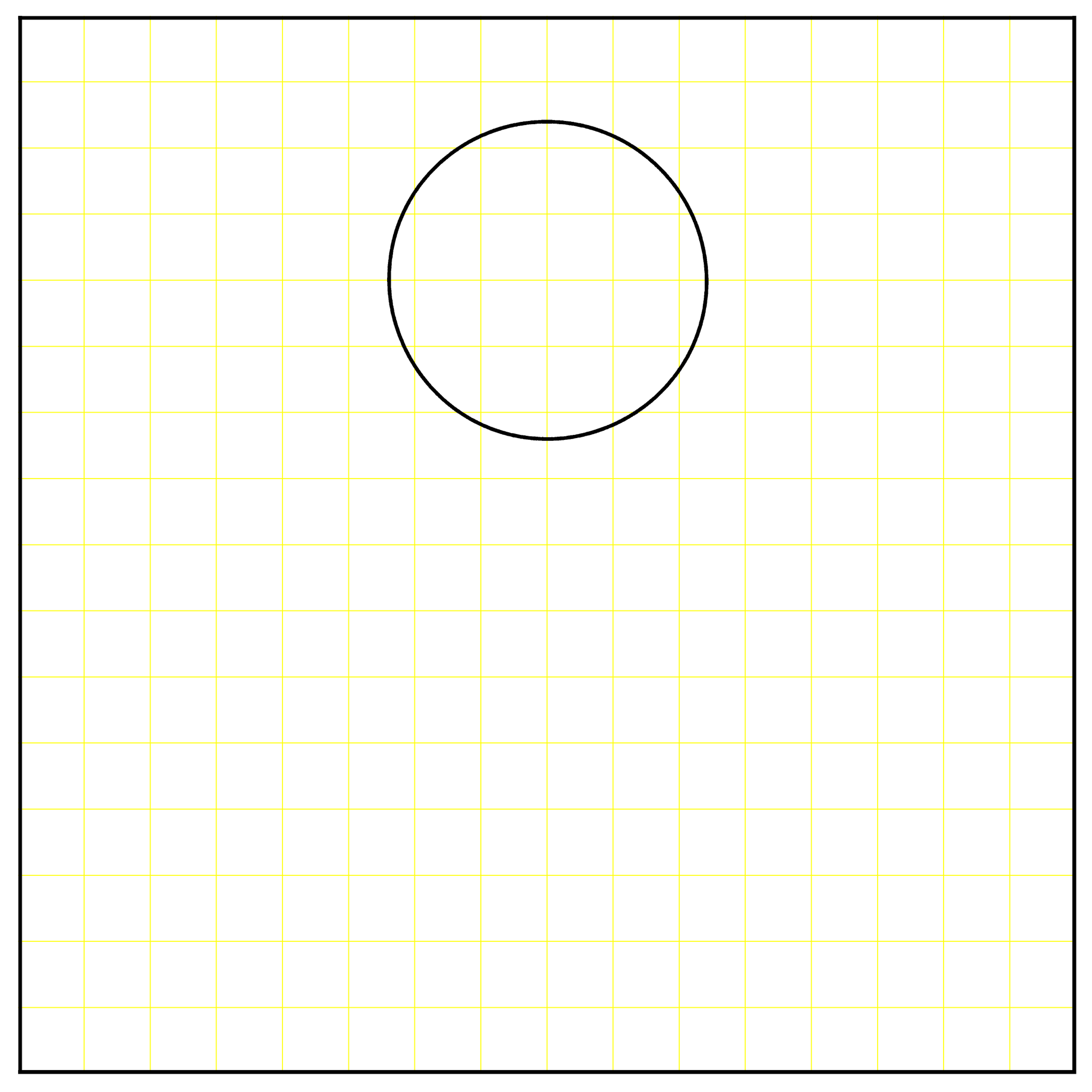}
			\caption{$t=16, N_f = 16$}	
		\end{subfigure}
%	\end{center}	
	\caption{2D swirl test using the characteristic mapping method with dynamic grid. Submaps are computed on a  $N_c = 32$ grid. Fine grid size is dynamic and capped at $N_f = 512$.}
	\label{fig:2Dswirl}
\end{figure}

We use a $N_c = 32$ grid for the advection and a dynamic-sized fine grid for the remapping. This fine grid starts at $N_f = 32$ and is refined to a maximum resolution of $N_f = 512$. The remapping and resizing tolerances were set to $\mathcal{E} _1=5\times10^{-6}$ and $\mathcal{E}_2 = 10^{-4}$. Figure \ref{fig:2Dswirl} shows the results.

The initial identity transformation is very well represented on the $32\times32$ grid, but as time advances, the flow creates very fine structures. The transformation at $t=8$ would not be well represented on a $32\times32$ grid. Indeed, the adaptive grid for $\vec{\chi}_0$ is refined to $512\times512$ at $t=8$, which is sufficient to represent correctly the sharp deformations induced by the vector field. Then as the swirl returns to its original shape, the grid coarsens. At $t = 16$, the transformation only needs a $16\times16$ grid to be correctly represented. We see that the initial circle is recovered at the final time, even though the advection was computed on the rather coarse $32\times32$ grid and only a few remapping steps involved finer grid calculations.

\begin{figure}[h]
	\begin{center}
		\begin{subfigure}{0.2\linewidth}
			\includegraphics[width=\textwidth]{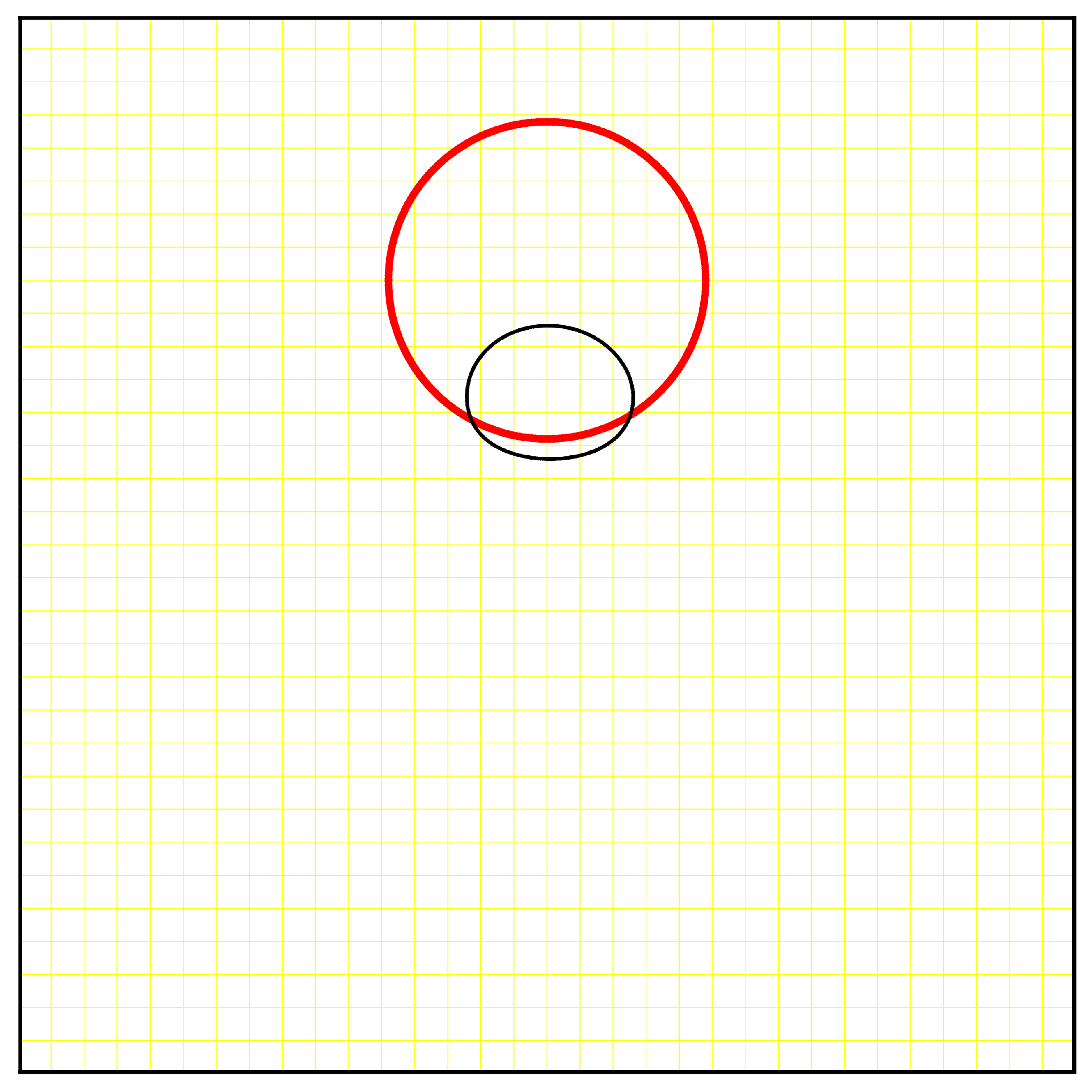}
			\caption{GALS \\ \hspace*{1.5em}$N_g = 32$}	
		\end{subfigure}
		\begin{subfigure}{0.2\linewidth}
			\includegraphics[width=\textwidth]{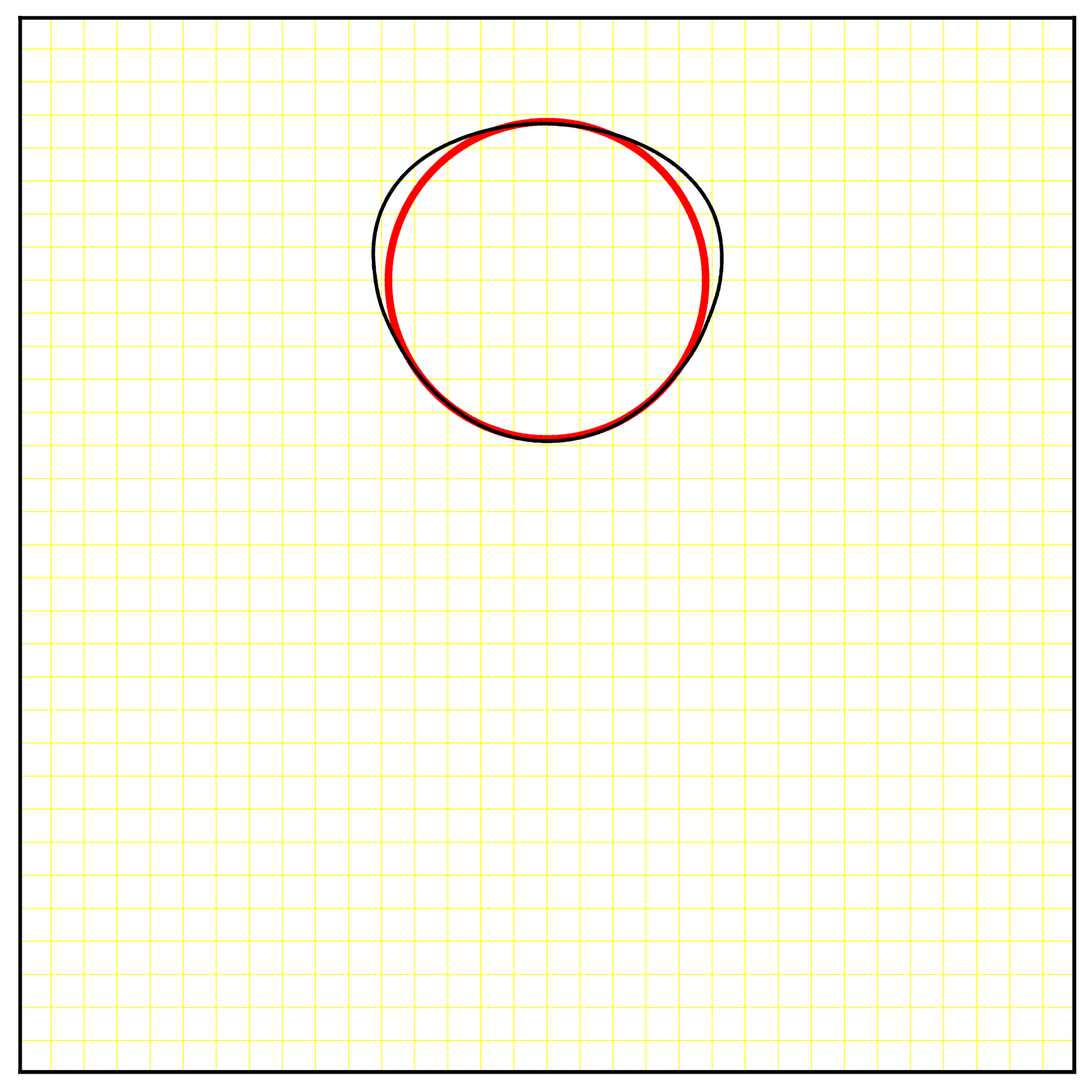}
			\caption{CM $N_c=32$ \\ \hspace*{1.5em} $N_f \leq 32$}	
		\end{subfigure} \hspace*{1.5em}
		\begin{subfigure}{0.2\linewidth}
			\includegraphics[width=\textwidth]{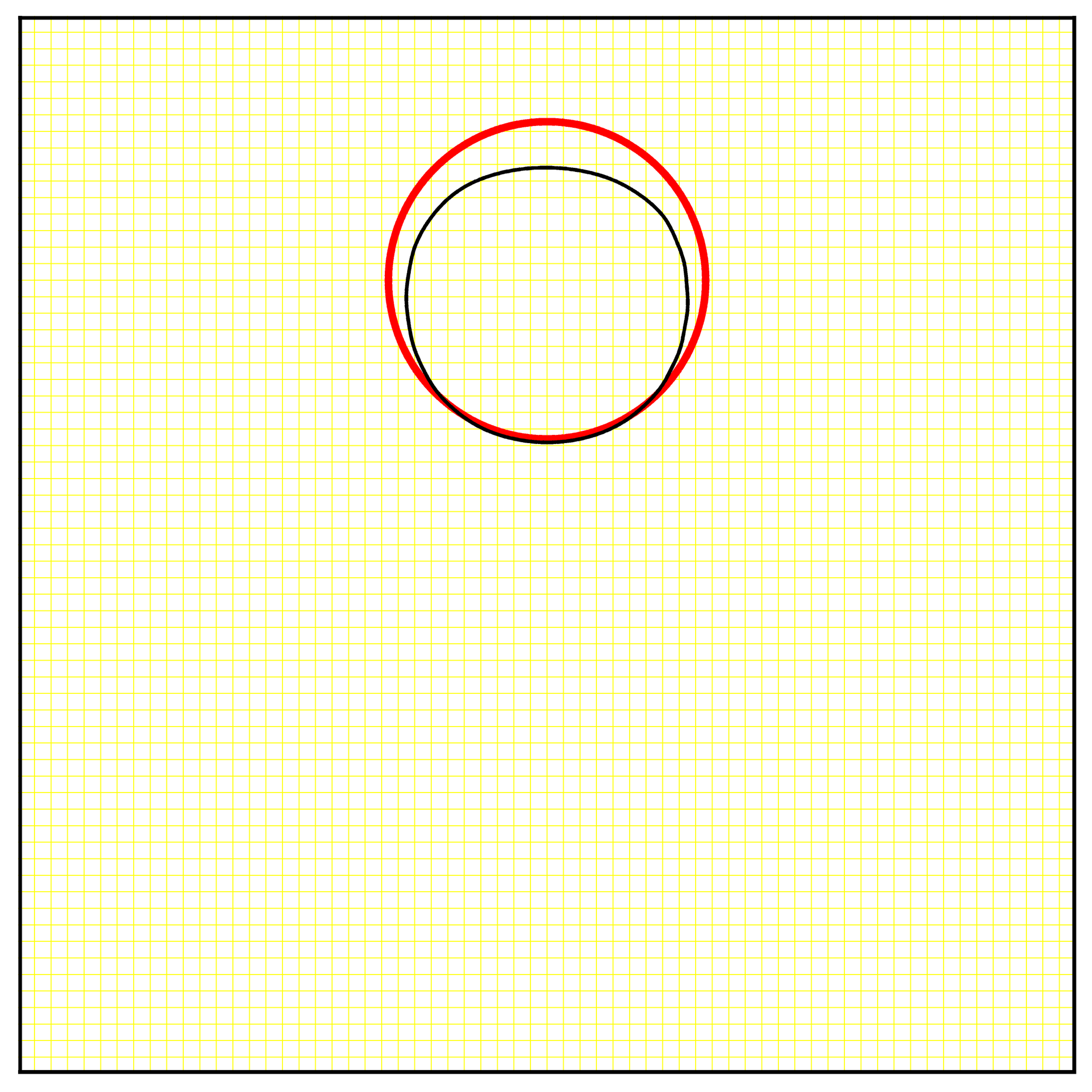}
			\caption{GALS \\ \hspace*{1.5em}$N_g = 64$}	
		\end{subfigure}
		\begin{subfigure}{0.2\linewidth}
			\includegraphics[width=\textwidth]{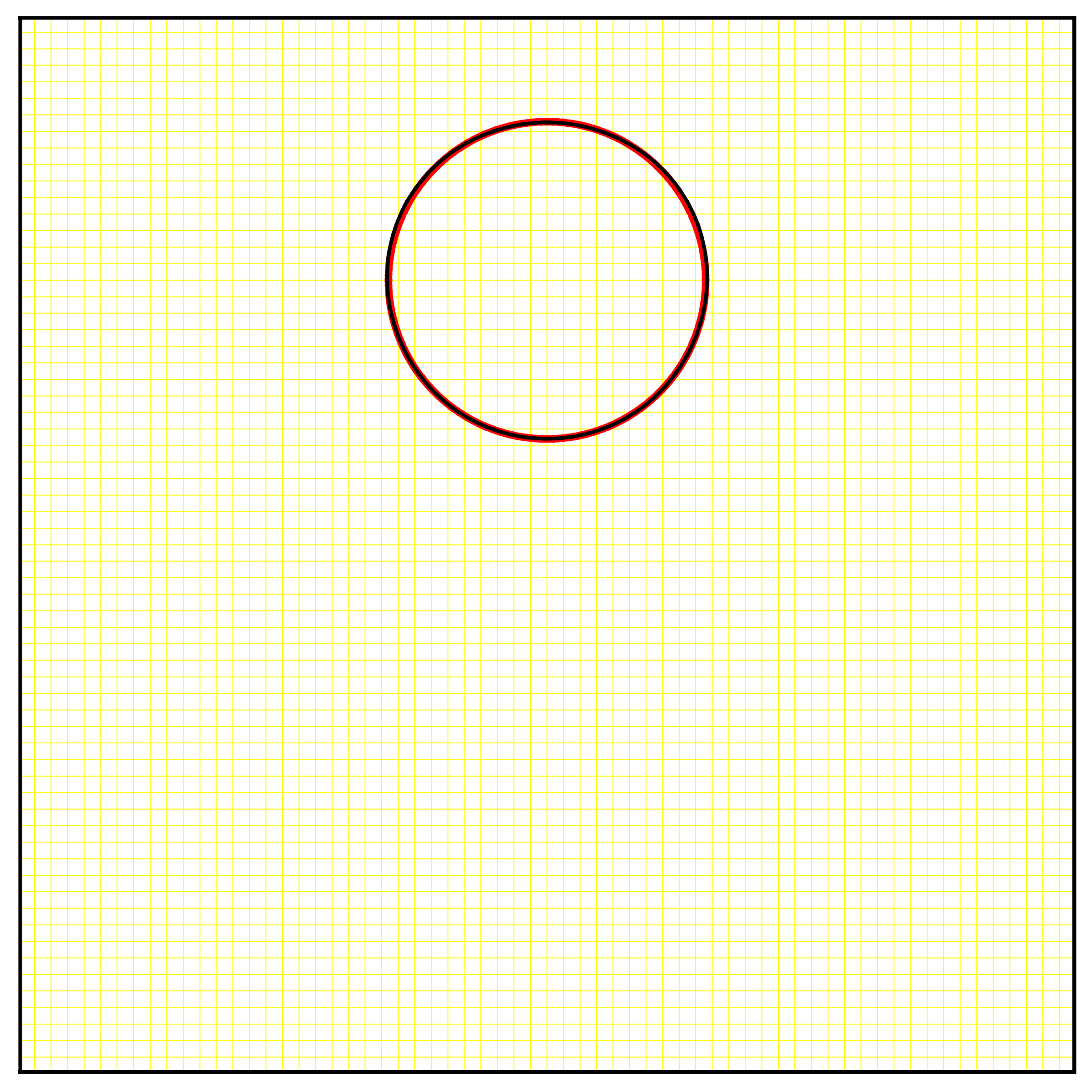}
			\caption{CM $N_c=32$\\ \hspace*{1.5em}$N_f \leq 64$}	
		\end{subfigure} 
		\begin{subfigure}{0.2\linewidth}
			\includegraphics[width=\textwidth]{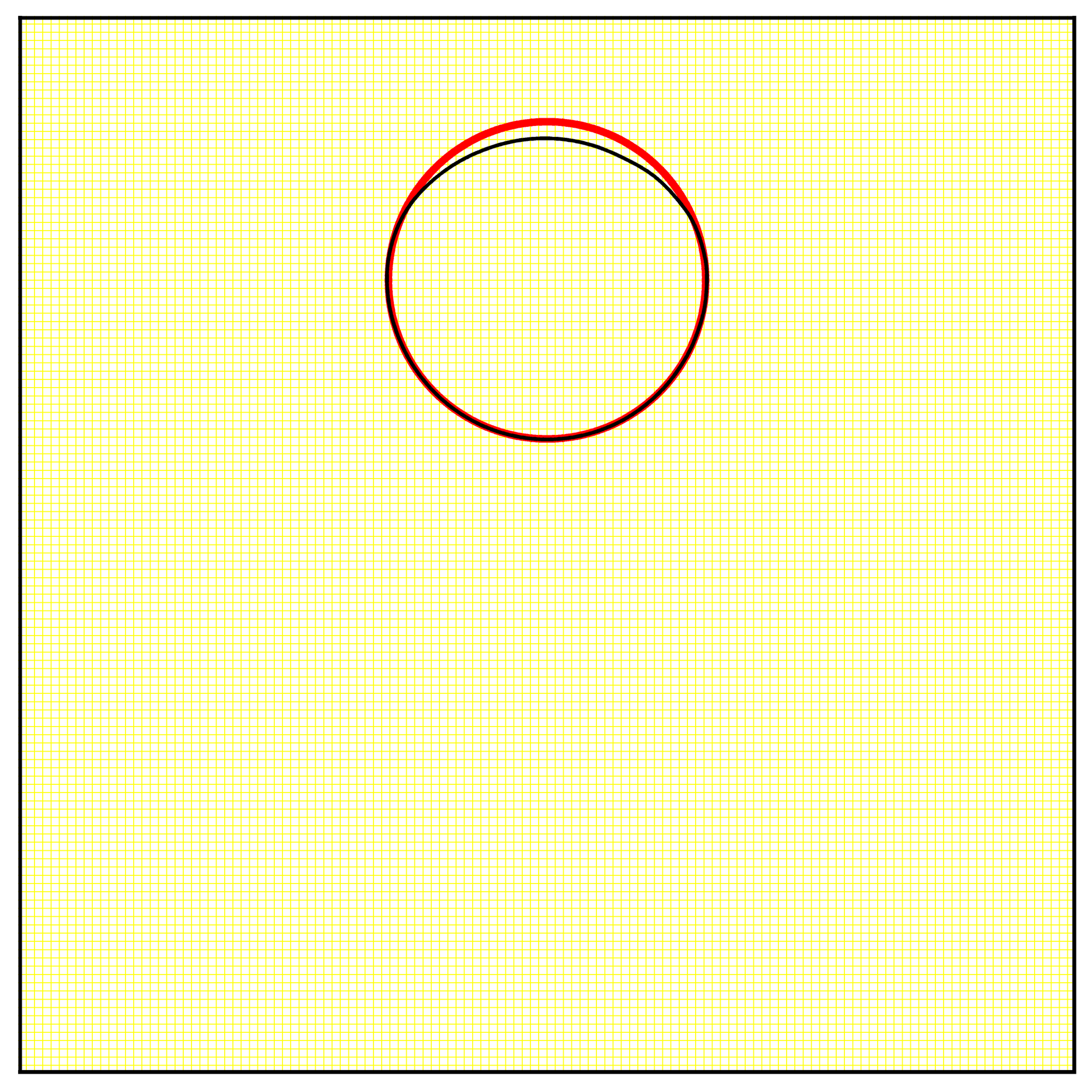}
			\caption{GALS \\ \hspace*{1.5em}$N_g = 128$}	
		\end{subfigure}
		\begin{subfigure}{0.2\linewidth}
			\includegraphics[width=\textwidth]{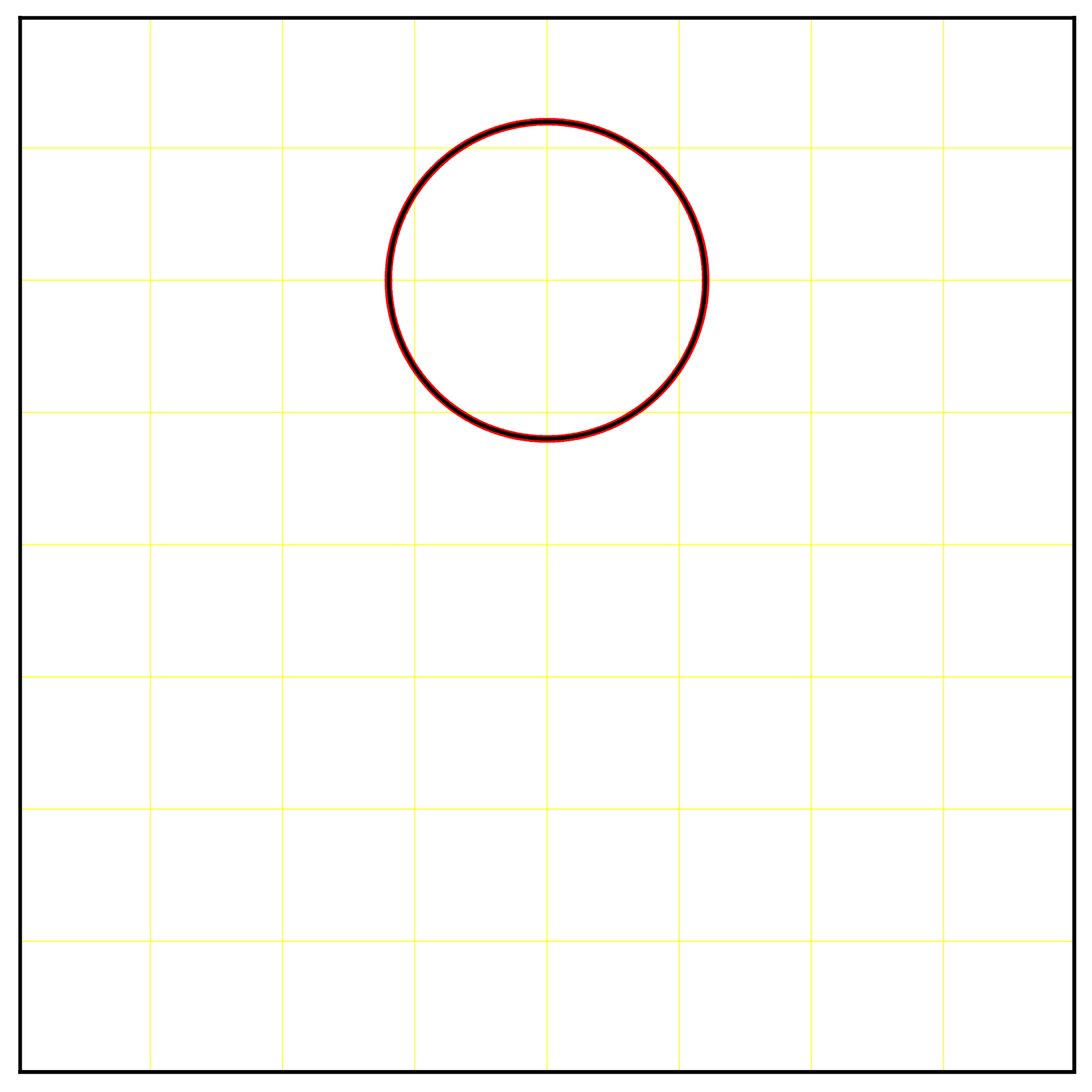}
			\caption{CM $N_c=32$\\ \hspace*{1.5em}$N_f \leq 128$}	
		\end{subfigure} \hspace*{1.5em}
		\begin{subfigure}{0.2\linewidth}
			\includegraphics[width=\textwidth]{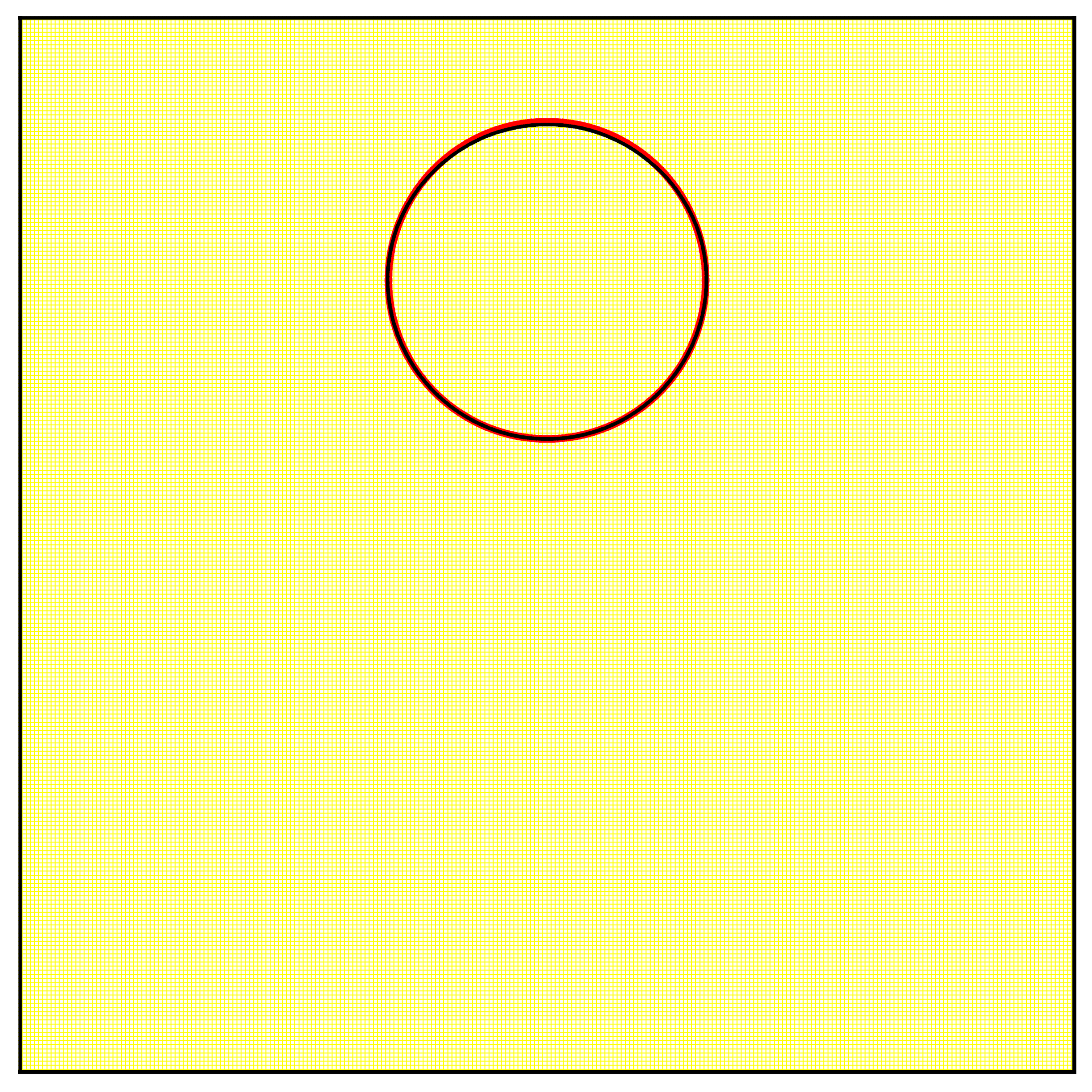}
			\caption{GALS \\ \hspace*{1.5em}$N_g = 256$}	
		\end{subfigure}
		\begin{subfigure}{0.2\linewidth}
			\includegraphics[width=\textwidth]{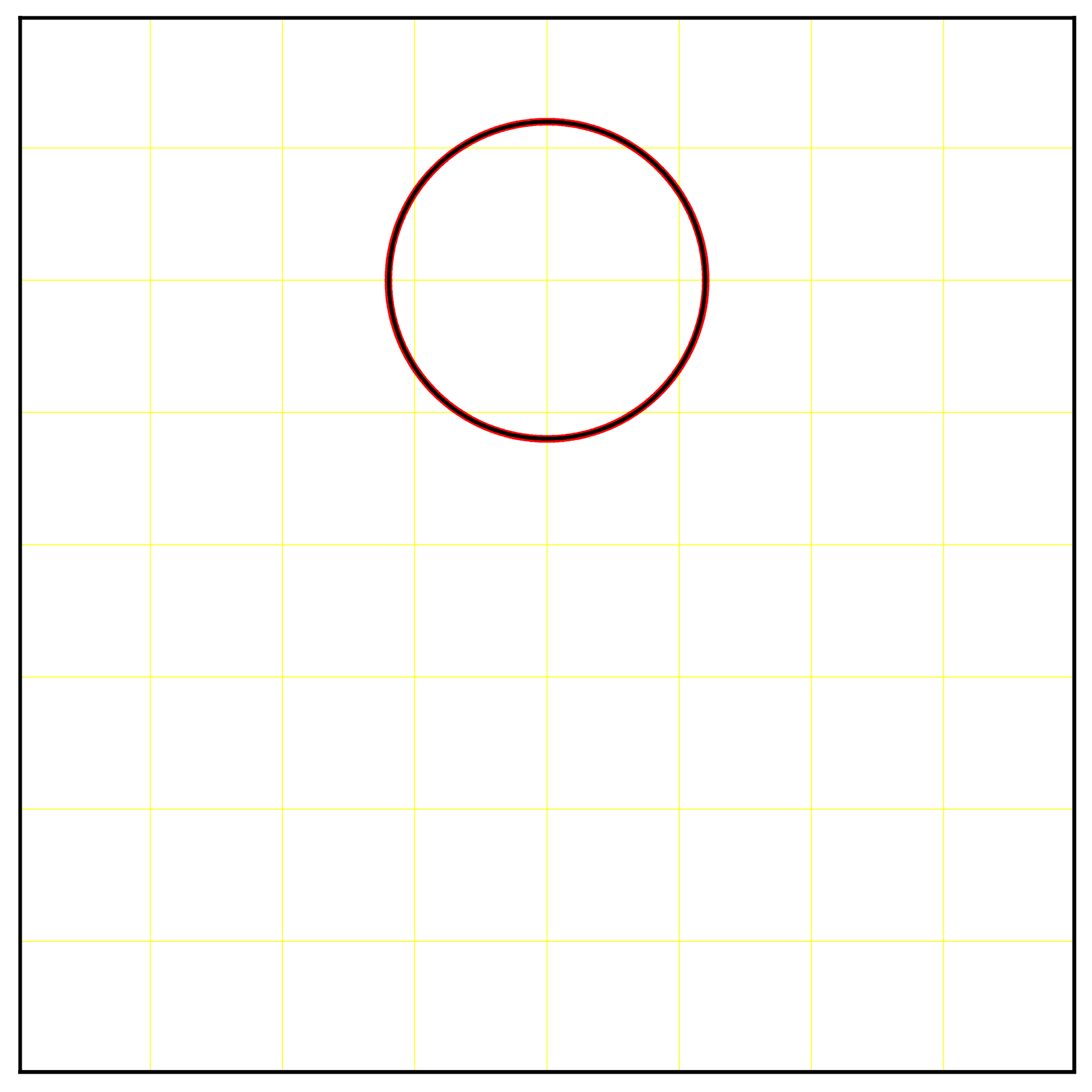}
			\caption{CM $N_c=32$\\ \hspace*{1.5em}$N_f \leq 256$}	
		\end{subfigure}
	\end{center}	
	\caption{Comparison between GALS and CM methods. The computed solution is shown in black and the exact solution in red.}
	\label{fig:GALSvsXYZ}
\end{figure}

We use a similar test to compare the characteristic mapping method with the standard GALS method \cite{nave2010gradient}. We use the same initial set and transport it in the vector field \eqref{eqGroup:swirl_def} with $A=8$ until $t=16$, which corresponds to the set being stretched and returned to its original shape twice. We make four different tests by modifying the grid sizes. For the CM method, we use a $N_c = 32$ for \emph{all} tests but we set the maximal $N_f$ to different values. The maximum allowed fine grid for each test are $N_f = 32, 64, 128, 256$ respectively. For the GALS method, we fix the grid resolution to be the same as for the finest possible remapping grid of the CM method, i.e. $N_g = 32, 64, 128, 256$. Results are shown in figure \ref{fig:GALSvsXYZ}.
\begin{table}[h]
\centering
\begin{tabular}{|c||c|c|c|c|}
\hline
\multicolumn{5}{ |c| }{Computational Time (sec.)} \\
\hline
\multirow{2}{*}{Method} & \multicolumn{4}{ c| }{$N_f$ or $N_g$} \\
\cline{2-5}
& \hspace{0.3cm}32\hspace{0.3cm} & \hspace{0.3cm}64\hspace{0.3cm} & \hspace{0.3cm}128\hspace{0.3cm} & \hspace{0.3cm}256\hspace{0.3cm} \\ 
\hline
\hline
GALS & 2 & 12 & 92 & 727    \\
\hline
CM with $N_c = 32$ & 6 & 7 & 9 & 11    \\ 
\hline
\end{tabular}
\caption{Time comparison (in seconds) of the GALS and CM methods for the swirl test (same test as in figure \ref{fig:GALSvsXYZ}). The GALS method is advected on the given grid sizes ($N_g$). The CM method uses a $N_c = 32$ grid for advection in all four cases, but the remapping is done on grids of the given sizes ($N_f$).}
\label{table:times}
\end{table}

We see that for a given grid size, the characteristic mapping method gives a better solution than the GALS method. This is in part due to the fact that the CM method starts from a linear initial condition $\vec{\chi}_0 = \vec{x}$ which can be exactly represented by Hermite cubics. Any feature of complexity higher than linear arise from the velocity field and accumulates at the rate of $\bigO (\incr{t} )$ per step, hence it takes longer for significant high frequencies to appear. In contrast, the GALS method has a more complicated initial condition. The level-set function already contains $\bigO(1)$ linear, quadratic and cubic features. Under advection effects, these can quickly generate sharp features that the grid cannot represent. Therefore, when the maximum resolution level $N_f = N_g$ is fixed, CM will be more accurate than GALS, as shown in figure \ref{fig:GALSvsXYZ}.

A particularly attractive feature of the CM method is that computation efficiency is optimized by the adaptive separation between coarse grid advection calculations and fine grid deformation representation.  Indeed, in a short time interval, deformations caused by the advection are small, and errors are dominated by velocity integration in time. It is therefore most efficient to perform these steps on a coarse grid: they happen frequently, but computation is cheap on a coarse grid, and the coarse grid is sufficient to represent the small deformation. In fact, doing submap updates on a fine grid would be wasteful since fine grid representation errors would be negligible compared to error from velocity integration. For the global time map, the error is dominated by difficulties in representing sharp features which requires a high grid resolution. However, global map updates do not happen often, and we can afford to use a fine grid representation.

Table \ref{table:times} compares the computational times for the GALS and CM method. For very coarse grids, the CM method is slower due to the overhead from processing the adaptivity. For the cases where the representation grid is finer than the advection grid, the CM method is clearly faster. The efficiency aspect is studied in more detail in section \ref{sec:CPU}.

\subsection{3D deformation field} \label{sec:3D}
The CM method generalizes directly to any number of dimentions. We apply the CM method to a 3D surface evolution problem. The initial surface is a sphere of radius $0.15$ centered at $(0.35, 0.35, 0.35)$ in the domain $[0,1]^3$. The deformation velocity is given by the following vector field taken from \cite{leveque1996high}
\begin{gather} \label{eqGroup:3D_def}
\vec{v}(x, y, z, t) = \left( \begin{matrix}
2\cos\left(\frac{\pi t}{2}\right) \left( \sin(\pi x) \right)^2 \sin(2 \pi y) \sin(2 \pi z)\\
-\cos\left(\frac{\pi t}{2}\right) \sin(2 \pi x) \left( \sin(\pi y) \right)^2 \sin(2 \pi z)\\
-\cos\left(\frac{\pi t}{2}\right) \sin(2 \pi x) \sin(2 \pi y) \left( \sin(\pi z) \right)^2
\end{matrix}  \right)
\end{gather}

\begin{figure}[h]
	\begin{center}
		\begin{subfigure}{0.26\linewidth}
			\includegraphics[width=0.9\textwidth]{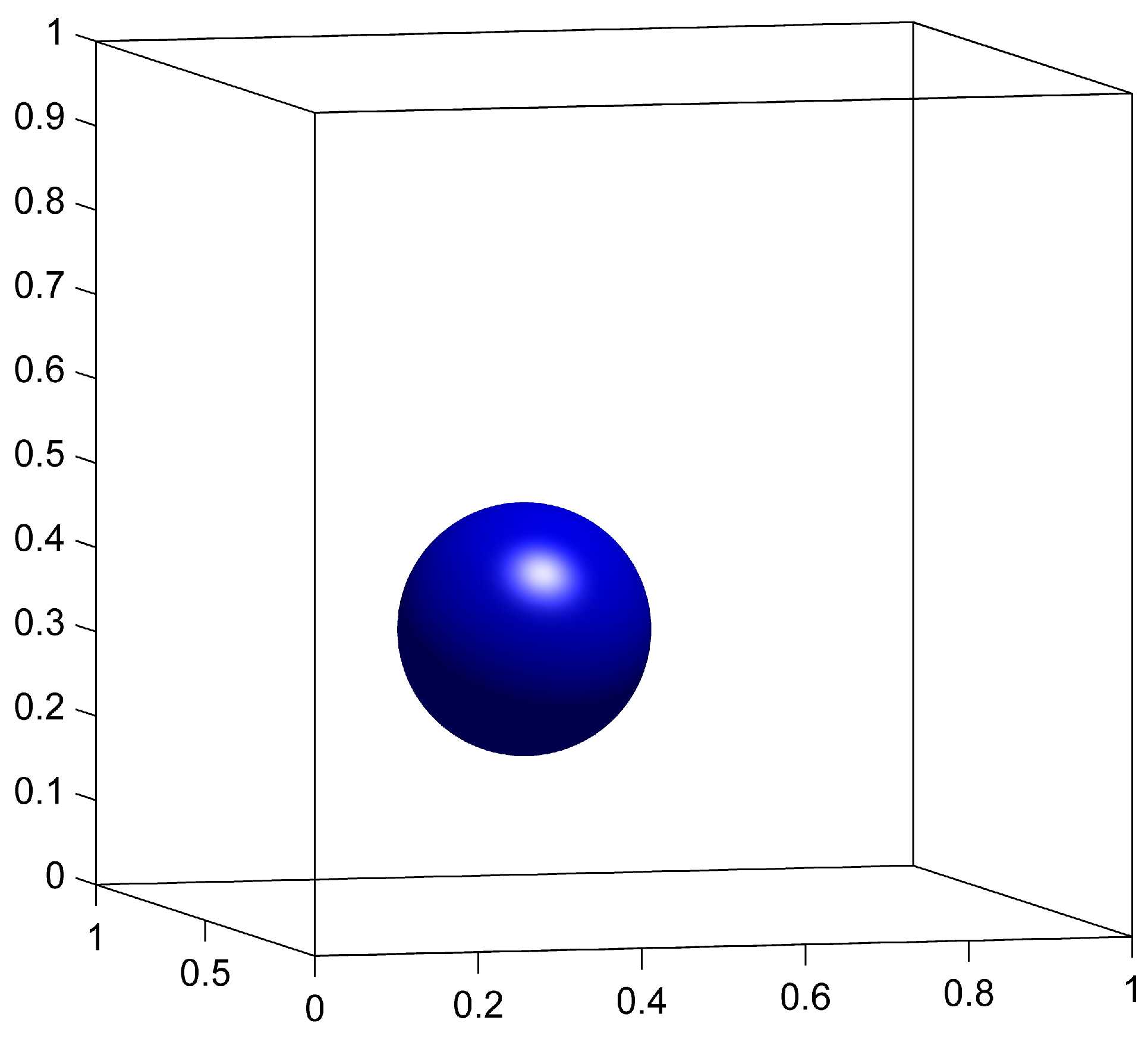}
			\caption{CM without remap, t=0}	 \label{fig:3Dcn0}
		\end{subfigure}
		\begin{subfigure}{0.26\linewidth}
			\includegraphics[width=0.9\textwidth]{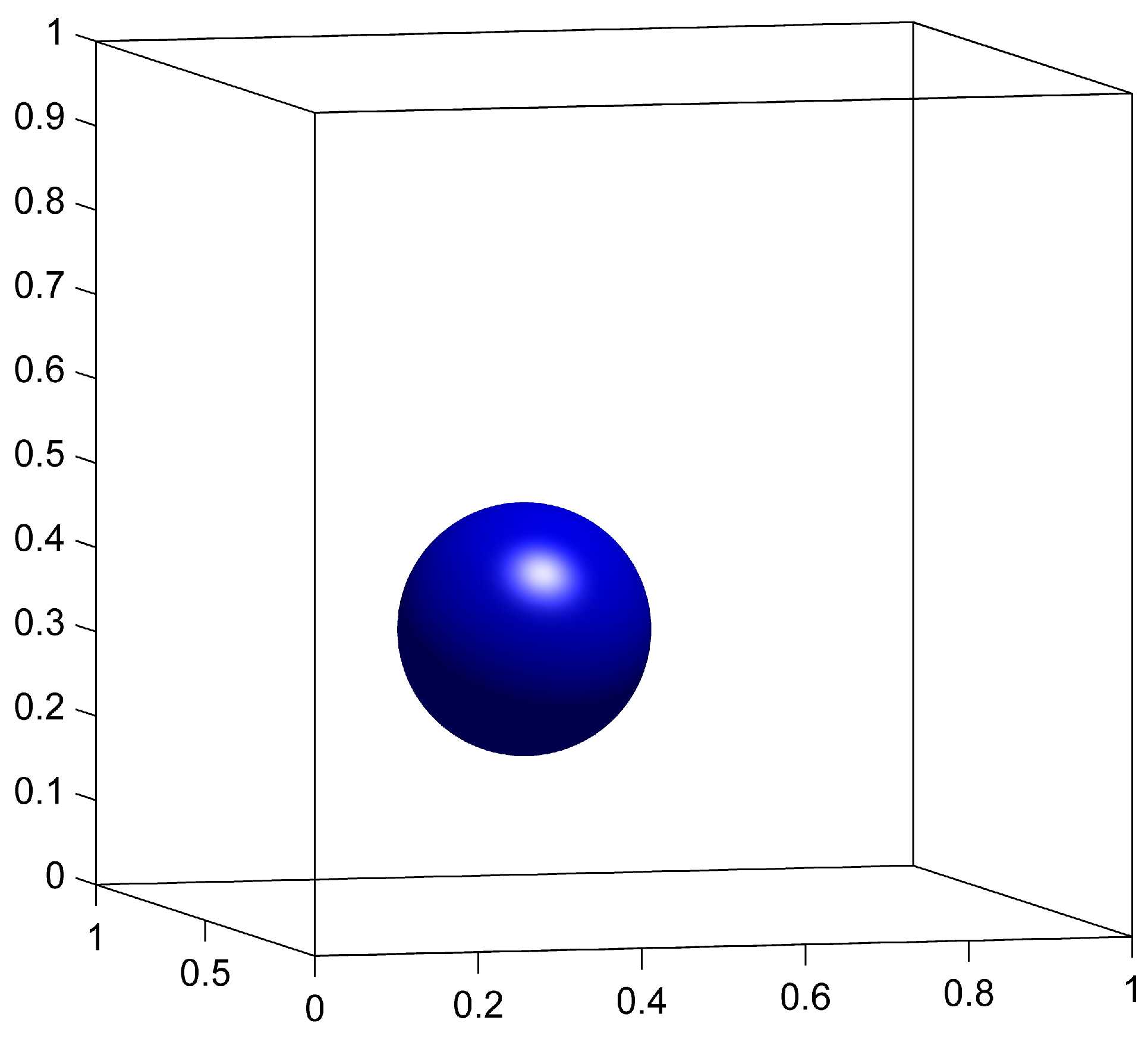}
			\caption{CM with remap, t=0}	\label{fig:3Dcr0}
		\end{subfigure}
		\begin{subfigure}{0.26\linewidth}
			\includegraphics[width=0.9\textwidth]{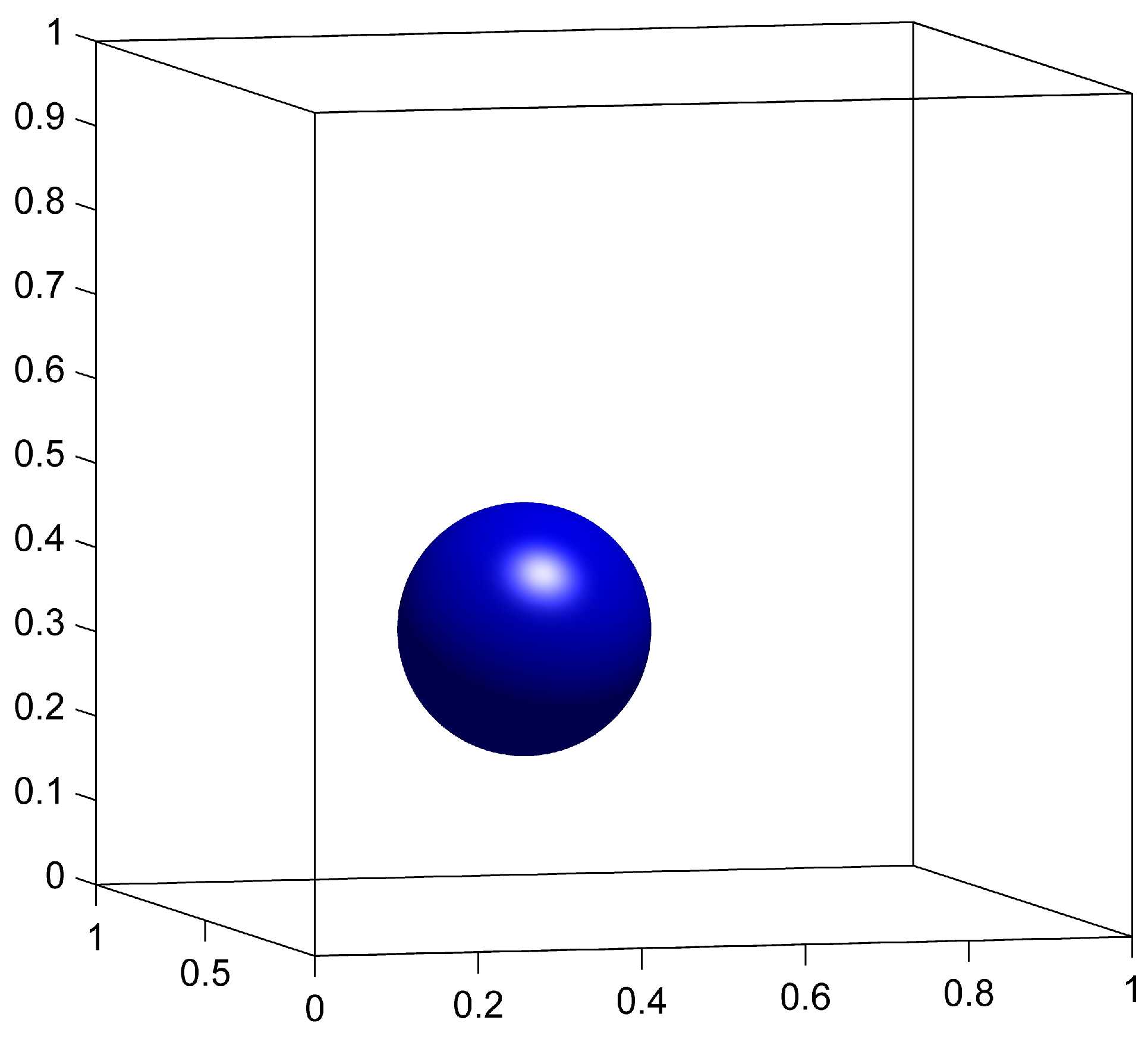}
			\caption{GALS, t=0}	 \label{fig:3Dga0}
		\end{subfigure}
		\begin{subfigure}{0.26\linewidth}
			\includegraphics[width=0.9\textwidth]{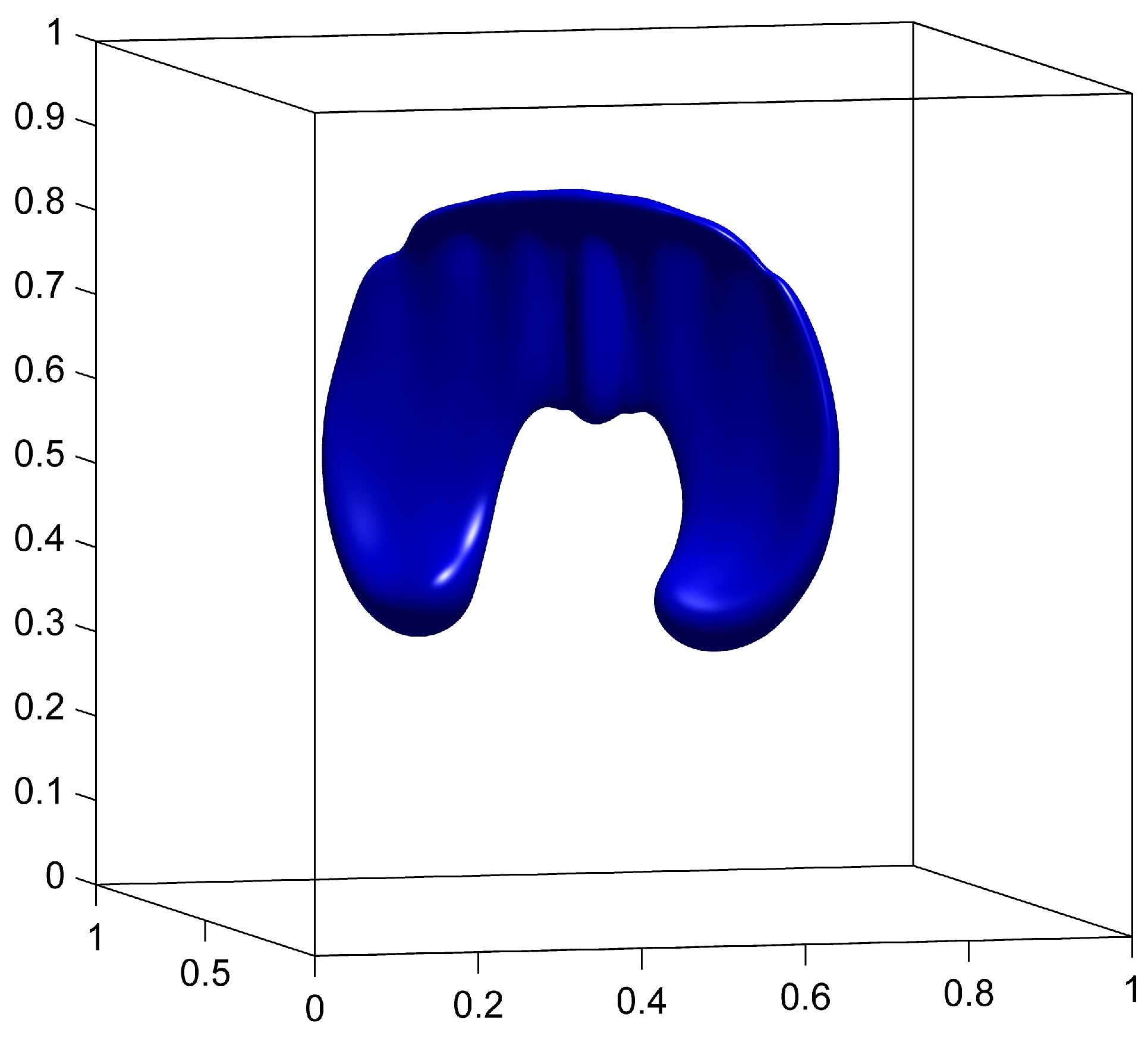}
			\caption{CM without remap, t=1}	\label{fig:3Dcn1}
		\end{subfigure}
		\begin{subfigure}{0.26\linewidth}
			\includegraphics[width=0.9\textwidth]{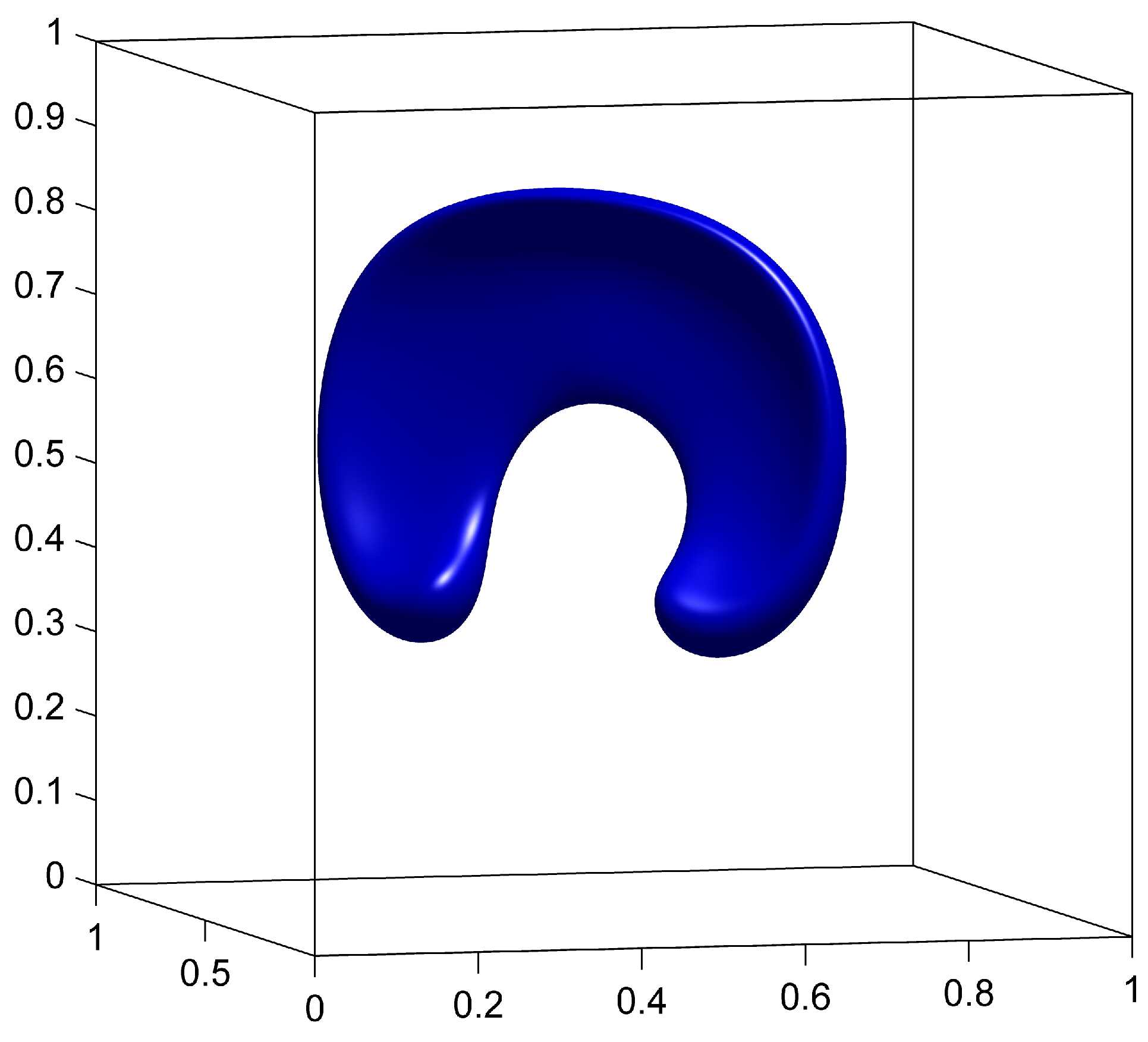}
			\caption{CM with remap, t=1}	\label{fig:3Dcr1}
		\end{subfigure}
		\begin{subfigure}{0.26\linewidth}
			\includegraphics[width=0.9\textwidth]{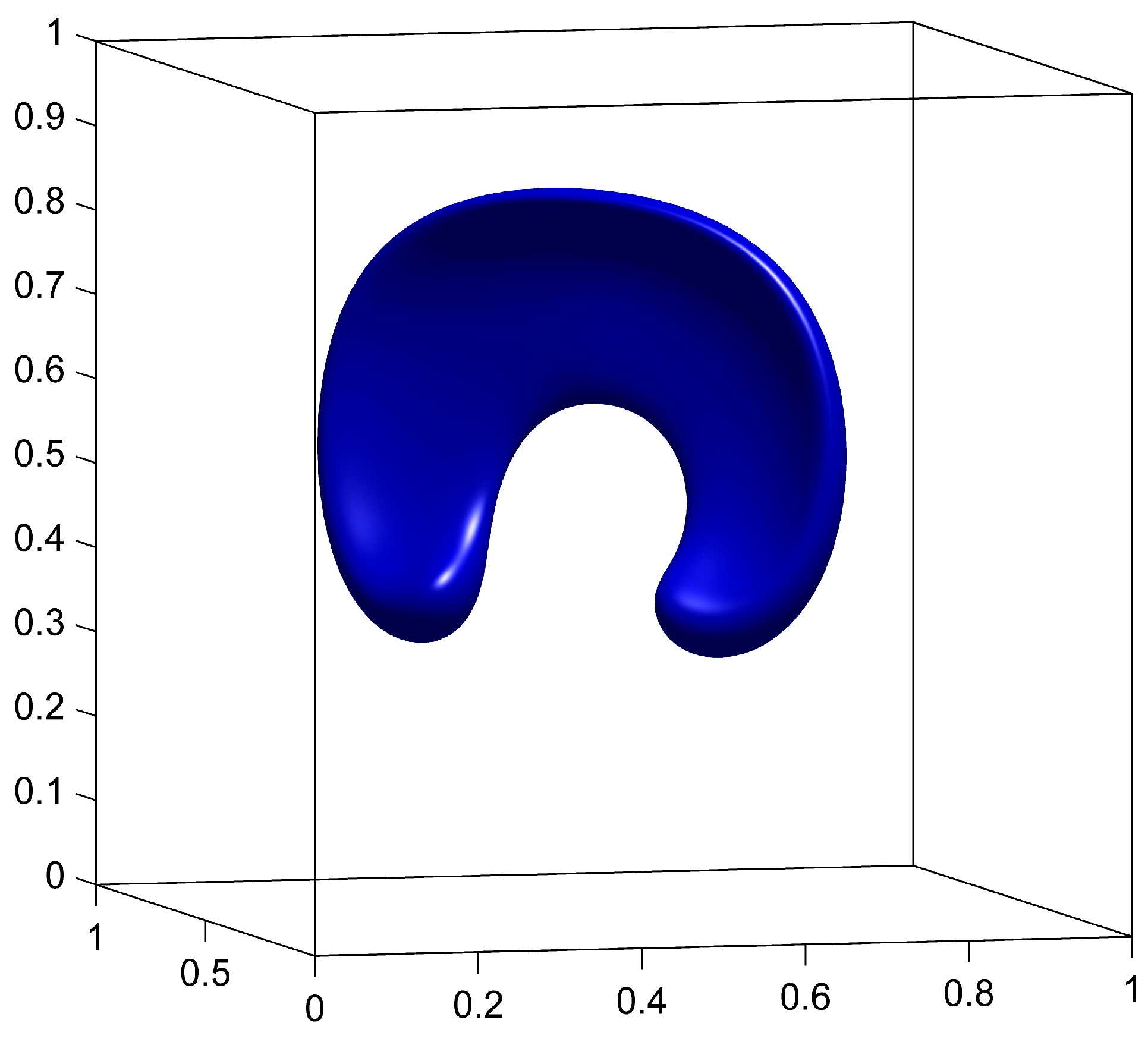}
			\caption{GALS, t=1}	\label{fig:3Dga1}
		\end{subfigure}
		\begin{subfigure}{0.26\linewidth}
			\includegraphics[width=0.9\textwidth]{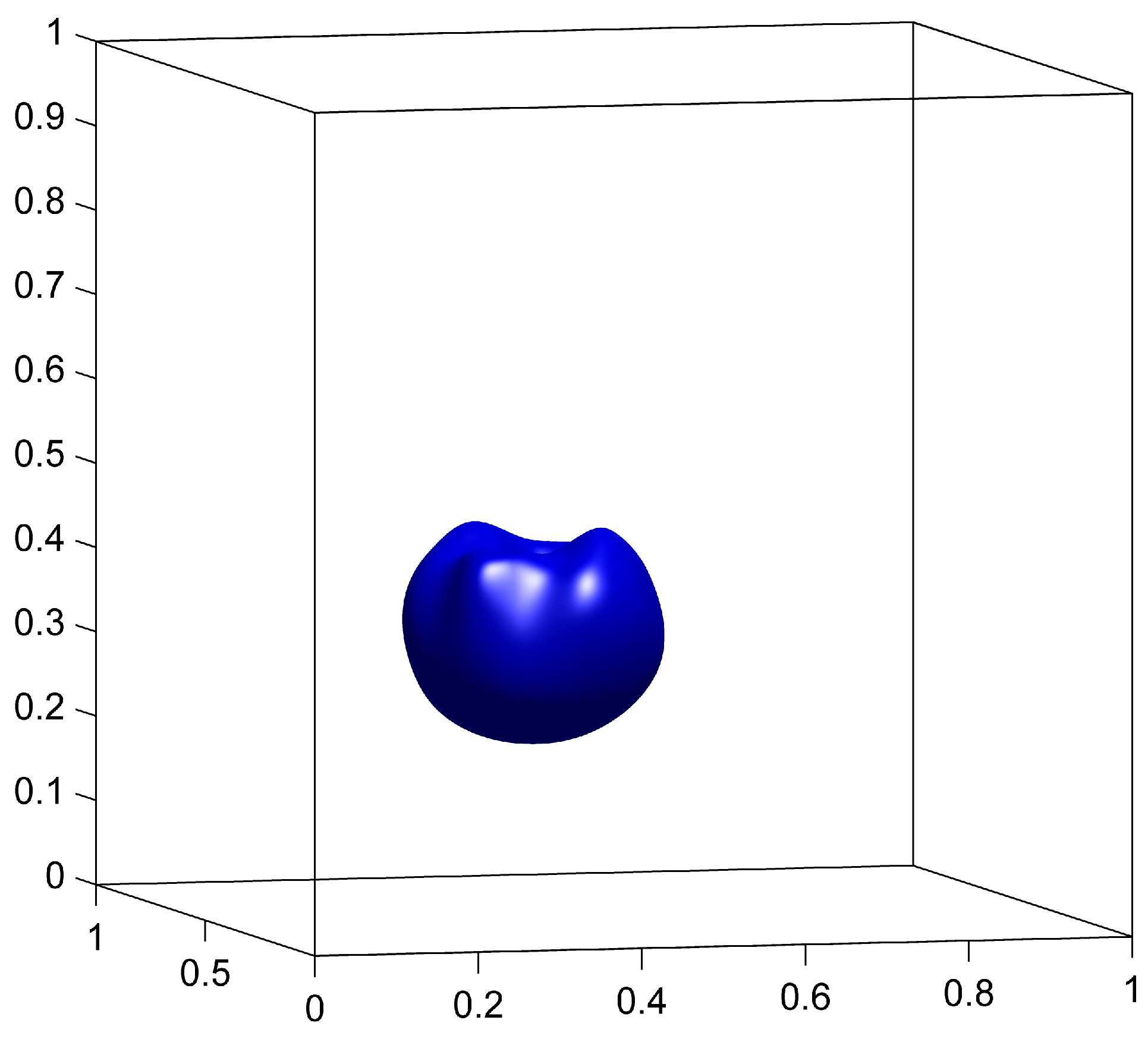}
			\caption{CM without remap, t=2}	\label{fig:3Dcn2}
		\end{subfigure}
		\begin{subfigure}{0.26\linewidth}
			\includegraphics[width=0.9\textwidth]{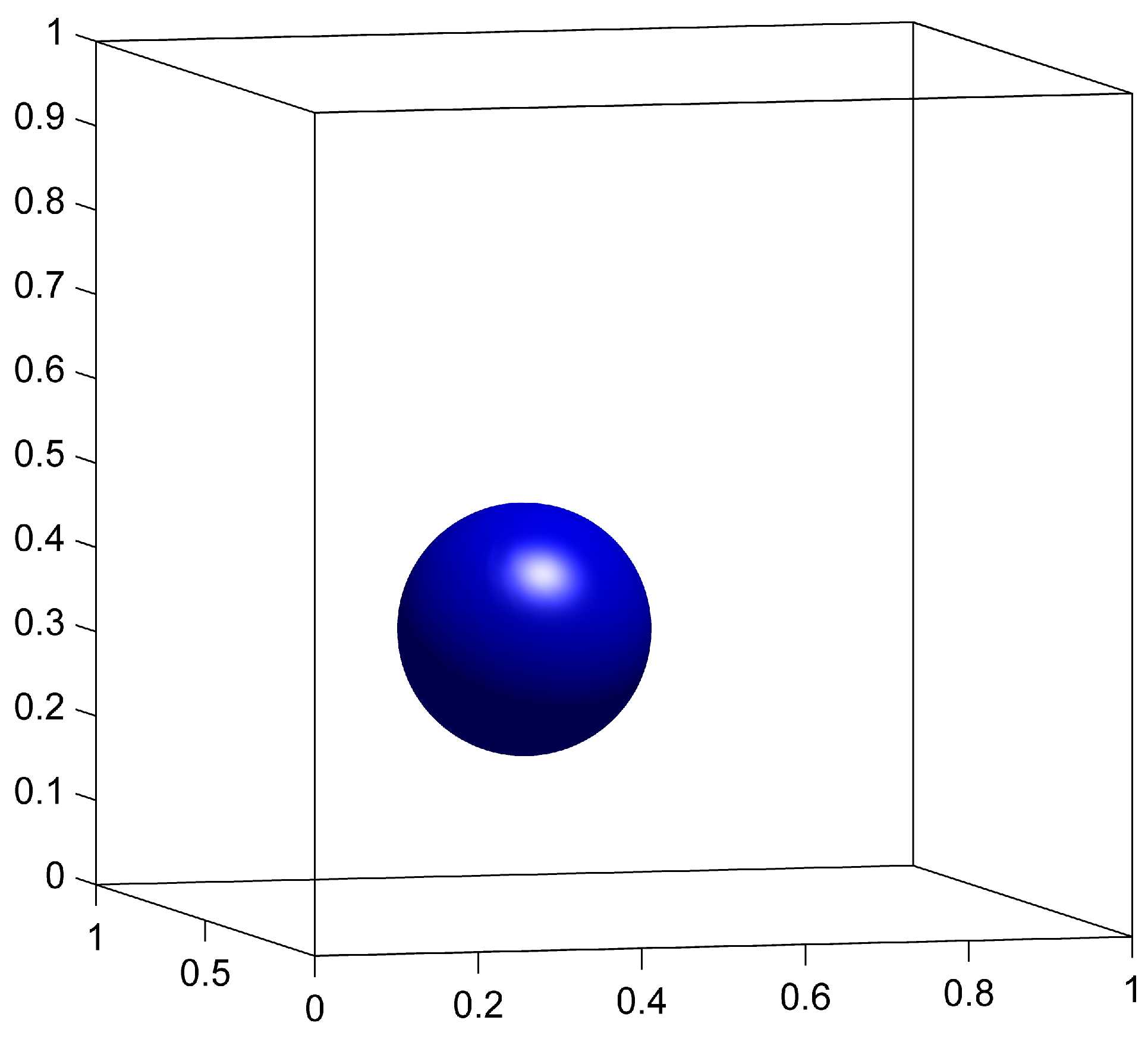}
			\caption{CM with remap, t=2}	\label{fig:3Dcr2}
		\end{subfigure}
		\begin{subfigure}{0.26\linewidth}
			\includegraphics[width=0.9\textwidth]{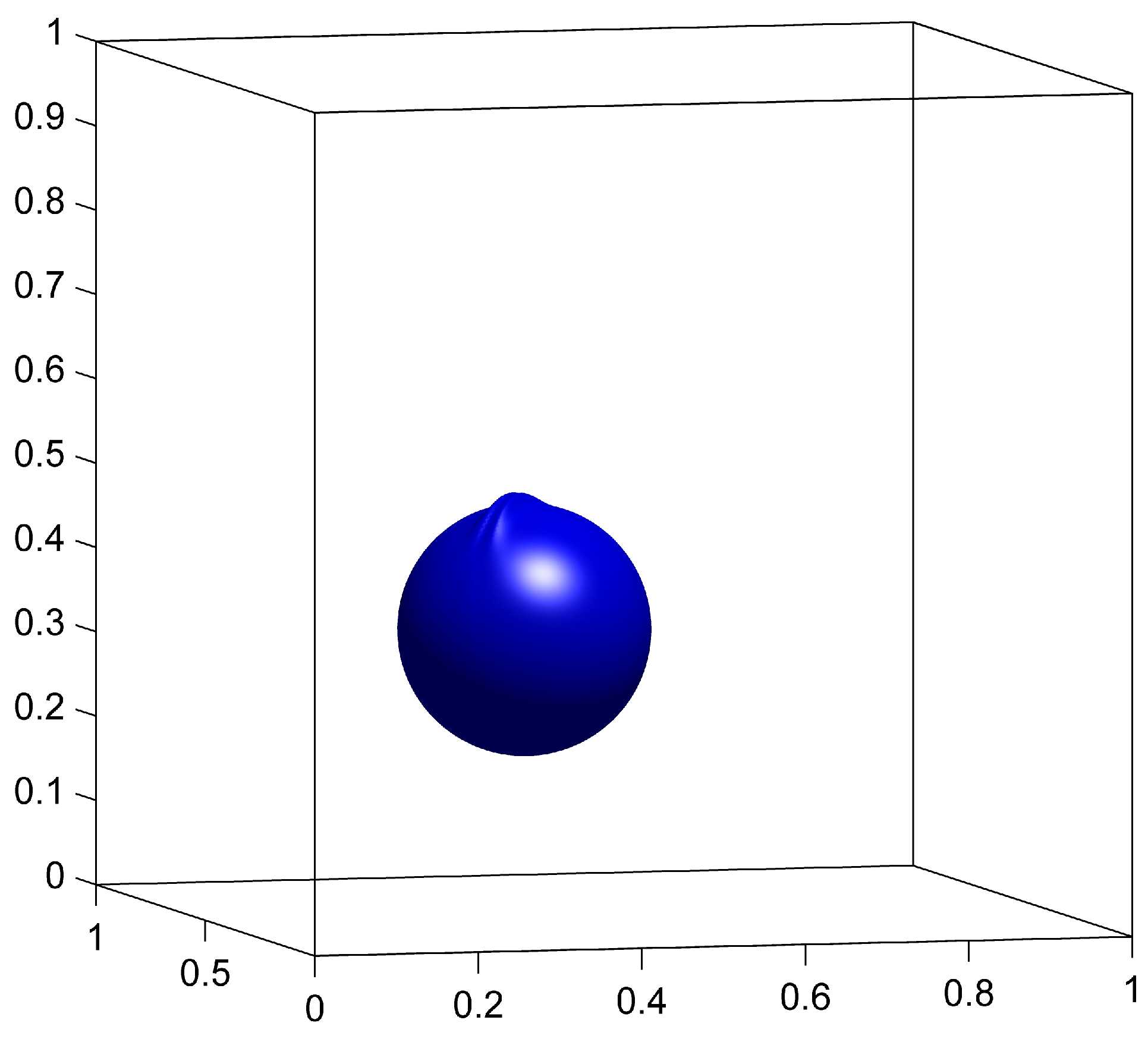}
			\caption{GALS, t=2}	\label{fig:3Dga2}
		\end{subfigure}
	\end{center}	
	\caption{3D deformation of a sphere. Left: CM on a $N_c = 16$ grid without remapping. Middle: CM on a $N_c = 16$ grid with remapping on a $N_f = 128$ grid. Right: GALS on a $N_g = 128$ grid.}
	\label{fig:3D}
\end{figure}

The advection of $\vec{\chi}$ is done on a $N_c = 16$ grid, and the remapping of $\vec{\chi}_0$ is done on a fixed $N_f = 128$ grid with a remapping tolerance $\mathcal{E}_1 = 10^{-4}$. We compare the results with and without remapping to demonstrate the benefits of the remapping step. We also compare the CM method to the GALS method computed on a $N_g = 128$ grid. Results for $t=0$, $t=1$ and $t=2$ for the three cases are shown in figure \ref{fig:3D}.

The surface at $t=2$ is expected to be identical to the surface at $t=0$ (sphere) due to symmetries in the velocity \eqref{eqGroup:3D_def}. This vector field causes the sphere to be stretched along the $y=1-x$ plane and creates very fine structures that cannot be represented on the coarse $16^3$ grid. We see from figure \ref{fig:3Dcn1} that without remapping, those structures are indeed not well represented and oscillations are observed on the scale of the coarse grid. These oscillations distort the surface for all subsequent times, and the final shape in figure \ref{fig:3Dcn2} is significantly different from the initial sphere.

When using the remapping on the fine grid, the fine structures caused by the deformation field can be accurately represented as a result of storing $\vec{\chi}_0$ on a $N_f = 128$ grid. In figure \ref{fig:3Dcr1}, the width of the stretched surface is of the order of the fine grid's cell width, therefore causing no major representation issues. At $t=2$ (figure \ref{fig:3Dcr2}), the surface is visually identical to the initial sphere. 

The GALS method also uses a $N_g = 128$ grid and therefore uses comparable computational resources for spacial resolution as CM with remapping. The main differences between the two methods is that CM method is initialized with a trivial identity map and does not discretize the more complicated initial condition of the level-set function. The effects of this can be seen later in the simulation in figures \ref{fig:3Dcr2} and \ref{fig:3Dga2}, where GALS accumulates more error than CM.

In addition to testing the accuracy of the CM method with remapping, this experiment also demonstrates the computational efficiency of the method. Table \ref{table:time3D} compares the time taken to compute the three tests of figure \ref{fig:3D}. As expected, the time taken by the CM method without remapping is small, but the results are not precise. The time taken by the CM method with remapping is approximately 5 times larger, but achieves the best result out of the three cases. The GALS method takes about 20 times longer than the CM method with remapping and fails to provide better results. This indicates that the separation of scales using the semi-group property of the characteristic maps provides significant improvements to computational efficiency.

\begin{table}[h]
\centering
\begin{tabular}{|c|c|}
\hline
Method & Computational Time (sec.) \\ 
\hline
\hline
CM without remapping & 267\\
\hline
CM with remapping & 1430\\ 
\hline
GALS & 29687\\
\hline
\end{tabular}
\caption{Time comparison (in seconds) of the GALS and CM methods with and without remapping for the 3D deformation test of figure \ref{fig:3D}. We use $N_c = 16$ and $N_f = N_g = 128$.}
\label{table:time3D}
\end{table}

\subsection{Complicated Sets} \label{sec:complex}
This section examines the application of the CM method for the advection of arbitrary sets. We take the deformation field $\vec{v}(\vec{x}, t)$ given by
\begin{gather} \label{eqGroup:mandelVec}
\vec{v}(x, y, t) = \left(  \begin{matrix}
\cos\left(\frac{\pi t}{A}\right) \left(  -  \frac{1}{4}L\left(x, y\right) \left( \sin\left(\frac{3}{2} \pi x\right) \right)^2 \sin\left(4 \pi y\right)  +  \frac{3}{4}R\left(x, y\right)\left(x- \frac12 \right)    \right) \\
\cos\left(\frac{\pi t}{A}\right) \left(  \frac{1}{4}L\left(x, y\right) \left( \sin\left(2 \pi y\right) \right)^2\sin\left(3 \pi x\right)  +  \frac{3}{4}R\left(x, y\right)\left(y- \frac12 \right)   \right)
\end{matrix}  \right)
\end{gather}
with A=16 and where $L$ and $R$ are smooth weight functions defined by
\begin{gather}
R\left(x, y \right) = \sin\left(\pi\left(4x^2 - 5x^3 + 2x^4\right)\right)\sin\left(\pi y\right)\\
L\left(x, y \right) = \sin\left(\pi \left(1-x\right)^3\right)\sin\left(\pi y\right)\sin\left(\frac{3}{2} \pi x\right) \left( \sin\left(2 \pi y\right) \right)^2
\end{gather}

\begin{figure}[h]    
	\begin{center}
		\begin{subfigure}{0.2\linewidth}
			\includegraphics[width=\textwidth]{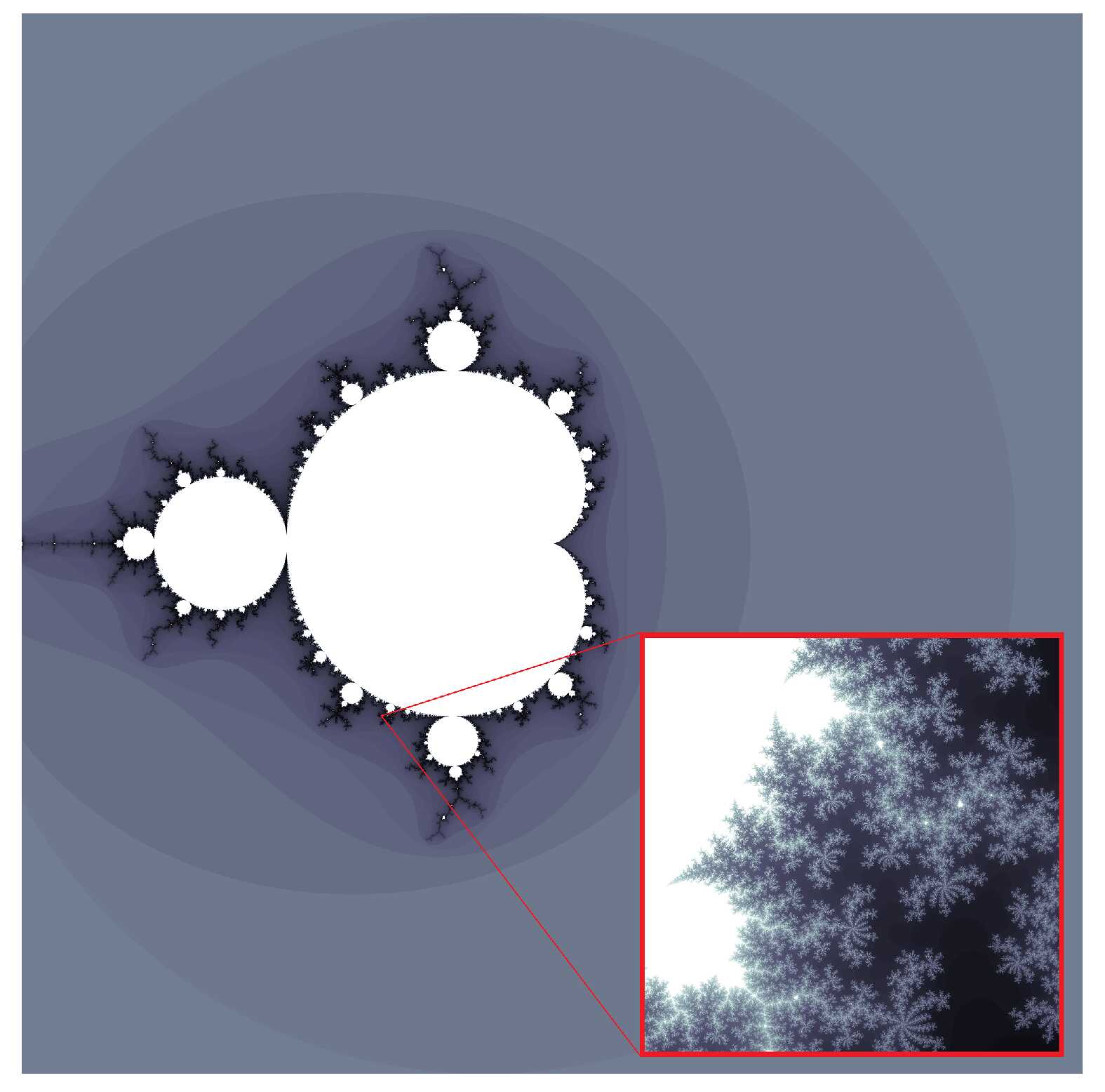}
			\caption{t=0}	\label{fig:mandelbrotstart}
		\end{subfigure}
		\begin{subfigure}{0.2\linewidth}
			\includegraphics[width=\textwidth]{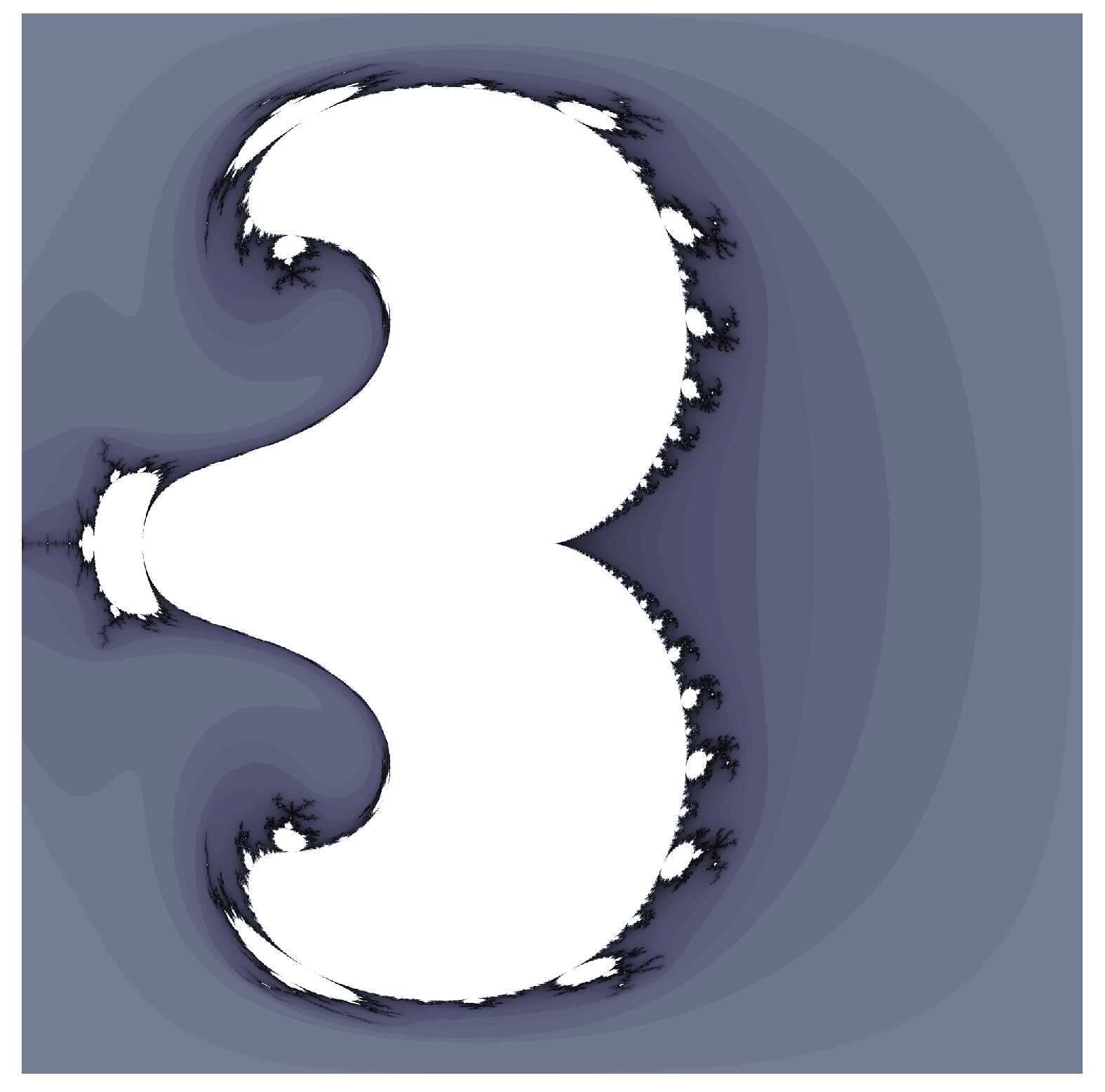}
			\caption{t=2}	
		\end{subfigure}
		\begin{subfigure}{0.2\linewidth}
			\includegraphics[width=\textwidth]{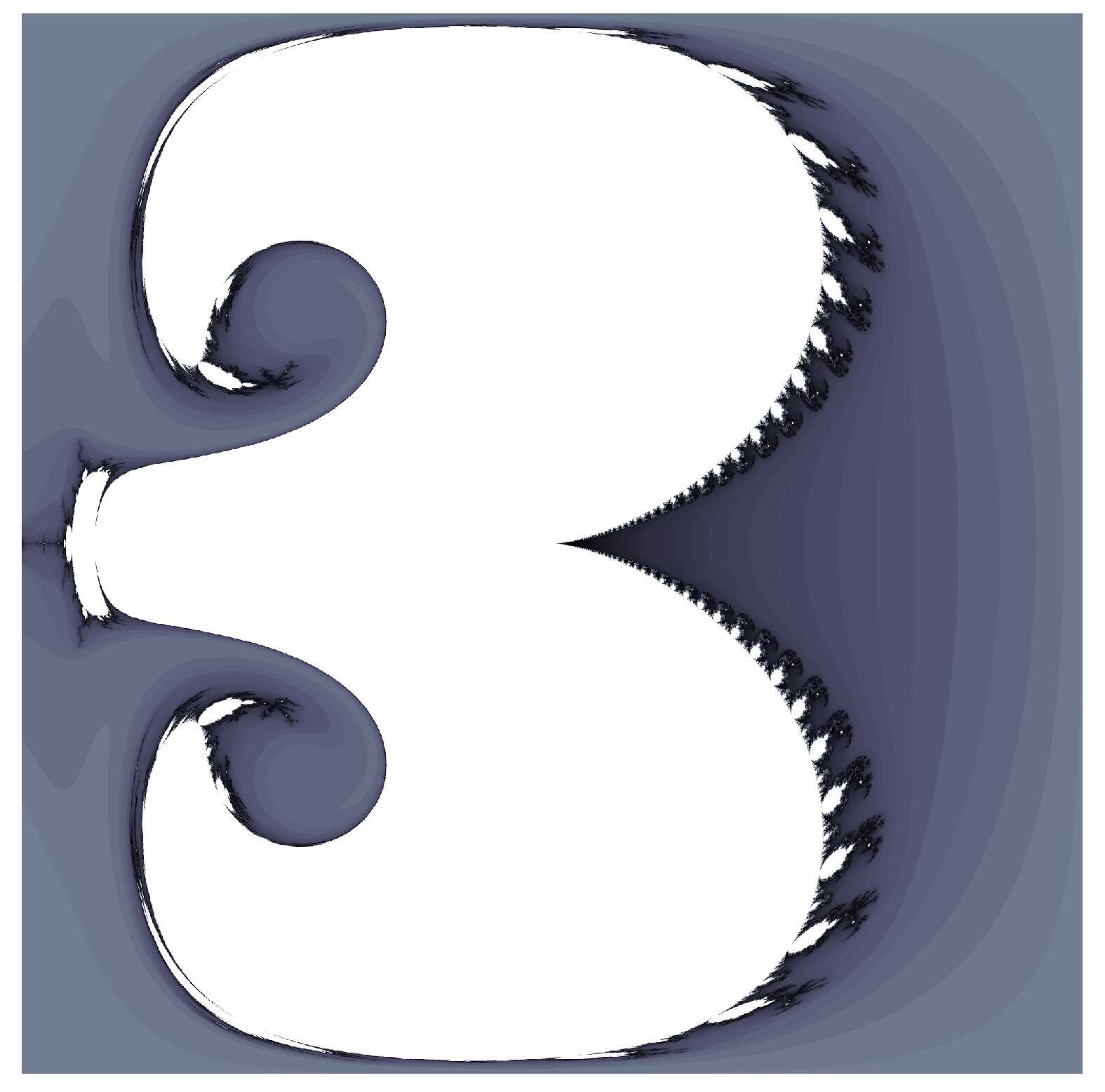}
			\caption{t=4}	
		\end{subfigure}\\
		\begin{subfigure}{0.2\linewidth}
			\includegraphics[width=\textwidth]{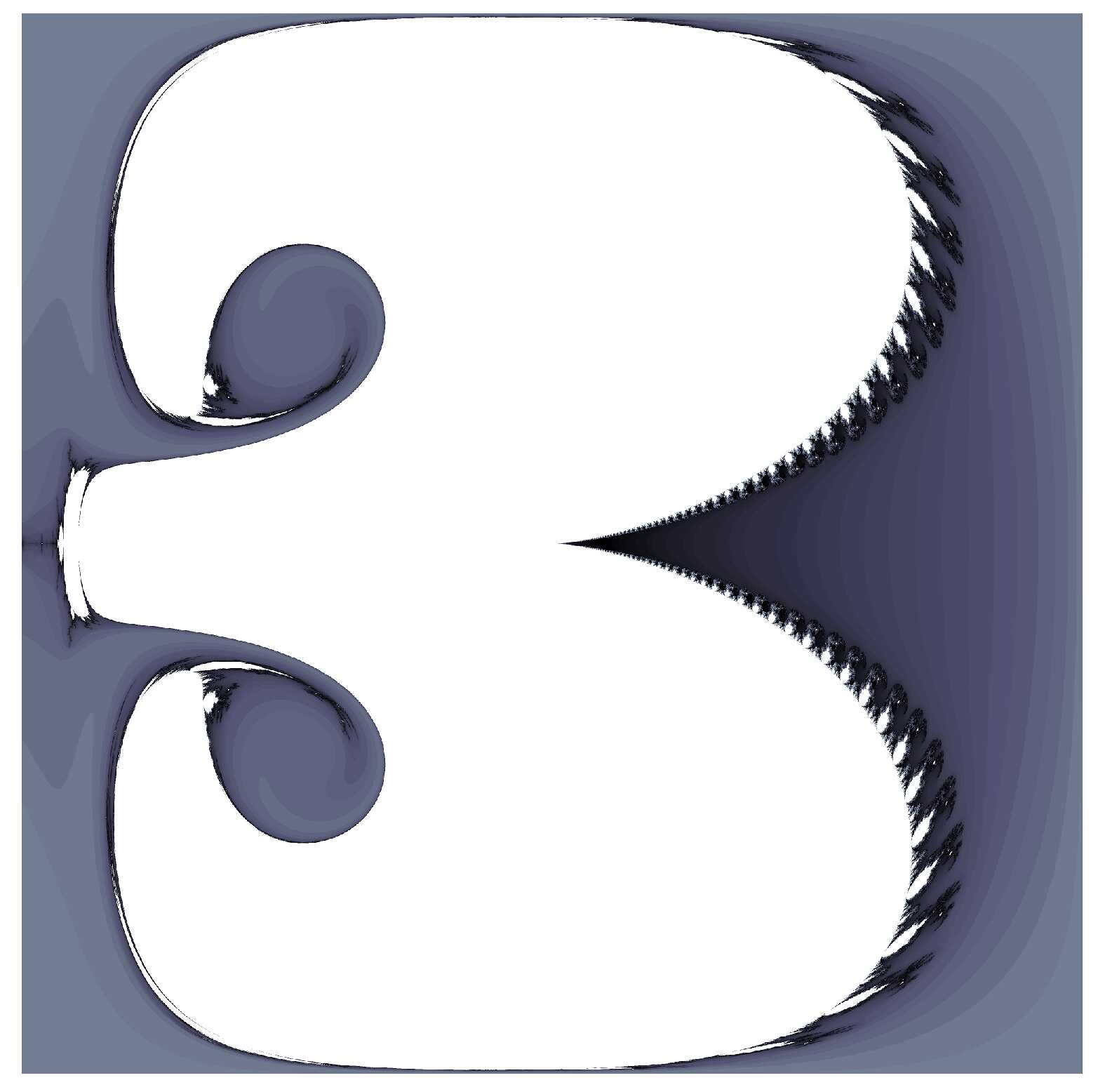}
			\caption{t=6}	
		\end{subfigure}
		\begin{subfigure}{0.2\linewidth}
			\includegraphics[width=\textwidth]{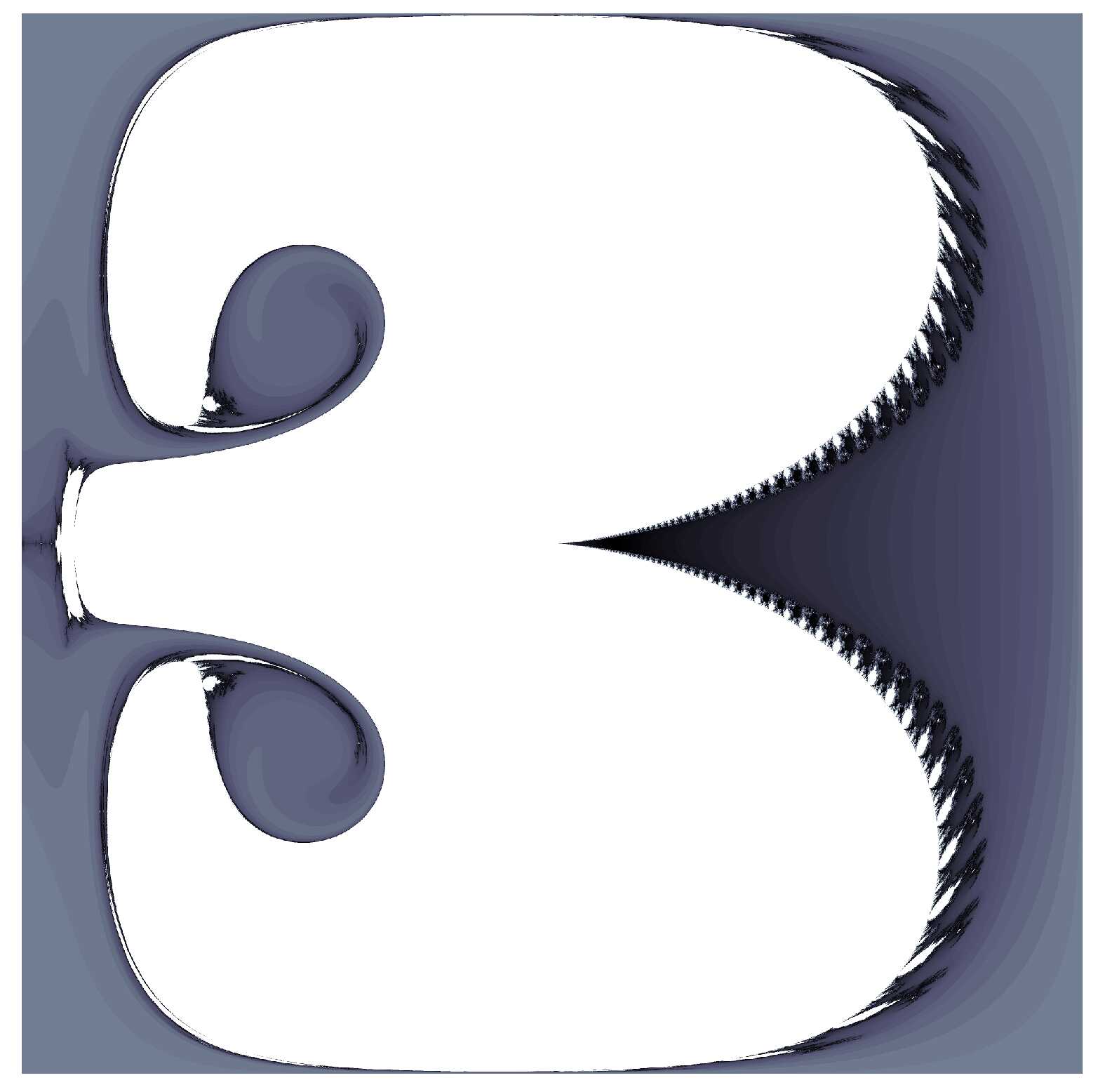}
			\caption{t=8}	
		\end{subfigure}
		\begin{subfigure}{0.2\linewidth}
			\includegraphics[width=\textwidth]{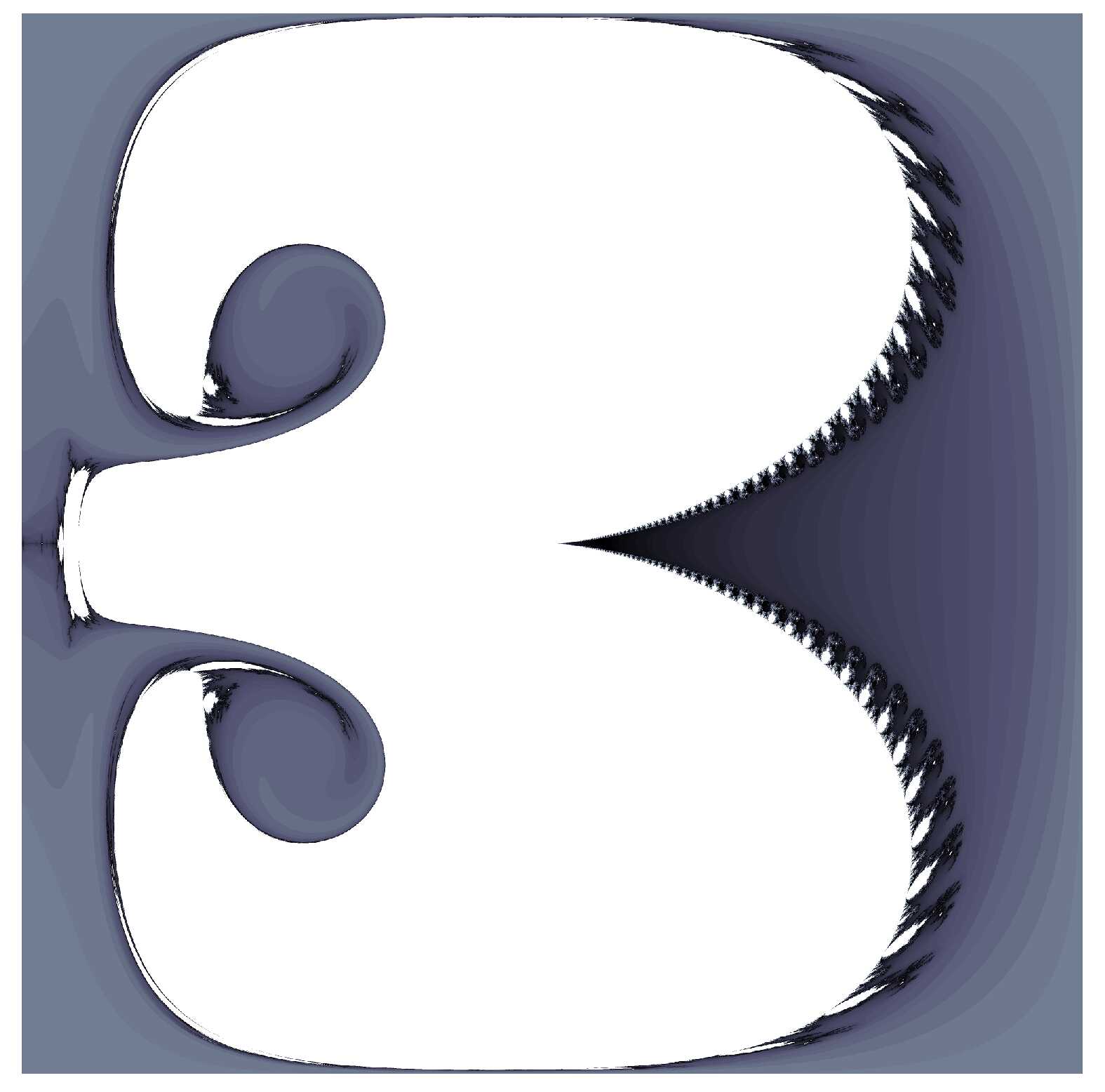}
			\caption{t=10}	
		\end{subfigure}\\
		\begin{subfigure}{0.2\linewidth}
			\includegraphics[width=\textwidth]{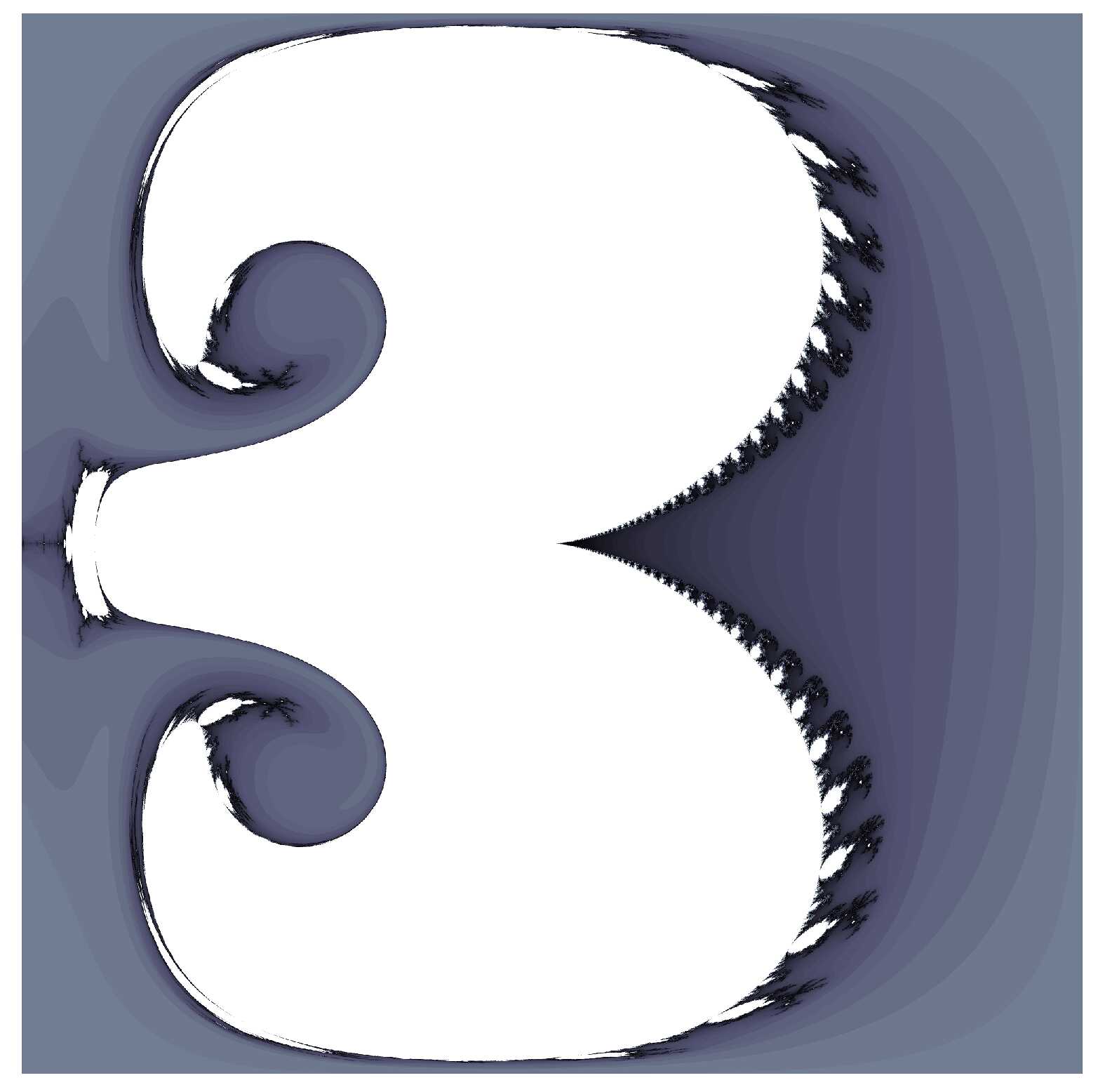}
			\caption{t=12}	
		\end{subfigure}
		\begin{subfigure}{0.2\linewidth}
			\includegraphics[width=\textwidth]{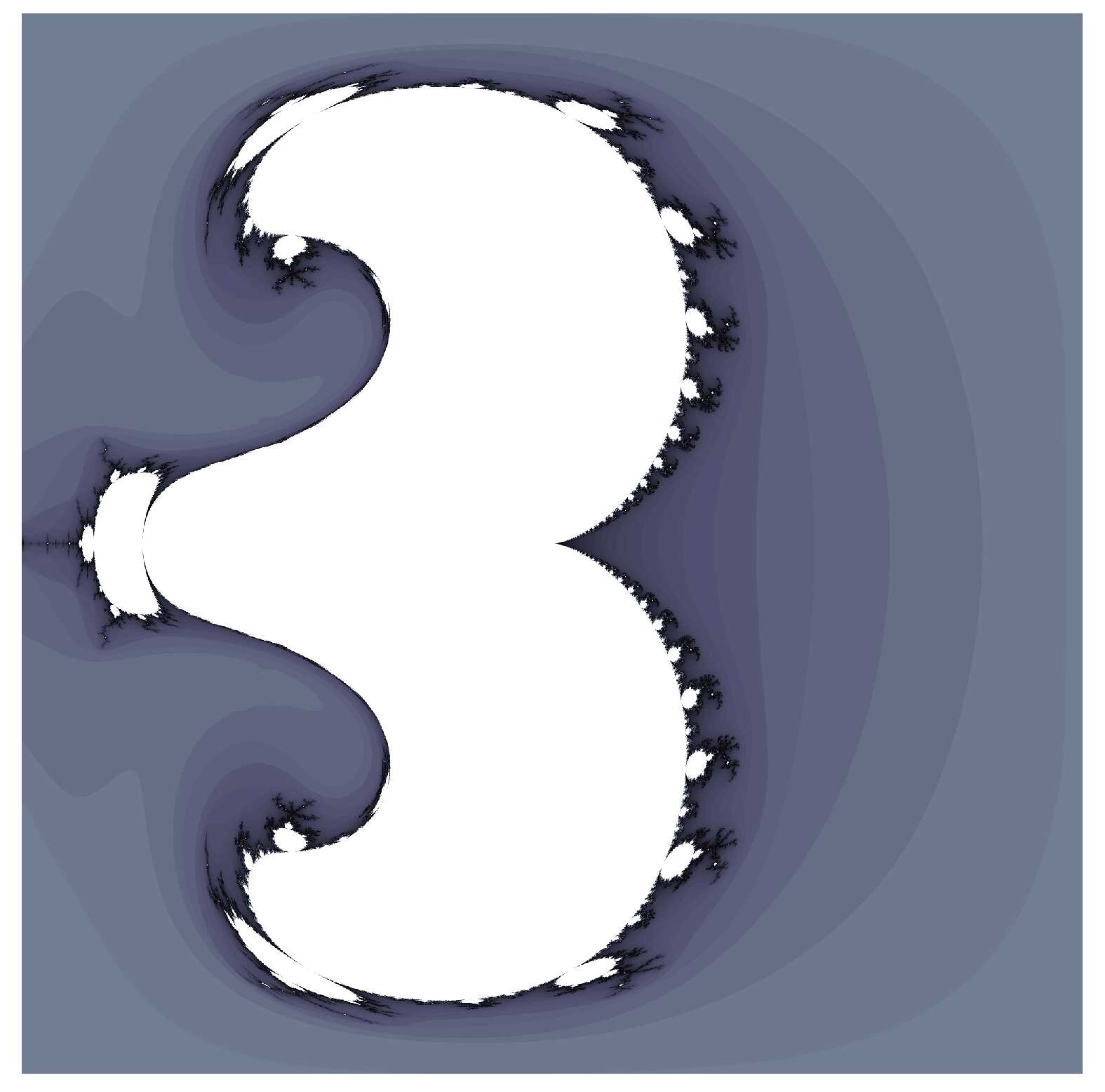}
			\caption{t=14}	
		\end{subfigure}
		\begin{subfigure}{0.2\linewidth}
			\includegraphics[width=\textwidth]{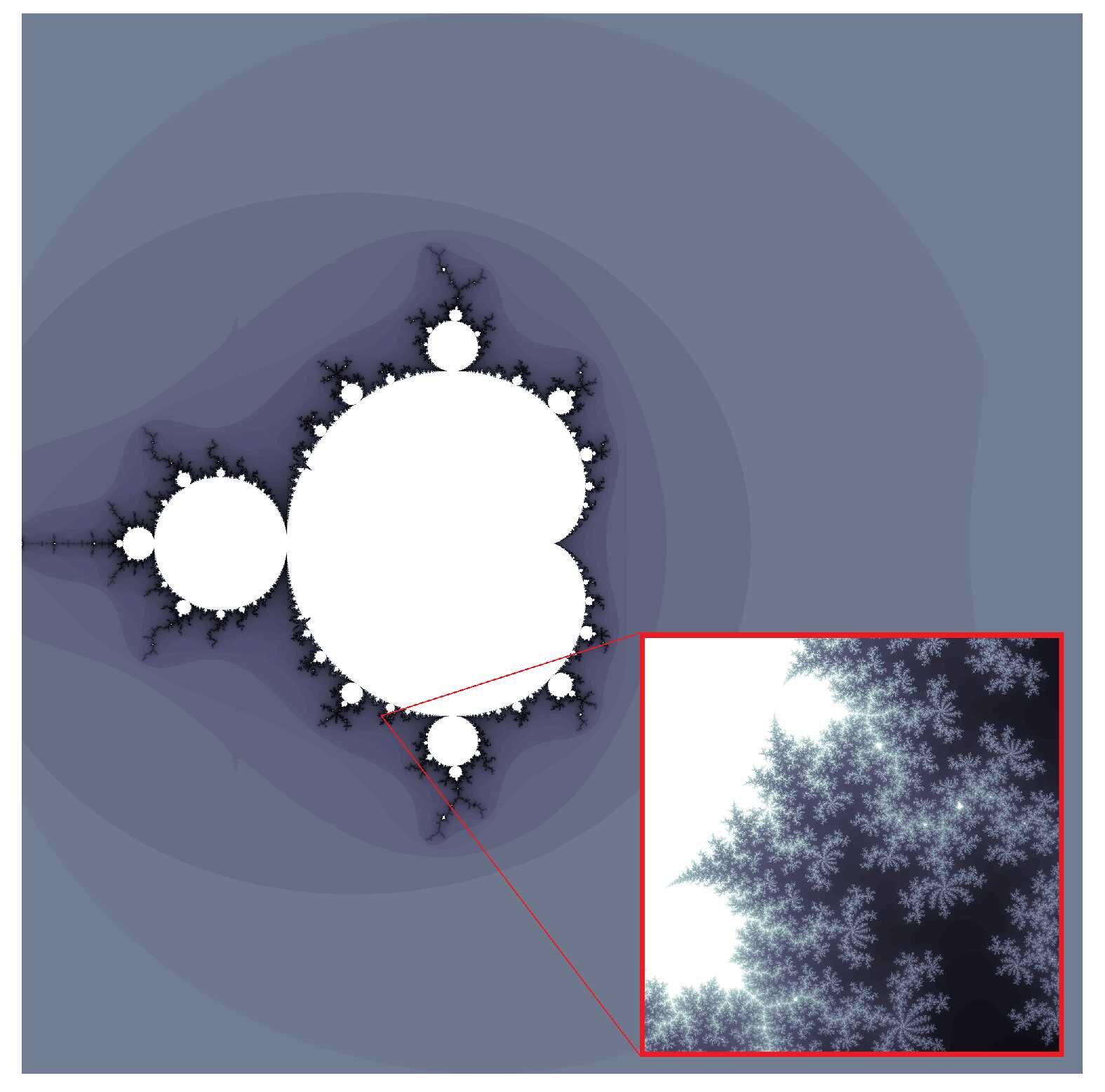}
			\caption{t=16}	\label{fig:mandelbrotend}
		\end{subfigure}
	\end{center}	
	\caption{2D deformation of the Mandelbrot set computed with the CM method with $N_c = 32$ and fixed $N_f = 1024$. The zoomed regions in figures $(a)$ and $(i)$ correspond to a single grid cell of the remapping grid.}
	\label{fig:mandelbrot}
\end{figure}

This vector field causes the left region to swirl into two vortices and the right region to expand around the center of the domain. We apply this vector field to the Mandelbrot set and compute the advection of the set with the CM method. We use $N_c = 32$ for the coarse grid, a fixed $N_f = 1024$ for the fine grid and a remapping tolerance $\mathcal{E}_1 = 10^{-7}$. The results at times $t=0$ to $t=16$ are shown in figure \ref{fig:mandelbrot}. 

Traditional advection schemes would have difficulties dealing with complex geometries such as the Mandelbrot set. With the CM method, the complexity of the initial set does not translate to difficulties in numerical simulation. Since the advected quantity is the pullback of the initial function by the deformation map and since we know the initial set to arbitrary precision, there is no difficulty capturing all the fine structures. This is a property of the interpolation structure of the method; since the characteristic map is evolved as a function defined everywhere in the domain, $\phi(\vec{x}, t) = \phi_0(\vec{\chi}(\vec{x}, t))$ can be evaluated directly for any $\vec{x}$. Therefore, we achieve arbitrary subgrid resolution and preserve the information contained in the initial condition $\phi_0$ throughout the simulation. We observe indeed in figures \ref{fig:mandelbrotstart} and \ref{fig:mandelbrotend} that the same resolution is available at the beginning and at the end of the simulation. In both figures, we show a zoomed view of a single cell in the fine grid to show that at the end we recuperate, with arbitrary subgrid resolution, the same set as the initial condition, as expected. This is a remarkable and unique property of the CM framework and has important implications in real-life applications. As illustrated by the above example, the spacial resolution of the solution to the advection equation is not affected by the resolution of the characteristic map. Indeed, the Mandelbrot set requires infinite resolution to represent whereas the map is computed on a $N_f = 1024$ grid; $N_f$ is selected to provide a good approximation to the characteristic map $\vec{X}$, it only needs to resolve the dynamics of the flow. In practice, one might have a very finely sampled discretized initial condition. Traditional methods such as GALS need to numerically represent the initial condition accurately in order to maintain resolution throughout the simulation. If the dynamics of the flow contain less spacial features, this would be a waste. The CM method is more efficient in that the characteristic map uses only the necessary computational resources to represent the deformation. The initial condition, no matter how finely resolved, will be pulled-back by the deformation map to generate solutions of the sample resolution as the initial condition.

\subsection{Triple and Quadruple Points Mosaic}
For a given velocity field $\vec{v}$, the corresponding characteristic map $\vec{\chi}$ is a solution operator which maps any initial condition to the solution of the advection under $\vec{v}$. This implies that multiple initial conditions can be advected simultaneously using a single map $\vec{\chi}$. Furthermore, these solutions remain mutually consistent throughout the entire simulation because they arise from the same deformation map. We present two tests showing this feature.

\begin{figure}[h]
	\begin{center}
		\begin{subfigure}{0.2\linewidth}
			\includegraphics[width=\textwidth]{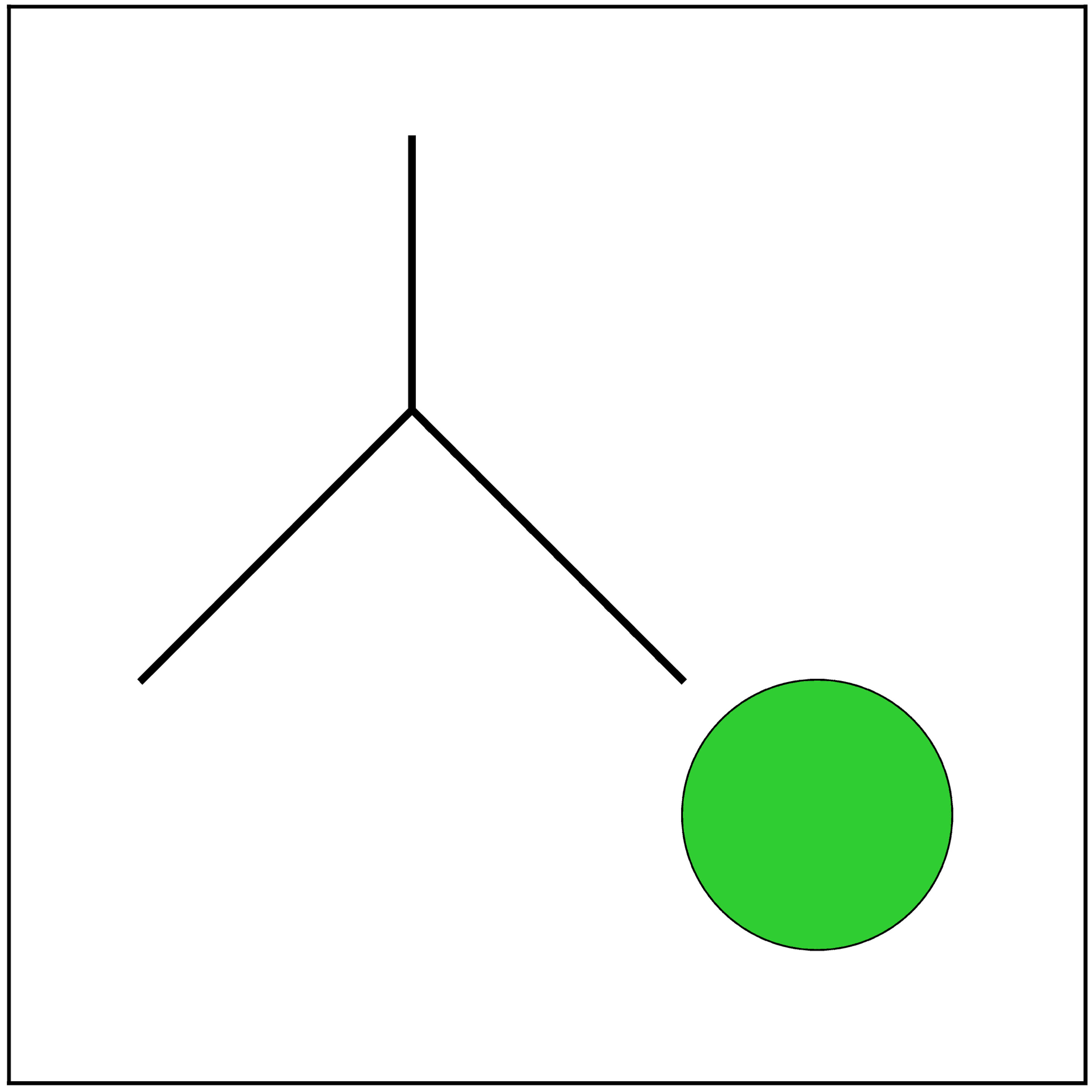}
			\caption{t=0}	
		\end{subfigure}
		\begin{subfigure}{0.2\linewidth}
			\includegraphics[width=\textwidth]{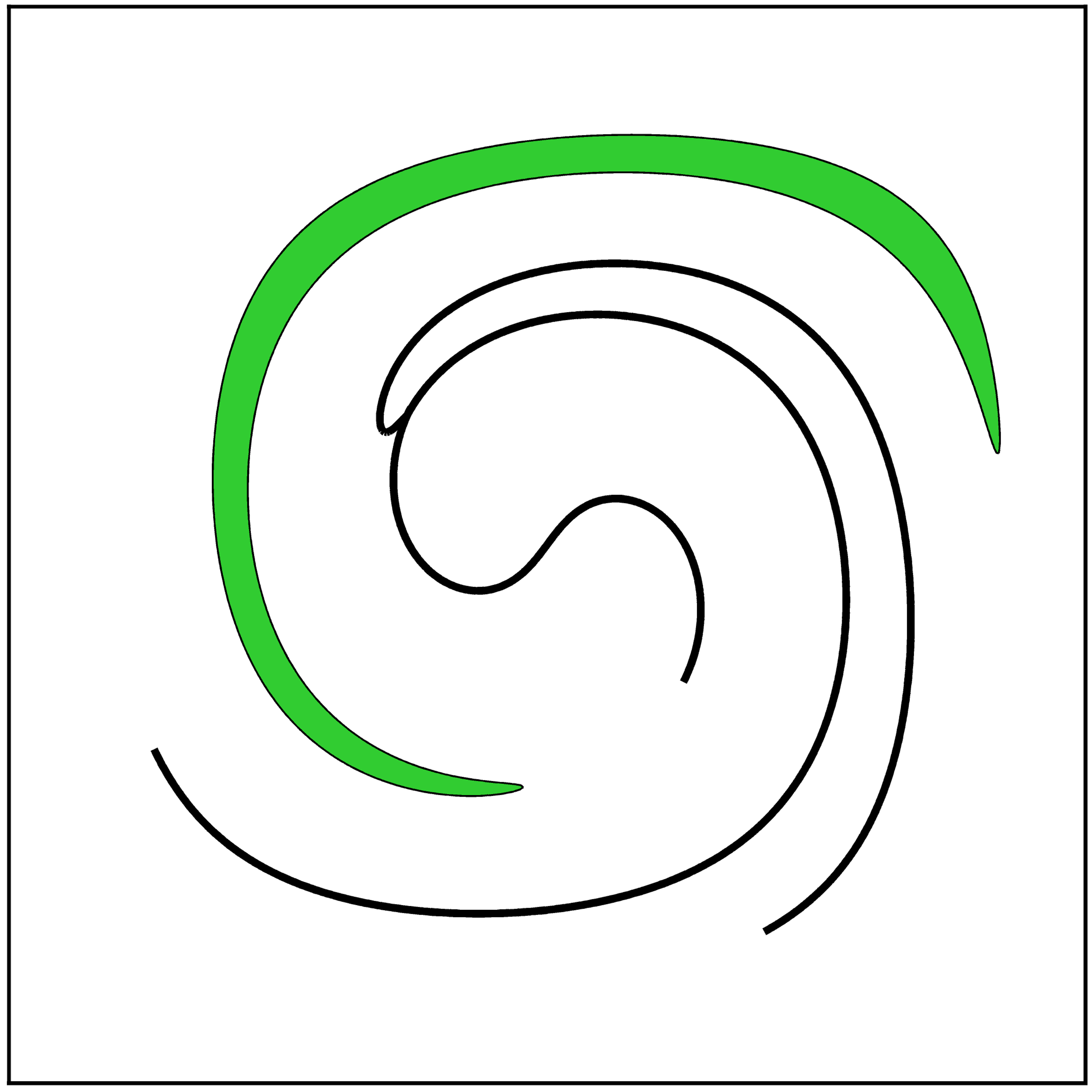}
			\caption{t=2}	
		\end{subfigure}
		\begin{subfigure}{0.2\linewidth}
			\includegraphics[width=\textwidth]{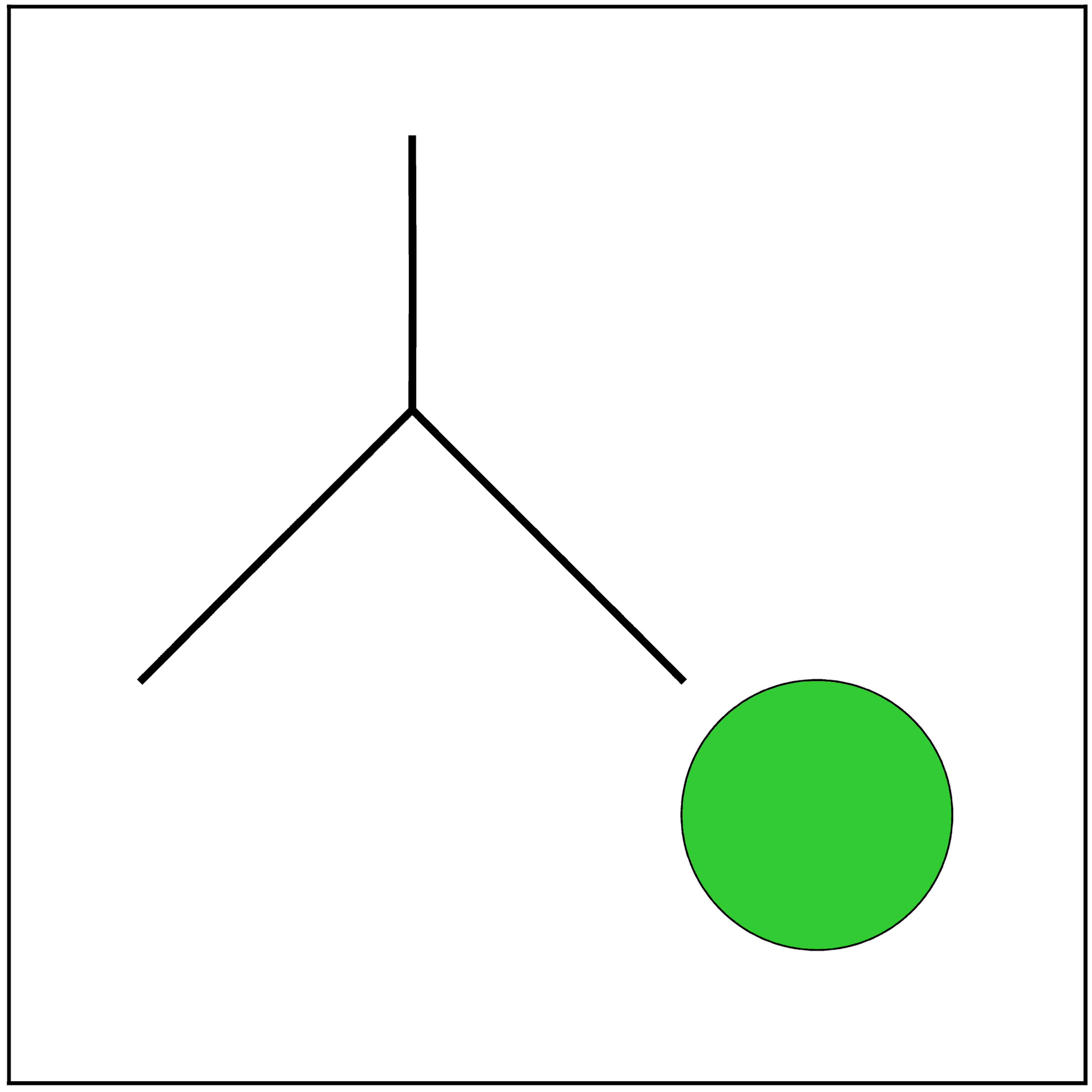}
			\caption{t=4}	
		\end{subfigure}
	\end{center}	
	\caption{2D deformation of open curves. Each of the three branches (black) is defined by a linear level-set function enclosed in a mask set, also defined by another level-set function. The circle is also represented by a level-set function. }
	\label{fig:openCurve}
\end{figure}

The first test involves open curves. These are represented implicitly using two functions: a level-set function and a mask function. For instance, a segment of a straight line can be represented using an affine function whose zero level-set is the line containing the segment, and a mask function which represents a region whose intersection with the line gives the desired segment. Figure \ref{fig:openCurve} shows an example where we advect a circle along with three independent line segments forming a triple point. These objects are transported in the swirling velocity field \eqref{eqGroup:swirl_def} with $A = 4$. Normally, these seven functions used to define the curves (three lines, three masks and one circle) need to be advected separately; one would need to solve seven advection equations independently, giving potentially mutually incoherent results. The CM method however, decouples the advection from the initial set, we can advect those seven functions at the same time using the same characteristic map. When plotting is required, the functions are evaluated at $\vec{\chi}(\vec{x},t)$, incurring only the extra cost of seven function evaluations and an evaluation of the map on the plotting grid. Since the solutions are obtained from pullback by the same characteristic map, the coherence of the results is guaranteed (e.g. triple points are guaranteed to stay triple points).

The second test shows an example of the transport of a Mosaic pattern. We take $[0,1]^2$ as a periodic domain and subdivide it into multiple regions. This kind of initial condition arises for instance in the simulation of multiphase flows. The difficulty of these simulations resides in the multiple intersection points caused by the junction of more than two different fluids. These points are hard to represent using a single level-set function, but can be represented piecewise as multiple level-sets. Saye and Sethian \cite{saye2012analysis,saye2011voronoi} propose an approach to multiphase flows using level-sets, while Da et al. \cite{da2013multimaterial} suggest to use multimaterial front tracking. In the context of the CM method, this is not a problem since we transform the whole domain and can deal with functions of any level of complexity afterwards. The different regions can simply be represented by a single function by multiplying a different constant to their respective indicator functions.

The velocity field used for this test is
\begin{gather}
\vec{v}(x, y, t) = \left( \begin{matrix}
\cos\left(\frac{\pi t}{2}\right)\cos\left(2y\pi + 2\sin\left( 2  \left( \cos \left(\frac{\pi t}{2}\right) \right)^2  \right)\right)\\
\cos\left(\frac{\pi t}{2}\right)\sin\left(2x\pi + 2\sin\left(  \left( \cos\left(\frac{\pi t}{2}\right) \right)^2   \right)\right)
\end{matrix}  \right)
\end{gather}
This velocity generates a periodic flow map with period 2. We used $N_c = 32$ and $N_f = 512$. It took $65$ seconds to compute the $2048$ steps of this simulation on a single $3.0$GHz CPU, taking $0.008$ seconds for a regular step and $0.45$ seconds when remapping, which was necessary on average approximately every $12$ step. We also highlight three regions in the domain using dashed circles. These regions represent two triple points and one quadruple point. The dashed circles are transported in the flow as passive particles with the initial intersection points as starting position. This shows that the intersections obtained from CM do follow the right path. 

\begin{figure}[h]
	\begin{center}
		\begin{subfigure}{0.2\linewidth}
			\includegraphics[width=\textwidth]{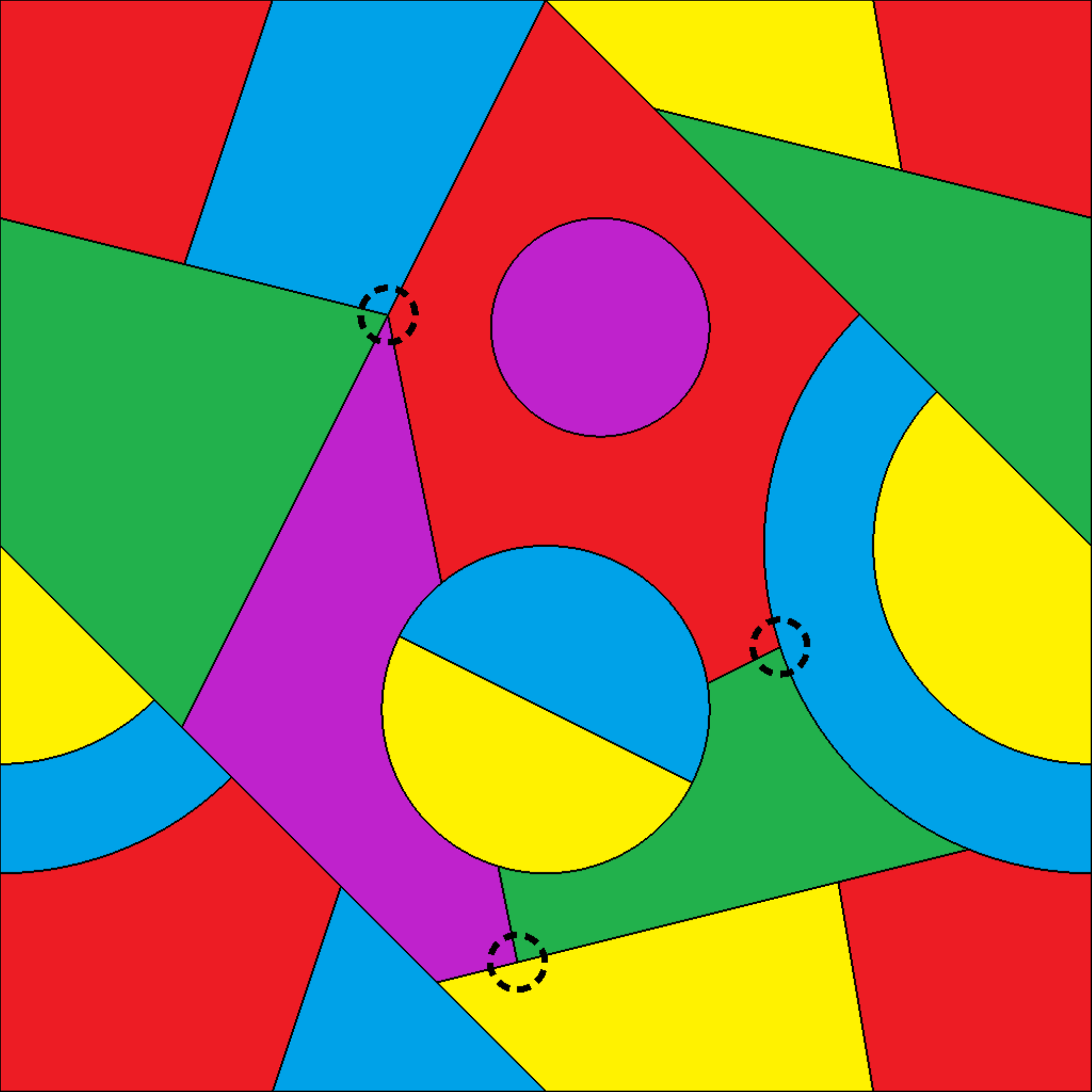}
			\caption{t=0}	
		\end{subfigure}
		\begin{subfigure}{0.2\linewidth}
			\includegraphics[width=\textwidth]{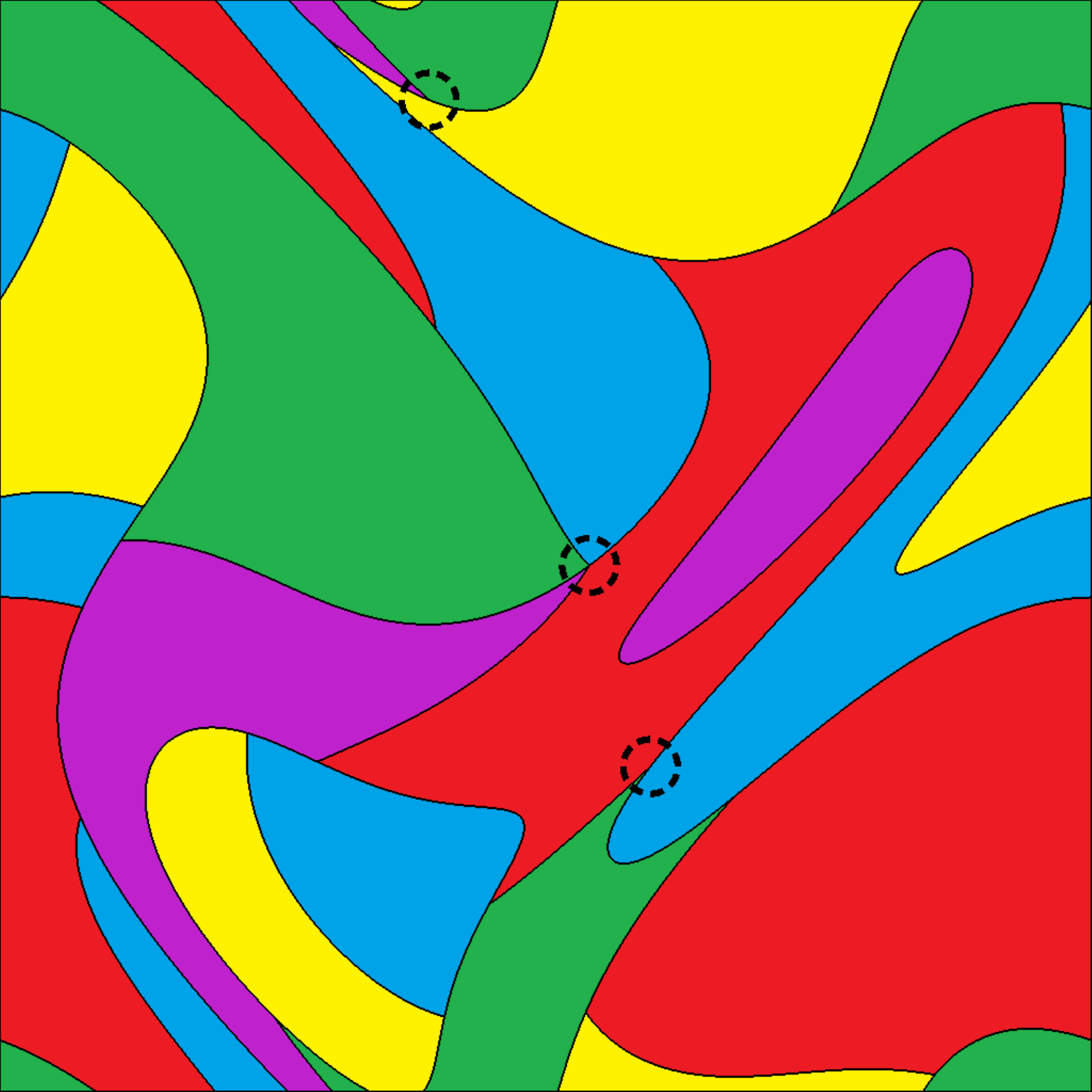}
			\caption{t=0.25}	
		\end{subfigure}
		\begin{subfigure}{0.2\linewidth}
			\includegraphics[width=\textwidth]{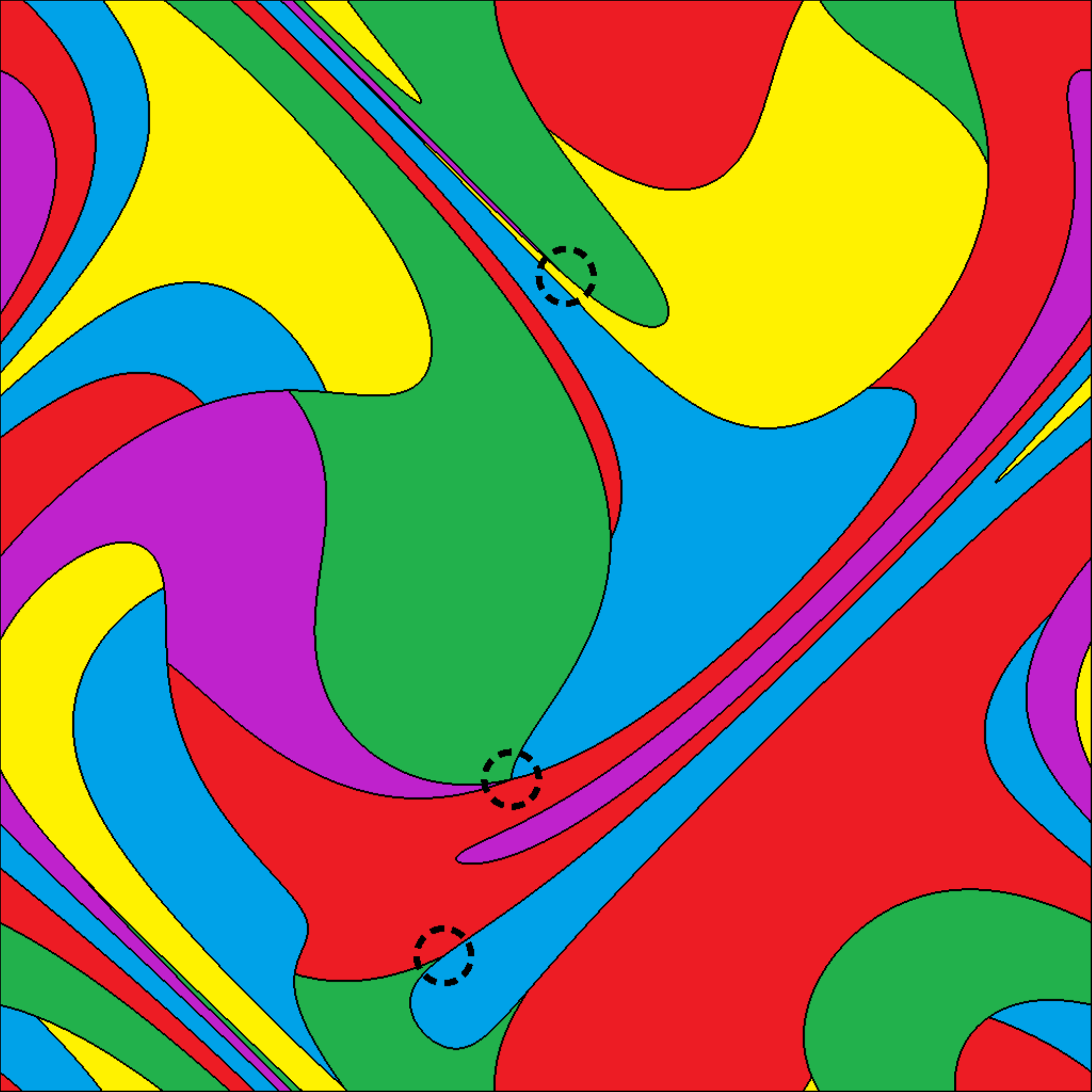}
			\caption{t=0.5}	
		\end{subfigure}\\
		\begin{subfigure}{0.2\linewidth}
			\includegraphics[width=\textwidth]{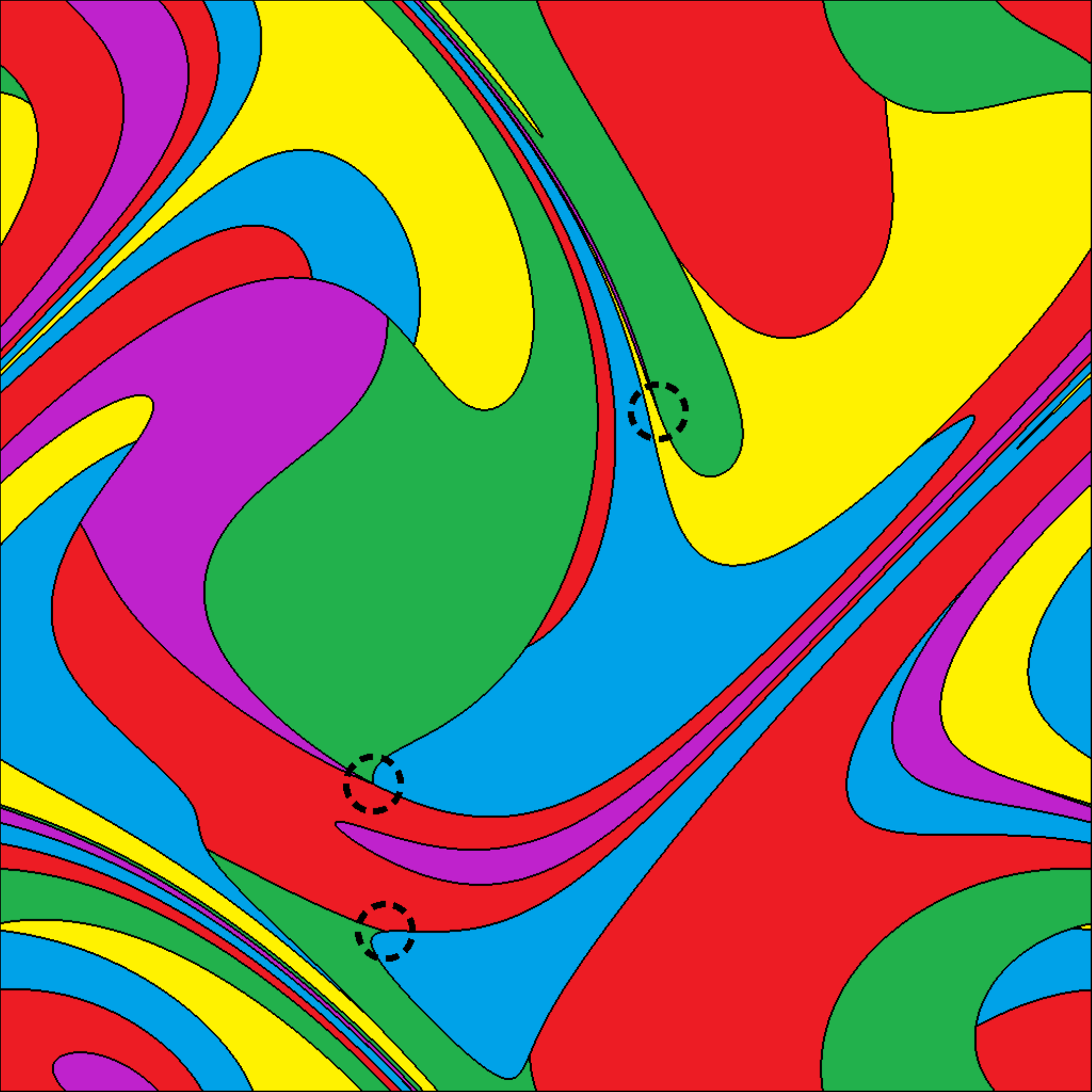}
			\caption{t=.75}	
		\end{subfigure}
		\begin{subfigure}{0.2\linewidth}
			\includegraphics[width=\textwidth]{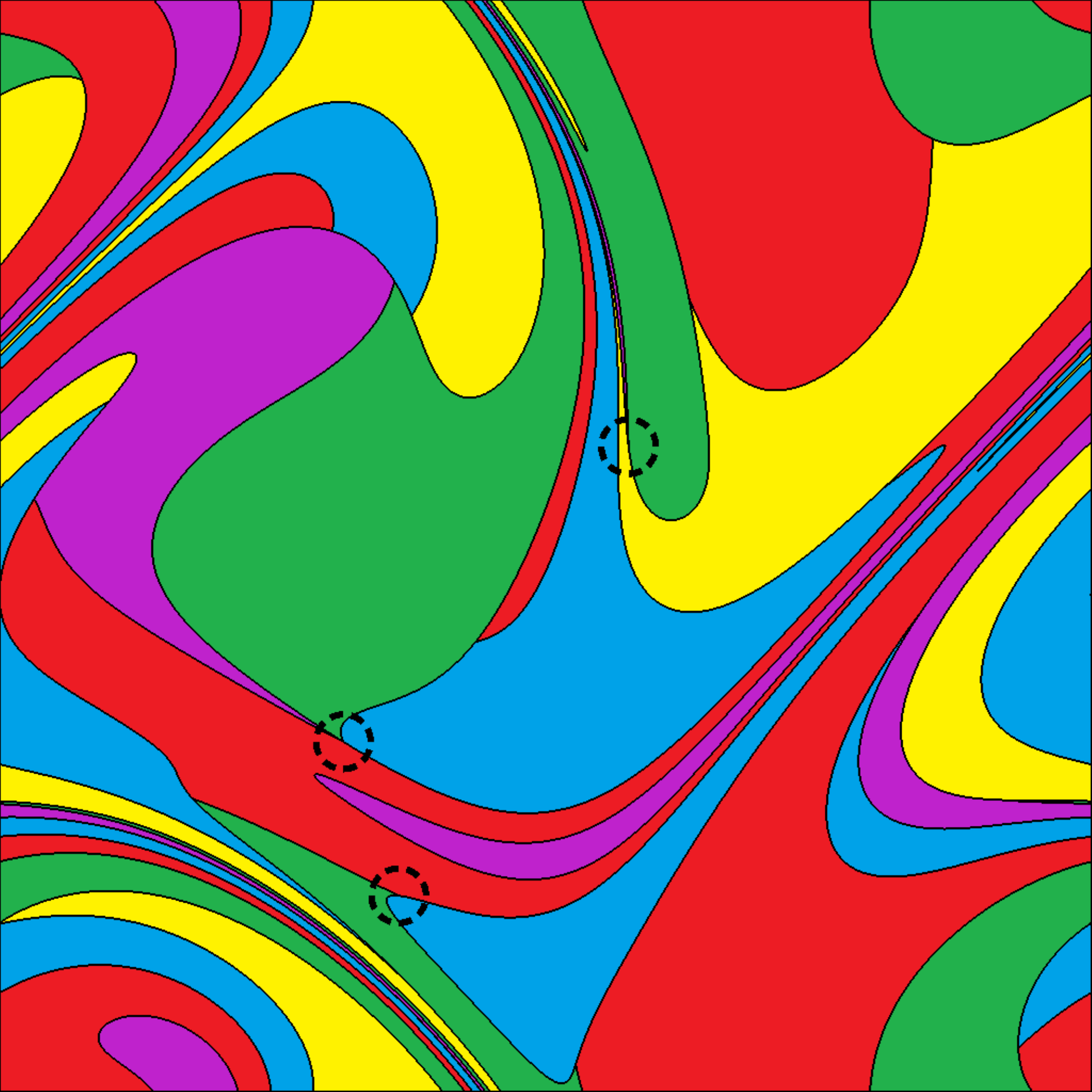}
			\caption{t=1}	
		\end{subfigure}
		\begin{subfigure}{0.2\linewidth}
			\includegraphics[width=\textwidth]{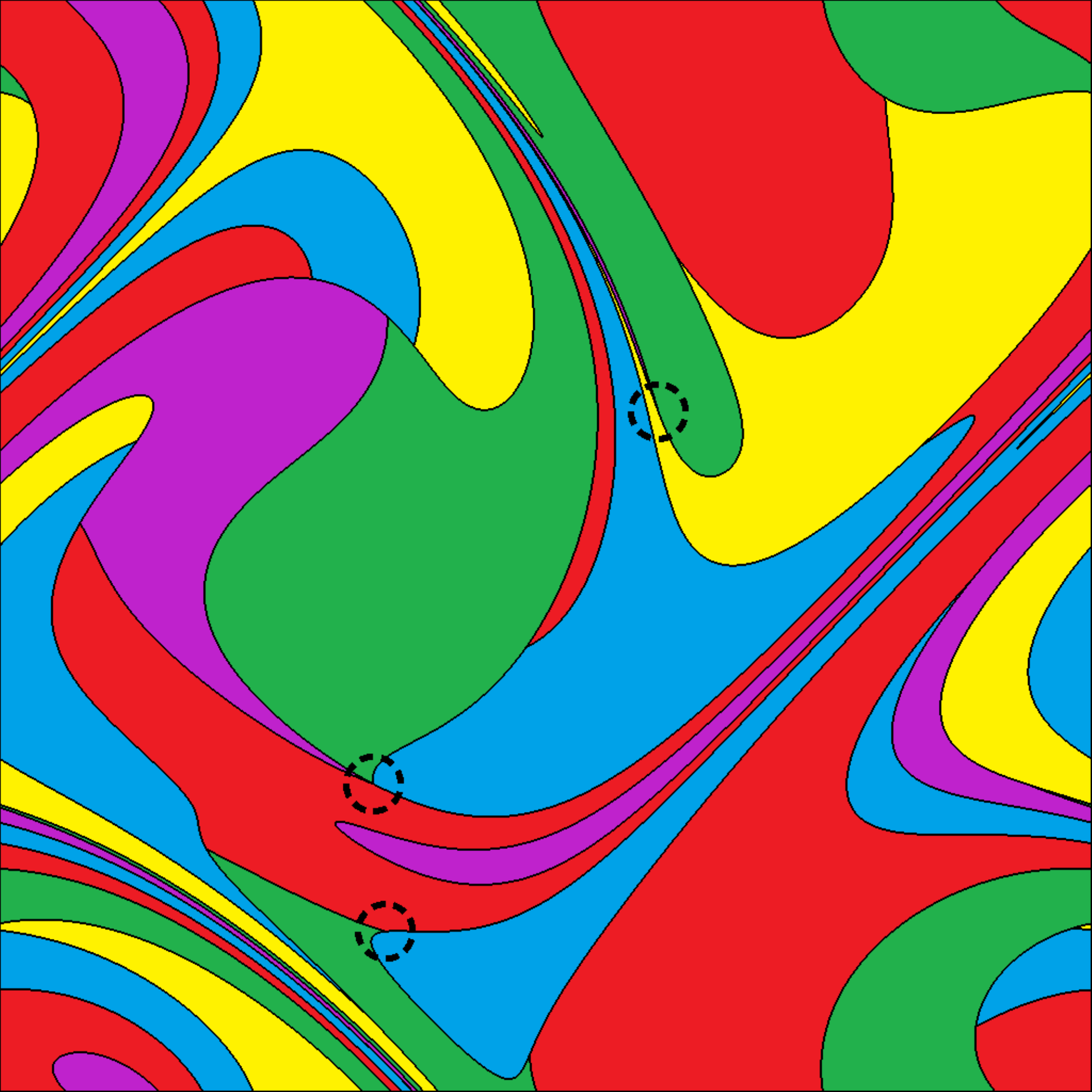}
			\caption{t=1.25}	
		\end{subfigure}\\
		\begin{subfigure}{0.2\linewidth}
			\includegraphics[width=\textwidth]{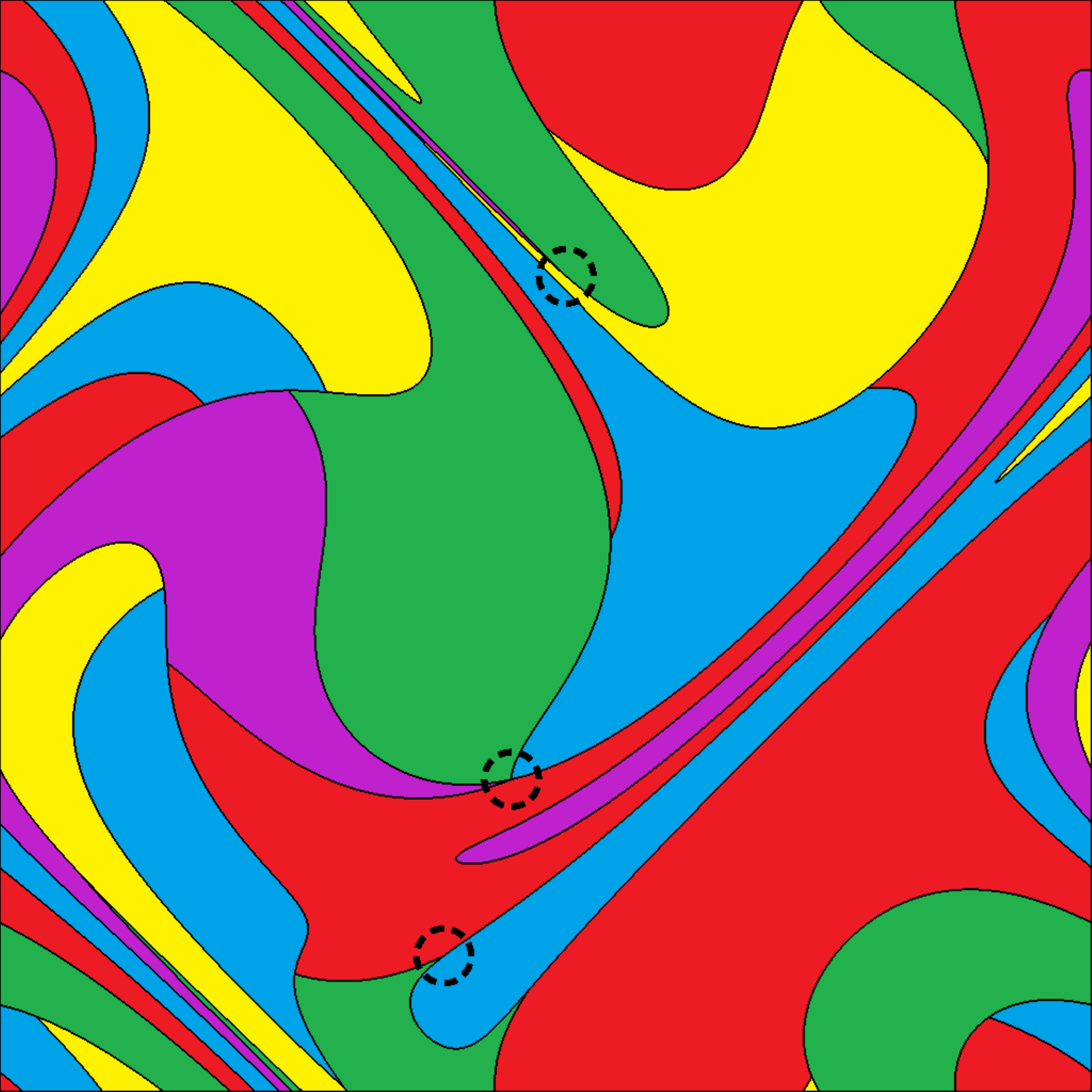}
			\caption{t=1.5}	
		\end{subfigure}
		\begin{subfigure}{0.2\linewidth}
			\includegraphics[width=\textwidth]{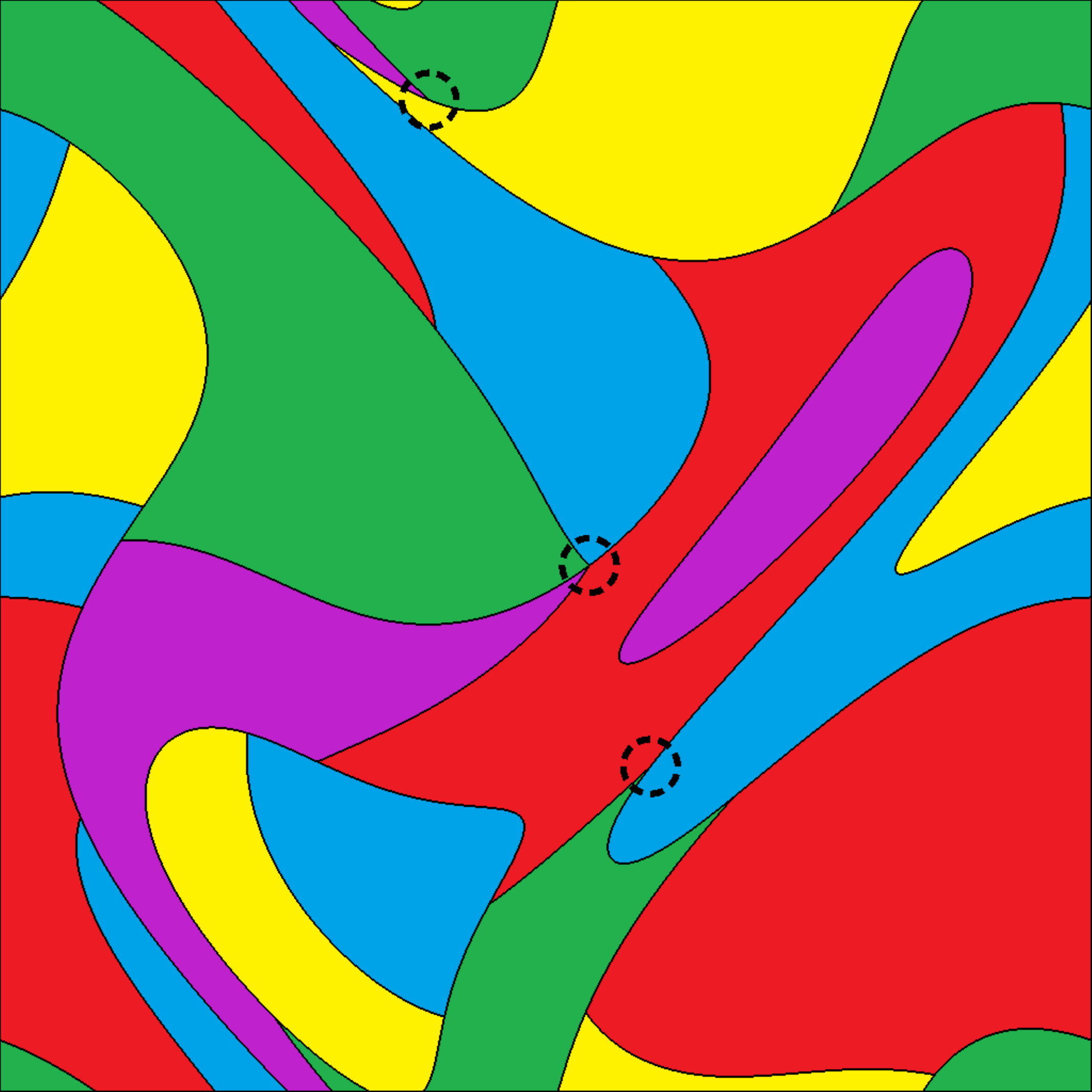}
			\caption{t=1.75}	
		\end{subfigure}
		\begin{subfigure}{0.2\linewidth}
			\includegraphics[width=\textwidth]{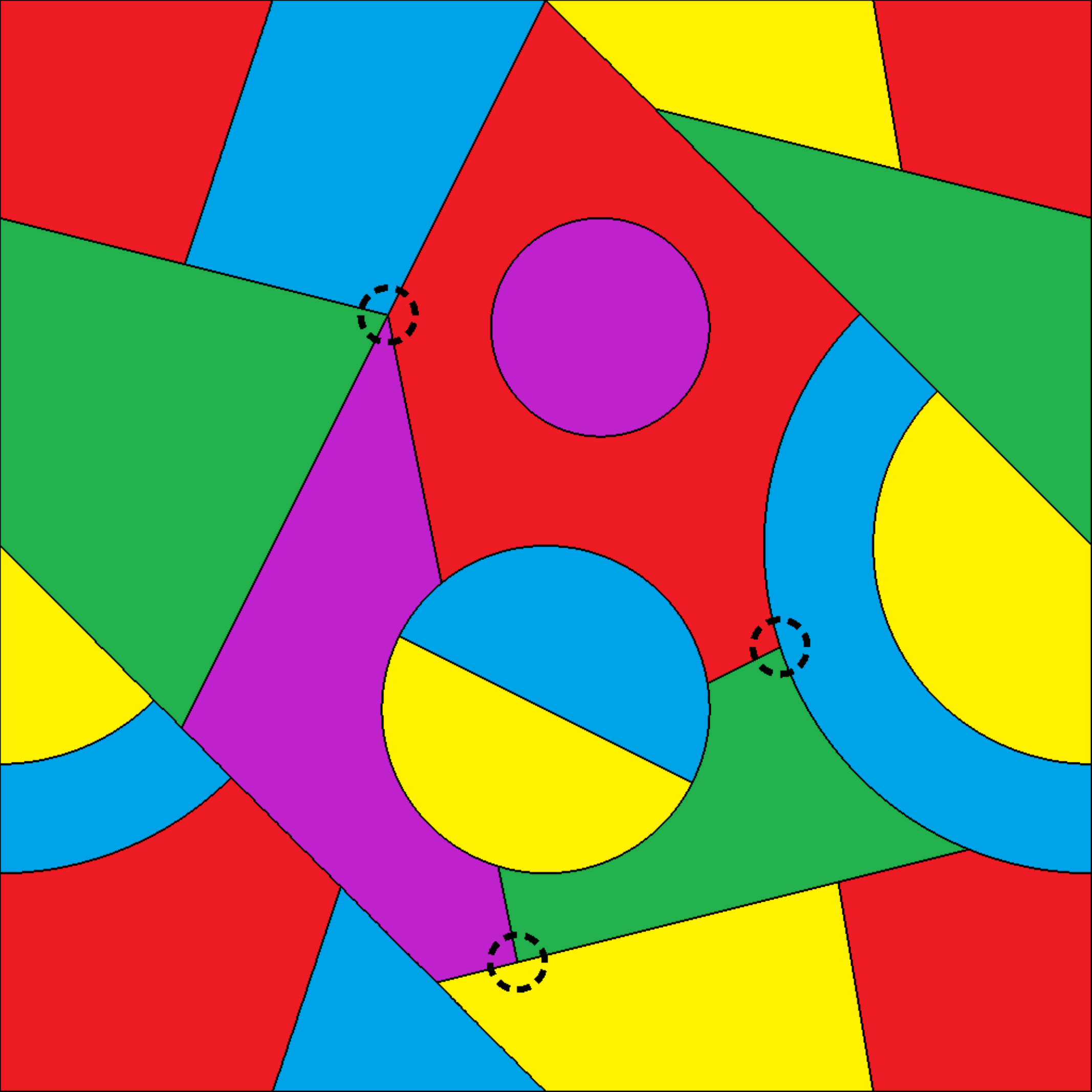}
			\caption{t=2}	
		\end{subfigure}
	\end{center}	
	\caption{2D deformation of a periodic mosaic pattern. $N_c = 32$ and $N_f = 512$. The three circled regions are tracked as passive Lagrangian particles.}
	\label{fig:mosaic}
\end{figure}

Resolution of thin filaments and accurate tracking of junction points are both important elements of a multiphase flow simulation. From figure \ref{fig:mosaic}, we see that the CM method manages to track very thin filaments. If each level-set was advected separately using traditional methods, the errors might not be consistent across individual level-sets. This could result in important qualitative errors on the presence and location of thin filaments, persisting throughout the entire simulation. However, since all level-sets in the CM method are transported using the same deformation map, relationships between different regions are naturally preserved, therefore resolving and maintaining the thin filaments. This also leads to an accurate and consistent tracking of triple and quadruple junctions. As seen in figure \ref{fig:mosaic}, even though we can see some error in the position of the quadruple junction when deformations are large, this error is not amplified and does not persist throughout the simulation. We see indeed that the junctions return to the correct positions at time $t=2$. This also indicates that the time symmetry present in the velocity is better preserved in the numerical solution: the characteristic map deforms until $t=1$, then returns to the identity map at $t=2$. This would be difficult to achieve with traditional methods, working directly with the level-set functions as the sharp and convoluted shapes of the regions will incur high numerical dissipation, thereby moving or destroying junction points. 

\subsection{Computational Efficiency} \label{sec:CPU}
Several elements factor into the efficiency of an algorithm, notably the number of computations per step and memory usage necessary to reach a certain accuracy, and the parallelizability of the method. In the following, we give some estimates on the error and running time of the method. For simplicity, we will restrict ourselves to fixed fine grid $N_f$ and try to examine the effect on efficiency from varying $N_c$, $N_f$ and $\mathcal{E}_1$.

For GALS, the error is know to be $\bigO( \incr{t} N_g^{-2} + N_g^{-4}) $. For the CM method, we proceed as follows.

First, we look at the error from each submap computation, $\vec{\chi}_{[\tau_{i-1}, \tau_i]}$. Let $L_i = \tau_i - \tau_{i-1}$. Each time step, the error from time integration and evaluation at foot point is of order $\bigO(\incr{t}^2N_c^{-2} + \incr{t}^4)$. Taking $\incr{t} \leq N_c^{-1}$, we simplify this error to roughly $(\mu_i + \sigma)\incr{t}^2N_c^{-2}$, where $\mu_i$ is a constant of the order of $D^4_x \vec{X}_{[\tau_{i-1}, \tau_i]}$ coming from interpolation errors and $\sigma \approx D^4_t \vec{v}$ from integration errors. The total submap grid errors is then
\begin{gather} \label{eq:submapGridError}
\text{submap grid value error} \approx \mu_i L_i \incr{t}N_c^{-2} + \sigma  L_i \incr{t}^3
\end{gather}

There is also a submap representation error when evaluating the Hermite cubic $\vec{\chi}_{[\tau_{i-1}, \tau_i]}$ at the remapping step. This error exists even if grid values are exact and is due to the limited resolution available for a Hermite cubic on a coarse grid:
\begin{gather} \label{eq:submapRepError}
\text{submap representation error} \approx \mu_i N_c^{-4}
\end{gather}

Finally, at each remapping step, we incur an error from evaluating the global map $\vec{\chi}_{[0, \tau_{i-1}]}$:
\begin{gather} \label{eq:globalmapComposeError}
\text{global map evaluation error} \approx \nu_i N_f^{-4}
\end{gather}
where $\nu_i$ is on the order of $D_x^4 \vec{X}_{[0, \tau_{i-1}]} $ (note that $\nu_1 = 0$).

We let $\mu$ and $\nu$ be the average of $\mu_i$ and $\nu_i$ over all remapping steps. Assuming there are $m$ remapping steps between $t=0$ and $t=T$, the global error is roughly
\begin{gather} \label{eq:globalError}
\text{global error} = \mu T \incr{t}N_c^{-2} + \sigma  T \incr{t}^3 + m \mu N_c^{-4} + m \nu N_f^{-4}
\end{gather}

Both $\mu$ and $m$ are influenced by the choice of $\mathcal{E}_1$. If $\mathcal{E}_1$ is small, remapping happens more often, making $L_i$ shorter. Since $m$ correlates inversely with $L_i$, it would increase. On the other hand, $\mu$ can be assumed to correlate directly with $L_i$, and hence would decrease. We make the following assumptions and estimates. 

We suppose that $\mu_i \approx B L_i$ for some constant $B$, that is, for short times, higher frequencies not resolved by the coarse grid accumulate linearly in time. It is also implied from this that $\nu \approx BT$, although this may be an overestimate. For small times, we can assume that sharp features in the deformation grow linearly with time. However, in global time, this might not always be true as seen in the swirl test, where the periodic velocity undoes the sharp features it created in the first half of the period. We will however assume that the estimate is generally of correct order of magnitude. From these assumptions, we reach the following estimate for the global error:
\begin{gather} \label{eq:globalErrorFull}
E \approx \nu \left( L\incr{t} N_c^{-2} + \sigma T \incr{t}^3 + N_c^{-4} + \frac{T}{L} N_f^{-4}  \right)
\end{gather}
where $L$ is a function which scales directly with $\mathcal{E}_1$ and inversely with $\incr{t}^3 + N_c^{-4}$. Details of this calculation can be found in Appendix \ref{sec:errorEstimates}.

Closer examination of this formula reveals three main regimes of the magnitude of the error depending on the choice of $\mathcal{E}_1$. When $\mathcal{E}_1$ is taken on the order of the local truncation error for submap evolution, i.e. $\mathcal{E}_1 = \bigO ( \incr{t}^4 +\incr{t} N_c^{-4} )$, the remapping step happens every or almost every submap update step. The resulting global error is $\bigO(\incr{t}^2 N_c^{-2} + \incr{t}^3 + N_c^{-4} + \incr{t}^{-1} N_f^{-4})$. The other extreme is when $\mathcal{E}_1$ is of the order of global trunction error, i.e. $\mathcal{E}_1 = \bigO ( \incr{t}^3 +\incr{t} N_c^{-2} + N_c^{-4} )$. In this case, barely any remapping occurs. The resulting error is $\bigO ( \incr{t}^3 +\incr{t} N_c^{-2} + N_c^{-4} + N_f^{-4} )$, with the $N_f$ term appearing if at least one remapping step occurred. Lastly, in the intermediate regime, with $\mathcal{E}_1$ roughly $\bigO ( \incr{t}^2 N_c^{-2} )$, we should observe errors of order $\bigO( \incr{t}^3  + N_c^{-4})$.

In order to compare the efficiency of GALS and CM, we need the following estimates on computational time.

For the GALS method, the cost consists of tracing back the footpoints, and then evaluating the interpolant at those locations.
\begin{gather}
\text{cost GALS } =  \underbrace{C_1 \frac{T}{\incr{t}} N_g^d}_\textrm{footpoints} + \underbrace{C_2 \frac{T}{\incr{t}} N_g^d}_\textrm{interpolation} \label{eq:CPU_GALS}
\end{gather}
for some constants $C_1$ and $C_2$, and where $d$ denotes the number of dimensions and $N_g$ is the number of grid cells for each spatial dimensions.

For the CM method, we also need to trace back the footpoints and interpolate the function, but that is done on the coarse $N_c^d$ grid, once for each of the $d$ coordinate functions. Additionally, we need to advect the test particles and do remapping steps when required. The computational cost is
\begin{gather}
\text{cost CM } = \underbrace{d C_1 \frac{T}{\incr{t}} N_c^d}_\textrm{footpoints} + \underbrace{d C_2 \frac{T}{\incr{t}} N_c^d}_\textrm{interpolation} + \underbrace{C_3 \rho \frac{T}{\incr{t}}  N_c^d}_\textrm{particles} + \underbrace{\frac{T}{L} d C_4 N_f^d}_\textrm{remapping} \label{eq:CPU_CM}
\end{gather}
where $\rho$ is the number of particles per cell. Recall that $L$ scales directly with respect to $\mathcal{E}_1$, so picking small $\mathcal{E}_1$ increases the frequency of remapping steps. This is the most costly operation of the method and hence strongly affects the efficiency. Roughly, the cost of CM is 
\begin{gather}
\text{cost CM } = \bigO \left(  \frac{N_c^d}{ \incr{t}} +   \frac{ \incr{t}^3 +  N_c^{-4}}{ \mathcal{E}_1 } N_f^d \right)
\end{gather}

\begin{figure}[h]
	\begin{center}
		\begin{subfigure}{0.45\linewidth}
		\includegraphics[width=\textwidth]{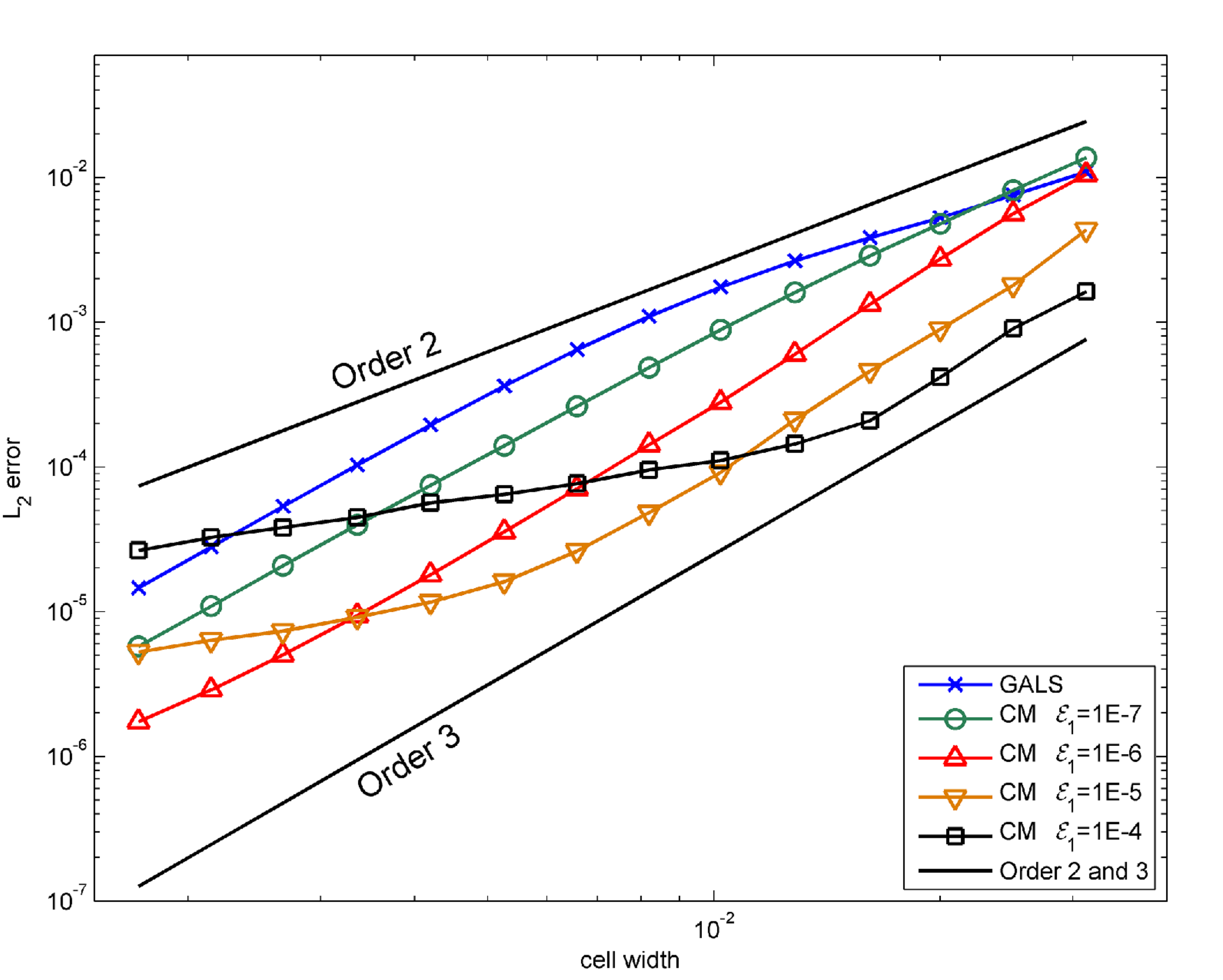} \caption{Error vs. cell width} \label{fig:CPU_error_cell}
		\end{subfigure}
		\begin{subfigure}{0.45\linewidth}
		\includegraphics[width=\textwidth]{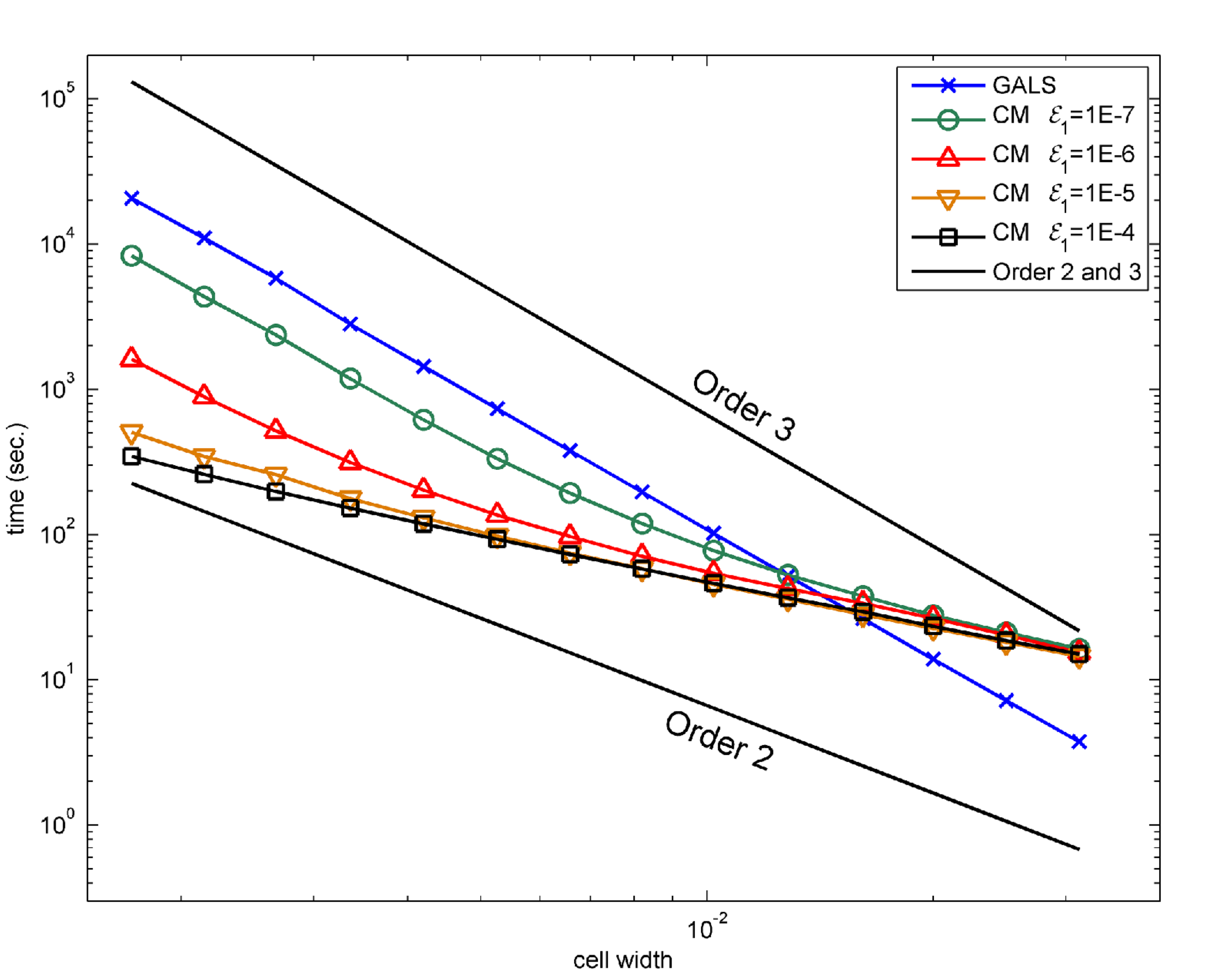} \caption{Time vs. cell width} \label{fig:CPU_time_cell}
		\end{subfigure}
		\begin{subfigure}{0.48\linewidth}
		\includegraphics[width=\textwidth]{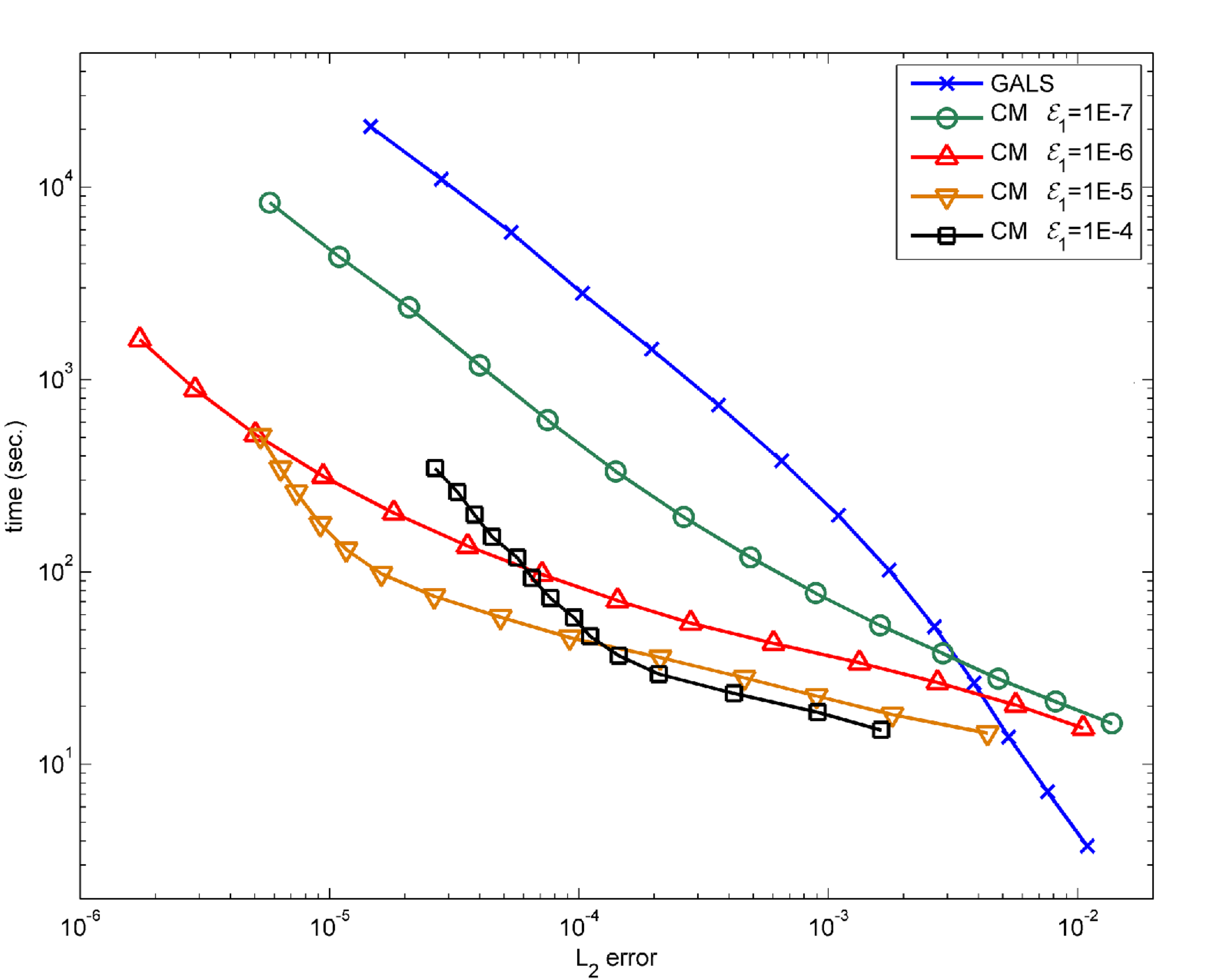} \caption{Time vs. error.} \label{fig:CPU_time_error}
		\end{subfigure}
	\end{center}	
	\caption{Error and computational time} \label{fig:all_cell}
\end{figure}

In order to illustrate the behavior of error, computational time and efficiency of the algorithms, we ran several numerical experiments using the swirl test presented in section \ref{sec:swirl} with $A=8$ and a final time $T=16$. We do not use a dynamic fine grid and set $N_f = N_g$, $\incr{t} = N_f^{-1}$, and fix $N_c = 32$ for all computations. For values of $\mathcal{E}_1$ ranging from $10^{-4}$ to $10^{-7}$, we test the error and running times for different values of $N_f$. The results are show in figure \ref{fig:all_cell}.

In figure \ref{fig:CPU_error_cell}, we observe that when $\mathcal{E}_1$ is much smaller than the global error, we see third order convergence with respect to $\incr{t}$ as described earlier ($N_c^{-2}$ was small enough for $N_c^{-2}\incr{t}^2$ terms not to appear). However, when $\mathcal{E}_1$ is close to the global error, almost no remapping step occurs and the error scales linearly in $\incr{t}$.

In figure \ref{fig:CPU_time_cell}, we can observe the effect of $\mathcal{E}_1$ on running time. When $\mathcal{E}_1$ is relatively large, the running time is dominated by submap time stepping: computational time scales linearly with $\incr{t}^{-1}$. However, when $\mathcal{E}_1$ is picked very small, the cost of remapping becomes more obvious. Indeed, we see that when picking $\mathcal{E}_1 = 10^{-7}$, for large $N_f$, running time scales like $\incr{t}^{-1} N_f^2$.

%\begin{figure}[h]
%	\begin{center}
%		\includegraphics[width=0.45\linewidth]{images/time_error.pdf}
%	\end{center}	
%	\caption{Time vs. error.}
%	\label{fig:CPU_time_error}
%\end{figure}
These results suggest that when a prescribed resolution on the velocity is chosen by fixing $N_c$, and a target error $\mathcal{E}_1$ for the advection by the approximated velocity is given, one can find an optimal fine grid size $N_f$ minimizing both error and computational time. 
%We can see from figure \ref{fig:CPU_time_error} that it is most efficient to pick $\mathcal{E}_1$ close to the desired global error. This in turn gives us the appropriate $N_f$ and $\incr{t}$ which we can estimate from \eqref{eq:globalError}.

Lastly, the CM method is easily parallelizable. Multiple transported interfaces can be advected by evolving a single characteristic map. Moreover, since the evolution of the characteristic map is done using the local GALS scheme, grid value updates are fully explicit and use compact localized stencils. The same applies to the remapping step due to it being essentially composed of two Hermite interpolant evaluations. It is not the aim of this paper to investigate the parallel implementation of the CM method, however work in this direction has the potential for creating a fast and accurate parallelized solver for the linear advection equation.
\section{Conclusion}
In this paper we have presented a new numerical approach for the linear advection of arbitrary sets. The new method relies on computing the solution as the pullback of the initial conditions by the deformation map. This method relies extensively on the GALS framework. The interpolation structure of GALS together with the deformation map approach of this paper allows for the representation of solutions to arbitrarily fine subgrid resolution. One key observation is that the deformation map can be decomposed in time and thus can be periodically reinitialized to the identity map whenever an error tolerance is reached. This led naturally to the remapping idea presented in section 2. Additionally, the observation that the submaps can be resolved on a coarser grid than the global map was used to devise the two-grid strategy presented in section 3. In section 4 we presented several benchmark tests and demonstrated the accuracy and efficiency of the method. Specifically, we showed that the proposed method is able to handle the advection of closed and open curves, fractal sets, and also more complicated sets containing multiple intersecting domains. The latter example may be of particular interest for the simulation of multiphase fluid flows. We also showed that depending on the error parameter in \eqref{eq:tolerance_1_1_1}, one may achieve $\bigO(N)$ time computation instead of the typical $\bigO(N^3)$ in 2D.

The current approach opens a wealth of possibilities for applications. However, since only the linear advection equation is considered, no topological changes are allowed. It is clear that the formalism presented in section 2 and used throughout the paper relies on existence of a diffeomorphism for all time and is thus incompatible with topological changes. Nonetheless, there may be extensions and modifications of the method that enable dealing with changes of topology, specifically when considering non-linear coupling of the velocity to the solution of the advected set. These questions are the focus of our current research in the subject as we believe the proposed method provides an appropriate framework upon which more challenging problems may be solved.

%\vspace*{1cm}
\newpage

\bibliographystyle{siamplain}
\bibliography{CMmethod}

\newpage
%\vspace*{1cm}
\appendix
\section{Detailed error estimates} \label{sec:errorEstimates}

The sup-norm error for each submap is approximately
\begin{gather} \label{eq:localMapError}
\mu_i L_i \incr{t}N_c^{-2} + \sigma  L_i \incr{t}^3 + \mu_i N_c^{-4}
\end{gather}

We assume a fixed $\sigma$ coming from the time step error of RK integration of the velocity. We suppose that $\mu_i \approx B L_i$ for some constant $B$, meaning that the sharp features not resolved by the coarse grid grows linearly with $t$ for small time intervals.

$L_i$ is the time it takes for the submap error to reach the threshold $\mathcal{E}_1$, we equate $\mathcal{E}_1$ with \eqref{eq:localMapError} to estimate $L_i$:
\begin{subequations} \label{eqGroup:LEstimate}
	\begin{align}
	& \mathcal{E}_1 \approx  B L_i^2 \incr{t}N_c^{-2} +  L_i (B N_c^{-4} + \sigma  \incr{t}^3) \\
	\Rightarrow & L_i \approx \left(  \left( \frac{B N_c^{-4} + \sigma \incr{t}^3 }{B \incr{t} N_c^{-2}} \right)^2 + \frac{\mathcal{E}_1}{B \incr{t} N_c^{-2}}  \right)^{1/2} -   \frac{B N_c^{-4} + \sigma \incr{t}^3 }{B \incr{t} N_c^{-2}}
	\end{align}
\end{subequations}

We drop the subscript $i$ to consider averaged values $L$ and $\mu$. We have $\mu = BL$ and let $m \approx T/L$ be the total number of remapping steps. We also assume $\nu \approx BT$, although this may be an overestimate.

From combining \eqref{eq:submapGridError}, \eqref{eq:submapRepError} and \eqref{eq:globalmapComposeError}, the global error then reads:
\begin{gather} \label{eq:globalErrorFullAppedix}
E \approx \nu \left( L\incr{t} N_c^{-2} + \sigma T \incr{t}^3 + N_c^{-4} + \frac{T}{L} N_f^{-4}  \right)
\end{gather}

We can more closely examine the effect of the magnitude of $\mathcal{E}_1$. Assuming it is small enough to justify a first order expansion of the square root in \eqref{eqGroup:LEstimate}, that is $L \approx \frac{\mathcal{E}_1}{\incr{t}^3 + N_c^{-4}} $, we can sample a few $(L, \mathcal{E}_1)$ pairs to get an idea of the behaviour of global error. Note that by construction, $L$ must be between $\incr{t}$ and $T$. This gives us the two extremes pairs $( \incr{t}, \incr{t}^4 + \incr{t} N_c^{-4} )$ and $(T, T \incr{t}^3 + T N_c^{-4})$ which correspond to taking $\mathcal{E}_1$ on the order of local and global truncation errors respectively. These give global errors $\bigO(\incr{t}^2 N_c^{-2} + \incr{t}^3 + N_c^{-4} + \incr{t}^{-1} N_f^{-4})$ and $\bigO ( \incr{t}^3 +\incr{t} N_c^{-2} + N_c^{-4} )$.

The global error is monotone with respect to $\incr{t}$ and $N_c^{-1}$. From the above two extreme estimates, the error should be between $\bigO (\incr{t}^2 N_c^{-2} + \incr{t}^3 )$ and $\bigO ( \incr{t} N_c^{-2} )$ (here we assume $\incr{t}$ is at most $N_c^{-1}$). The transition happens as $L$ ranges from $\bigO ( \incr{t})$ to $\bigO (1)$. For instance, for the $(L, \mathcal{E}_1) = ( \incr{t}^2 N_c^2 , \incr{t}^5 N_c^2 + \incr{t}^2 N_c^{-2} )$ pair, and for $( \incr{t}^{-1} N_c^{-2}, \incr{t}^2 N_c^2 + \incr{t}^{-1} N_c^{-6})$ we get $\bigO (\incr{t}^3 + N_c^{-4} )$ error.

%%\end{thebibliography}
\end{document}